\documentclass[onefignum,onetabnum]{siamart190516}

\usepackage{lipsum}
\usepackage{amsfonts}
\usepackage{graphicx,subfigure}
\usepackage{epstopdf}
\usepackage{algorithmic}
\usepackage{amsopn}
\usepackage{wrapfig}

\usepackage{amsmath, amsxtra, amscd,amssymb,mathtools,bm}
\usepackage{multirow,booktabs,makecell}

\usepackage{tikz}
\usetikzlibrary{shapes,arrows}
\tikzstyle{process} = [rectangle,minimum width=2cm,minimum height=1cm,text centered,text width =4cm,draw=black]

\ifpdf
  \DeclareGraphicsExtensions{.eps,.pdf,.png,.jpg}
\else
  \DeclareGraphicsExtensions{.eps}
\fi

\newsiamremark{example}{Example}
\newsiamremark{remark}{Remark}
\newsiamremark{assumption}{Assumption}
\newsiamremark{hypothesis}{Hypothesis}
\crefname{hypothesis}{Hypothesis}{Hypotheses}
\newsiamthm{claim}{Claim}

\headers{Fast Algorithm and Convergence Analysis for Dynamic MFP}{J. Yu, R. Lai, W. Li, and S. Osher}

\title{A Fast Proximal Gradient Method and Convergence Analysis for Dynamic Mean Field Planning\thanks{J. Yu and R. Lai's work are supported in part by an NSF Career Award DMS–1752934. W. Li and S. Osher's work are supported in party by AFOSR MURI FP 9550-18-1-502.}}

\author{Jiajia Yu\thanks{Department of Mathematics, Rensselaer Polytechnic Institute, United States 
  (\email{yuj12@rpi.edu}).}
\and Rongjie Lai\thanks{Corresponding author. Department of Mathematics, Rensselaer Polytechnic Institute, United States 
  (\email{lair@rpi.edu}).}
\and Wuchen Li\thanks{Department of Mathematics, University of South Carolina, United States
  (\email{wuchen@mailbox.sc.edu}).}
\and Stanley Osher\thanks{Department of Mathematics, University of California, Los Angeles, United States
  (\email{sjo@math.ucla.edu})}
}

\newcommand{\nt}{{n_0}}
\newcommand{\nx}{{n_1}}

\newcommand{\Id}{\text{Id}}

\newcommand{\Rho}{P}

\newcommand{\Rhobar}{\overline{\Rho}}

\newcommand{\Mbar}{\overline{M}}
\newcommand{\Phibar}{\overline{\Phi}}
\newcommand{\Jbar}{\overline{J}}

\newcommand{\half}{\frac{1}{2}}
\newcommand{\deri}{\text{d}}
\newcommand{\upk}{{(k)}}
\newcommand{\kp}{{(k+\frac{1}{2})}}
\newcommand{\kpp}{{(k+1)}}

\newcommand{\calA}{\mathcal{A}}
\newcommand{\calC}{\mathcal{C}}
\newcommand{\calD}{\mathcal{D}}
\newcommand{\calE}{\mathcal{Y}}
\newcommand{\calF}{\mathcal{F}}
\newcommand{\calG}{\mathcal{G}}
\newcommand{\calL}{\mathcal{L}}
\newcommand{\bbR}{\mathbb{R}}
\newcommand{\bbV}{\mathbb{V}}
\newcommand{\bbRplus}{\bbR^+}
\newcommand{\bbRbar}{\overline{\bbR}}
\newcommand{\bigO}{\mathcal{O}}
\newcommand{\proj}{\text{Proj}}
\newcommand{\prox}{\text{Prox}}

\newcommand{\bmx}{\bm{x}}
\newcommand{\bmy}{\bm{y}}
\newcommand{\bmz}{\bm{z}}
\newcommand{\bmv}{\bm{v}}
\newcommand{\bmb}{\bm{\beta}}
\newcommand{\bmm}{\bm{m}}
\newcommand{\bmM}{\bm{M}}
\newcommand{\bmn}{\bm{n}}
\newcommand{\bme}{\bm{e}}
\newcommand{\bmE}{\bm{E}}
\newcommand{\bmi}{\bm{i}}
\newcommand{\bmj}{\bm{j}}
\newcommand{\bmJ}{\bm{J}}
\newcommand{\bmO}{\bm{0}}

\DeclareMathOperator{\Res}{Res}
\DeclareMathOperator{\Pro}{Pro}
\DeclareMathOperator{\Lap}{Lap}
\DeclareMathOperator{\Divg}{Div}
\DeclareMathOperator{\Grad}{Grad}

\DeclareMathOperator{\divg}{div}

\DeclareMathOperator{\Solve}{\cref{alg: fista disct}}
\DeclareMathOperator*{\argmin}{argmin}

\begin{document}

\maketitle

\begin{abstract}
    In this paper, we propose an efficient and flexible algorithm to solve dynamic mean-field planning problems based on an accelerated proximal gradient method. Besides an easy-to-implement gradient descent step in this algorithm, a crucial projection step becomes solving an elliptic equation whose solution can be obtained by conventional methods efficiently. By induction on iterations used in the algorithm, we theoretically show that the proposed discrete solution converges to the underlying continuous solution as the grid size increases. Furthermore, we generalize our algorithm to mean-field game problems and accelerate it using multilevel and multigrid strategies. We conduct comprehensive numerical experiments to confirm the convergence analysis of the proposed algorithm, to show its efficiency and mass preservation property by comparing it with state-of-the-art methods, and to illustrates its flexibility for handling various mean-field variational problems. 
\end{abstract}



\section{Introduction}
\label{sec:intro}
Mean field planning (MFP) problems study how a large number of similar rational agents make strategic movements to minimize their cost in a process satisfying given initial and terminal density distributions \cite{achdou2012mean,gomes2014mean,porretta2014planning}. 
On the one hand, MFP can be viewed as a generalization of optimal transport (OT) \cite{benamou1999numerical,benamou2000computational,peyre2019computational,villani2008optimal} where no interaction cost is considered in the process. On the other hand, MFP is also a special case of mean field game (MFG) problems where the terminal density is often provided implicitly \cite{gueant2011mean,huang2007large,huang2006large,lasry2007mean}.  
MFP, MFG and OT have wide applications in economics \cite{achdou2014partial,achdou2017income, galichon2018optimal}, engineering \cite{de2019mean,gomes2018mean,yang2017mean}, quantum chemistry \cite{buttazzo2012optimal,cotar2013density}, image processing  \cite{haker2004optimal,papadakis2015optimal} as well as machine learning \cite{arjovsky2017wasserstein, ruthotto2020machine, weinan2019mean, yang2018mean}. 

More specifically, the dynamic MFP problem has the following optimization formulation:
\begin{equation}
    \min_{(\rho,\bmm)\in\calC(\rho_0,\rho_1)} \int_0^1\int_\Omega L(\rho(t,\bmx),\bmm(t,\bmx)) \deri \bmx \deri t + \int_0^1 \calF(\rho(t,\cdot))\deri t
\label{eqn: mfp opt prob}
\end{equation}
where $\rho(t,\bmx)$ is the densities of agents, $\bmm:=\rho\bmv$ with $\bmv$ representing the strategy(control) of this agent, and any pair of $(\rho,\bmm)\in \calC(\rho_0,\rho_1)$ satisfies mass conservation and zero boundary flux conditions with initial and terminal densities of $\rho$ being $\rho_0,\rho_1$ provided as:
\begin{equation}
\begin{aligned}
\calC(\rho_0,\rho_1) 
:&= \left\{(\rho,\bmm): \partial_t\rho + \divg_{\bmx} \bmm=0, \right.\\
 &\quad\quad\quad\quad\quad\left. \bmm\cdot\bmn = 0\text{ for }\bmx\in\partial\Omega, 
                                  \rho(0,\cdot) = \rho_0, \rho(1,\cdot) = \rho_1,\right\}.
\end{aligned}
\label{eqn: mfp opt constraint}
\end{equation}
In this variational problem, $L(\rho,\bmm)$ denotes the dynamic cost, $\calF$ models the interaction cost.
Specially, with $\calF= 0$ and a specific choice of $L$, variational problem \cref{eqn: mfp opt prob} reduces to the dynamic formulation of optimal transport (OT) proposed in \cite{benamou1999numerical,benamou2000computational}. By relaxing the given terminal density as an implicit condition regularized by a functional $\calG$, one can retrieve a class of MFG as the following formulation~\cite{benamou2017variational,gueant2011mean,lasry2007mean}:
\begin{equation}
\begin{aligned}
    \min_{(\rho,\bmm)\in\calC(\rho_0)} \int_0^1\int_\Omega L(\rho(t,\bmx),\bmm(t,\bmx)) \deri \bmx \deri t + \int_0^1 \calF(\rho(t,\cdot))\deri t + \calG(\rho(1,\cdot))\\
    \text{where }\calC(\rho_0):=\{(\rho,\bmm):\partial_t\rho+\divg_{\bmx} \bmm = 0, \bmm\cdot\bmn=0\text{ for } \bmx\in\partial\Omega, \rho(0,\cdot)=\rho_0 \}.
\end{aligned}
\label{eqn: mfg opt prob}
\end{equation}

Several numerical methods have been established to solve dynamic MFP, MFG and OT problems. One class of methods is based on solving partial differential equations (PDEs) corresponding to the KKT system of the variational problem~\cite{achdou2012mean,achdou2013mean,achdou2010mean,benamou2014numerical}, where conventional numerical methods in nonlinear PDEs can be applied. In principle, this class of methods can also be applied to handle general MFP and MFG problems that may not come from variational formulas. However, one obvious challenge of solving PDEs directly is their nonlinearity. This also limits solvers to handle broader choices of the dynamic cost $L$ and interaction cost $\calF$. 

Another class of methods focuses on the variational formulas of dynamic MFP, MFG and OT problems. By naturally combining with recent advances from optimization, existing methods include several first-order optimization algorithms to solve dynamic OT problems such as augmented Lagrangian \cite{benamou2015augmented,benamou2016augmented,papadakis2014optimal}, primal-dual \cite{lee2020generalized} and G-prox \cite{jacobs2019solving}, etc. These methods work on either the Lagrangian or the dual problem of the original optimization problem, particularly for dynamic OT where $\calF\equiv0$. These algorithms work very well since the involved sub-optimization problems have closed-form solutions. Besides, OT also has a static linear programming formulation and can be (approximately) solved by many algorithms, for example, Sinkhorn algorithm \cite{cuturi2013sinkhorn}. However, in the case of MFP or MFG problems, a broader choice of the interaction cost and the regularization term for the terminal density will make related sub-problems harder to solve than those in the OT problem. Meanwhile, it is unclear if a static linear programming formula still exists for MFP or MFG problems due to the contributions from $\calF$ and $\calG$. Besides, all of these algorithms may not preserve mass in the evolution very well as the mass constraint term is not explicitly forced. 

We would like to propose a method that can efficiently compute the mean-field type of problems with mass preservation property and great flexibility on a broad range of objective functions. Note that the mass conservation constraint in MFP is a convex set. A straightforward calculation shows that projection to this convex set can be obtained from solving a standard Poisson equation. This motivates us to propose another algorithm to solve MFP problems based on the proximal gradient descent method \cite{rockafellar1970convex,bauschke2011fixed}. This method is simply composed of a gradient descent step and a projection step. For MFP problems with a smooth objective function, the gradient descent step can be evaluated very efficiently. It also enjoys the flexibility to handle a broader range of $L$ and $\calF$. More importantly, the projection step leads to mass preservation in each iteration. This step can also be computed very efficiently by conventional fast algorithms for solving a Poisson equation. In this work, we use an accelerated version of the proximal gradient descent method, referred to as the fast iterative soft threshold algorithm (FISTA)~\cite{beck2009fast}, to solve the MFP problems. After that, we further generalize our algorithm to handle MFG problems. In addition, inspired by \cite{achdou2012iterative,li2020simple,liu2021multilevel}, we also apply multigrid and multilevel methods to speed up the proposed algorithm. Our numerical experiments illustrate the efficiency, mass preservation and flexibility of the proposed algorithm to different MFP problems as well as MFG problems. The vanilla version of our algorithm performs comparable with state-of-the-art methods, while the multigrid and multilevel accelerated versions are much more efficient than state-of-the-art methods.

\begin{figure}
\centering
\subfigure[Proof in \cite{achdou2012mean}]{
\includegraphics[width=6cm]{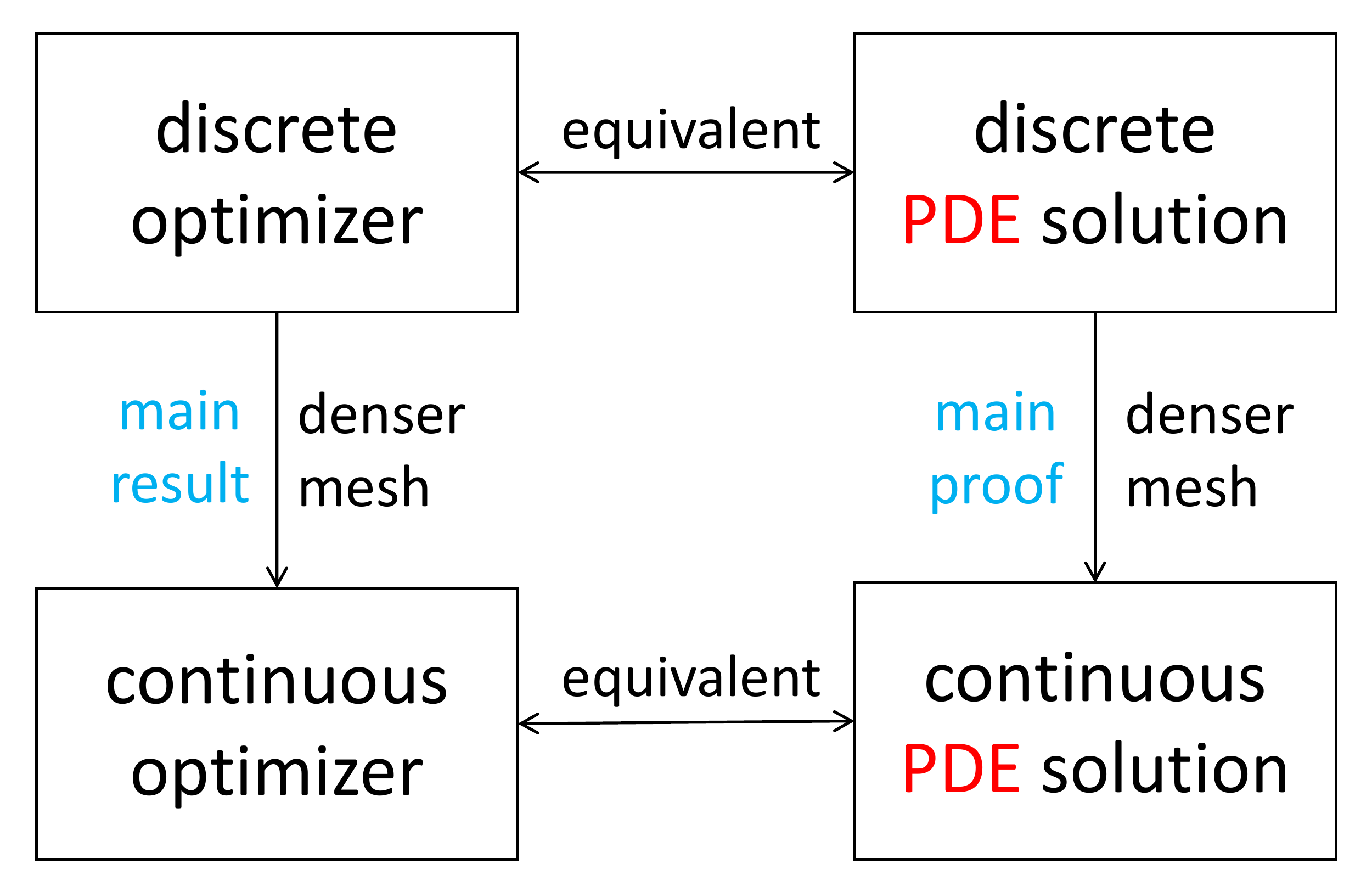}
}\hfill
\subfigure[Our proof]{
\includegraphics[width=6cm]{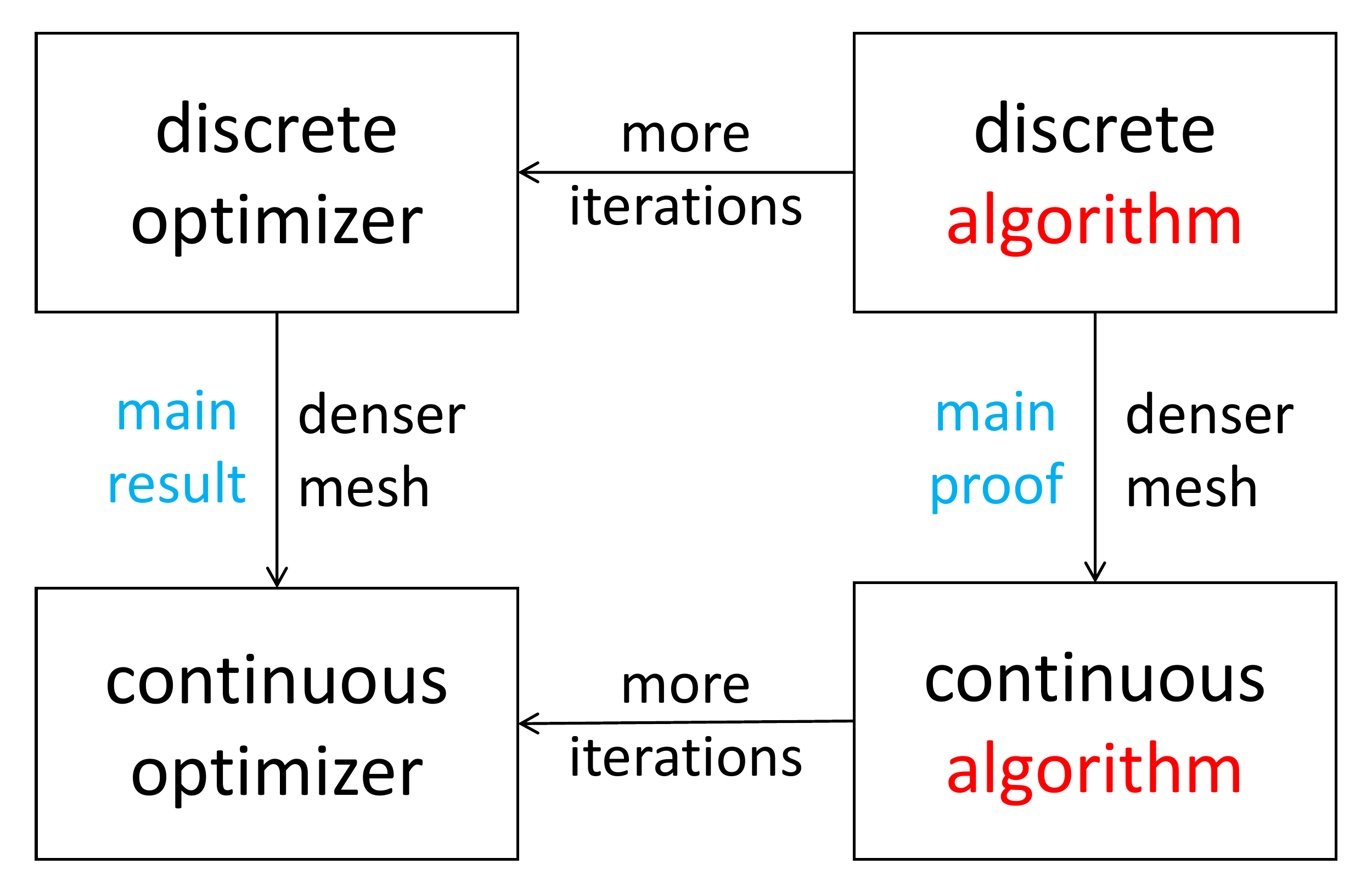}
}
\label{fig: proof sketch}
\caption{Sketch of two approaches of convergence proof.}
\end{figure}

Besides proposing a new algorithm for MFP problems, we also analyze errors between the discrete solution and the continuous solution. Since MFP is a functional optimization problem, all numerical methods on a given mesh grid only provide approximated solutions to the continuous problem.  It is important to understand how close the discrete numerical solution is to the continuous solution on a given mesh grid. Our analysis is from the algorithm perspective. We first derive an algorithm to optimize the variational problem and discretize each step of our algorithm. Our main effort is to prove that at each iteration, the discrete values are not far from the underlying continuous function values on grid points. Therefore we can show that the discrete algorithm converges to the continuous optimizer on grid points under certain conditions. Similar types of analysis may not be conveniently conducted in the existing methods including augmented Lagrangian, primal-dual and G-Prox since it could be difficult to have desired perturbation analysis of solving cubic equations involved in these three methods. 
To the best of our knowledge, this is for the first time to examine the discretization error based on variational MFP and its optimization algorithm. 
We remark that the convergence analysis has also been studied in \cite{achdou2012mean,achdou2013mean,achdou2020mean} from the PDE perspective, where the authors argue solution of discrete KKT converges to the continuous solution based on the equivalence of continuous systems and discrete systems. 
We indicate the major difference between our error analysis based on optimization perspective and error analysis based on PDEs perspective in \cref{fig: proof sketch}.

\textbf{Contributions:} We summarize our contributions as follows:
\begin{enumerate}
    \item We propose to use an accelerate proximal gradient method to solve the  MFP problem \cref{eqn: mfp opt prob}. 
    \item We analyze the error between the each iteration of discrete optimization  and its continuous counter part. We prove that the discrete solution converges to continuous optimizer on grid points as the mesh size converges. 
    \item We apply multilevel and multigrid strategies to to accelerate our algorithm.  We also generalize our algorithm to solve MFG problems.
    \item We conduct comprehensive numerical experiments to illustrate the efficiency and flexibility of our algorithms. 
    
\end{enumerate}

\textbf{Organization: } Our paper is organized as follows. In \cref{sec: review}, we give a brief review of the MFP problem and provide several MFP examples and a MFG example.
After that,  we describe our algorithm and the implementation details in \cref{sec: alg}. 
In \cref{sec: convergence}, we analyze the discretization error in our algorithm and prove the main theoretical result on the convergence of discrete solutions to the continuous solution as the grid size goes to zero. Furthermore, we generalize our algorithm to solve variational MFG problems and accelerate our algorithm by multilevel and multigrid methods in \cref{sec: acc and general}. 
Numerical experiments are provided in \cref{sec: num eg} to demonstrate the convergence order and to illustrate the efficiency and flexibility of our algorithm. At last we conclude this paper in \cref{sec: summary}.


\section{Review}
\label{sec: review}
In this section, we briefly review MFP problem and provide several examples which will be computed in the experiment section. 

Consider the model on time interval $[0,1]$ and space region $\Omega\in\bbR^D.$
Let $\rho:[0,1]\times\Omega\to\bbRplus:=\{a\in\bbR| a\geq0\}$ be the density of agents through $t\in[0,1]$, and $\bmm = (m_1,\cdots,m_D):[0,1]\times\Omega\to\bbR^D$ be the flux of the density which models strategies (control) of the agents. 
We are interested in $\rho$ with given initial and terminal density $\rho_0,\rho_1$ and
$(\rho,\bmm)$ satisfying zero boundary flux and mass conservation law, which gives the constraint set $\calC(\rho_0,\rho_1)$ defined in \cref{eqn: mfp opt constraint}.
We denote $L:\bbRplus\times\bbR^D\to\bbRbar:=\bbR\cup\{\infty\}$ as the dynamic cost function (e.g. \cref{eqn: dynamical term L} in this paper) and $\calF:\Omega^*\to\bbRbar$ as a functional modeling interaction cost.
The goal of MFP is to minimize the total cost among all feasible $(\rho,\bmm)\in\calC(\rho_0,\rho_1).$ Therefore the problem can be formulated as 
\begin{equation}
    \min_{(\rho,\bmm)\in\calC(\rho_0,\rho_1)} \int_0^1\int_\Omega L(\rho(t,\bmx),\bmm(t,\bmx)) \deri \bmx \deri t + \int_0^1 \calF(\rho(t,\cdot))\deri t.
\end{equation}

It is clear to see $\calC(\rho_0,\rho_1)$ is convex and compact. In addition, the mass conservation law $\partial_t\rho+\divg_{\bmx}\bmm=0$ and zero flux boundary condition $\bmm\cdot\bmn=0,\bmx\in\partial\Omega$ imply that $\calC(\rho_0,\rho_1)\neq\emptyset$ if and only if $\int_\Omega\rho_0=\int_\Omega\rho_1.$ 
Once $\calC(\rho_0,\rho_1)$ is non-empty, the existence and uniqueness of the optimizer depends on $L$ and $\calF$.

There are many different choices of $\calF$. In this paper, we consider
\begin{equation}
    \calF(\rho(t,\cdot)):=\lambda_E\int_\Omega  F_E(\rho(t,\bmx)) \deri\bmx + \lambda_Q\int_\Omega \rho(t,\bmx)Q(\bmx) \deri\bmx.
\end{equation}
where $\lambda_E,\lambda_Q \geq 0$ are two parameters,  $F_E:\bbRplus\to\bbR$ serves as a function to regularize $\rho$, and
$Q(\bmx):\Omega\to\bbRbar$ provides a moving preference for density $\rho$. Consider an illustrative example by choosing $\Omega_0\subset\Omega$ and  assuming $Q(\bmx) = \begin{cases}
0,\quad \bmx\in \Omega_0\\
+\infty,\quad \bmx \not\in \Omega_0
\end{cases}$,
then the mass has to move within $\Omega_0$ in order to keep the cost finite. In more general choice of $Q$, $\rho(t,\bmx)$ tends to be smaller at the place where $Q(\bmx)$ is larger and vice versa.

We then briefly discuss several concrete examples which will be considered in our numerical experiments. 

\begin{example}[Optimal transport~\cite{benamou1999numerical}]
In this paper, we consider a typical dynamic cost function $L$ by \begin{equation}
    L(\beta_0,\bmb) :=  \begin{cases} 
    \frac{\|\bmb\|^2}{2\beta_0} &\text{ if } \beta_0>0\\
    0              &\text{ if } \beta_0=0,\bmb=\bmO\\
    +\infty        &\text{ if } \beta_0=0, \bmb\neq\bmO.
    \end{cases},
\label{eqn: dynamical term L}
\end{equation}
If $\lambda_E=\lambda_Q=0$, the MFP becomes the dynamic formulation of optimal transport problem:
\begin{equation}
  \text{(OT) }  \min_{\rho,\bmm\in\calC(\rho_0,\rho_1)} \int_0^1\int_\Omega L(\rho(t,\bmx),\bmm(t,\bmx)) \deri \bmx \deri t. 
\label{eqn: mfp eg ot}
\end{equation}
Since $\bmm=\rho\bmv$, this definition of $L$ makes sure that $\bmm=\bmO$ wherever $\rho=0$. 
Because $\lambda_E=\lambda_Q=0$, OT can be viewed as a special case of MFP where masses move freely in $\Omega$ through $t\in[0,1]$.
\end{example}

\begin{example}[Crowd motion~\cite{ruthotto2020machine}]
Consider $F_E:\bbRplus\to\bbR, \rho\mapsto\begin{cases}
\rho\log(\rho),& \rho>0\\
0,& \rho=0
\end{cases}$,
and write $\Omega^+:=\Omega\cap\{\bmx\in\Omega:\rho(t,\bmx)>0\}$, 
we have the crowd motion model
\begin{equation}
  \min_{\rho,\bmm\in\calC(\rho_0,\rho_1)}
  \left\{\begin{aligned}
    &\int_0^1\int_\Omega L(\rho(t,\bmx),\bmm(t,\bmx)) \deri \bmx\deri t   \\
  + \lambda_E&\int_0^1\int_{\Omega^+} \rho(t,\bmx)\log(\rho(t,\bmx))\deri\bmx\deri t+ \lambda_Q\int_0^1\int_\Omega\rho(t,\bmx)Q(\bmx)\deri\bmx\deri t
  \end{aligned}\right\}.
\label{eqn: mfp eg model1}
\end{equation}
With $F_E$ decreasing on $[0,e^{-1}]$ and increasing on $[e^{-1},+\infty)$, $\rho(t,\bmx)$ tends to be close to $e^{-1}$ everywhere. So we expect to have the density $\rho(t,\bmx)$ to be not sparse and not very large everywhere. 
\end{example}

\begin{example}
If $\displaystyle F_E:\bbRplus\to\bbR,\rho\mapsto\begin{cases}
\frac{1}{|p|}\rho^p,\quad \rho>0\\
0,\quad \rho=0
\end{cases}$, where $p=2$ or $-1$, then we have the following two models.
\begin{equation}
  \min_{\rho,\bmm\in\calC(\rho_0,\rho_1)}
  \left\{\begin{aligned}
    &\int_0^1\int_\Omega L(\rho(t,\bmx),\bmm(t,\bmx)) \deri \bmx \deri t   \\
  + \lambda_E&\int_0^1\int_\Omega \frac{\rho^2(t,\bmx)}{2}\deri\bmx\deri t+ \lambda_Q\int_0^1\int_\Omega\rho(t,\bmx)Q(\bmx)\deri\bmx\deri t
  \end{aligned}\right\}
\label{eqn: mfp eg model2}
\end{equation}

\begin{equation}
  \min_{\rho,\bmm\in\calC(\rho_0,\rho_1)}
  \left\{\begin{aligned}
    &\int_0^1\int_\Omega L(\rho(t,\bmx),\bmm(t,\bmx))\deri \bmx \deri t  \\
    + \lambda_E&\int_0^1\int_{\Omega^+} \frac{1}{\rho(t,\bmx)}\deri\bmx\deri t+ \lambda_Q\int_0^1\int_\Omega\rho(t,\bmx)Q(\bmx)\deri\bmx\deri t
  \end{aligned}\right\}
\label{eqn: mfp eg model3}
\end{equation}
In \cref{eqn: mfp eg model2}, by Cauchy-Schwartz inequality, we have
\begin{equation}
    \left(\int_\Omega\rho(t,\bmx)\deri\bmx\right)^2\leq \int_\Omega \rho^2(t,\bmx)\deri\bmx \int_\Omega 1 \deri\bmx,
\end{equation}
therefore $\int_\Omega\rho^2(t,\bmx)\deri \bmx$ has a lower bound and achieves the lower bound when $\rho(t,\cdot)$ is a constant over $\Omega$. Therefore, model \cref{eqn: mfp eg model2} guides the solution density uniformly distributed over $\Omega$ .
In \cref{eqn: mfp eg model3}, since the total mass $\int_\Omega\rho(t,\bmx)\deri\bmx$ is fixed and $\displaystyle \frac{1}{\rho}$ is larger when $\rho$ is smaller, the value of regularization term $\displaystyle \lambda_E \int_{\Omega^+} \frac{1}{\rho(t,\bmx)}\deri\bmx$ is smaller if $\rho(t,\bmx)$ accumulates at several sites and vanishes at other regions. Therefore model \cref{eqn: mfp eg model3} pursues a sparse optimizer $\rho(t,\bmx)$.
\end{example}

\begin{example}[A MFG model~\cite{benamou2017variational, gueant2011mean,lasry2007mean}] We provide a MFG example to complete this section.
In the cases, the terminal density $\rho_1$ is not explicitly provided but it satisfies a given preference. This preference can be imposed by regularizing $\rho(1,\cdot)$ in the same spirit as $\int_\Omega\rho(t,\bmx)Q(\bmx)\deri\bmx$ and obtain the following MFG model,
\begin{equation}
  \min_{(\rho,\bmm)\in\calC(\rho_0)} 
  \left\{\begin{aligned}
    &\int_0^1\int_\Omega L(\rho(t,\bmx),\bmm(t,\bmx)) \deri \bmx \deri t   \\
    +\lambda_E&\int_0^1\int_\Omega \rho(t,\bmx)\log(\rho(t,\bmx))\deri\bmx\deri t+ \lambda_Q\int_0^1\int_\Omega\rho(t,\bmx)Q(\bmx)\deri\bmx\deri t \\
     + \lambda_G&\int_\Omega \rho(1,\bmx)G(\bmx)\deri\bmx.
\end{aligned}\right\}
\label{eqn: mfp eg mfg}
\end{equation}
Here $\lambda_G>0$ is a parameter,  $G:\Omega\to\bbRbar$ gives a preference of the distribution of $\rho(1,\bmx)$, and $ \calC(\rho_0):=\{(\rho,\bmm):\partial_t\rho+\divg_{\bmx} \bmm = 0, \bmm\cdot\bmn=0\text{ for } \bmx\in\partial\Omega, \rho(0,\cdot)=\rho_0 \}$.
\end{example}


\section{Algorithm}
\label{sec: alg}
In this section, we first briefly review FISTA algorithm proposed in \cite{beck2009fast}. Using a first-optimize-then-discretize approach, we describe the FISTA algorithm on variational problem \cref{eqn: mfp opt prob}. After that, we provide the details of our discretization and implementation for the MFP. In the end of this section, we discuss a different approach based on first-discretize-then-optimize strategy which turns out leading to same discrete algorithm.


To solve general nonsmooth convex model
\begin{equation*}
    \min_{\bmz} f(\bmz)+g(\bmz),
\end{equation*}
where $f$ is a smooth convex function and $g$ is convex but possibly nonsmooth, one can apply proximal gradient method~\cite{bauschke2011fixed,rockafellar1970convex}.
\begin{equation*}
\begin{aligned}
    &\bmz^\kpp:=\prox_{\eta^\upk g}\left(\bmz^\upk-\eta^\upk \nabla f\left(\bmz^\upk\right)\right).
\end{aligned}
\end{equation*}
Here $\eta^\upk>0$ is the stepsize and the proximal operator is defined as:
\begin{equation}
    \prox_{\eta g}(\bmz):=\argmin_{\bmy} \left\{ g(\bmy)+\frac{1}{2\eta}\|\bmy-\bmz\|_2^2 \right\}.
\end{equation}
In particular, for an indicator function
$\chi_\calC(\bmz)=\left\{\begin{array}{cc}
    0, & \bmz \in\calC\\
    +\infty, & \bmz\not\in\calC
    \end{array}\right.$ of a convex set $\calC$,  
its proximal operator is exactly the projection operator to $\calC$, i.e.
\begin{equation*}
    \prox_{\eta \chi_\calC}(\bmz)= \proj_\calC(\bmz)=\argmin_{\bmy\in\calC} \frac{1}{2}\|\bmy-\bmz\|_2^2 , \quad\forall \eta>0.
\end{equation*}
FISTA is essentially an accelerated proximal gradient algorithm~\cite{beck2009fast}. It introduces $\widehat{\bmz}^\upk$ as a linear combination of $\bmz^\upk$ and $\bmz^{(k-1)}$ in each iteration, and conducts proximal gradient on $\widehat{\bmz}^\upk$ to obtain $\bmz^\kpp$. The algorithm is summarized in \cref{eqn:fista no backtrack}, where the stepsizes $\eta^\upk$ can either be a constant or be obtained by a backtracking line search. 
\begin{equation}
\label{eqn:fista no backtrack}
    \left\{\begin{aligned}
        \bmz^\kpp &=\prox_{\eta^\upk g}\left(\widehat{\bmz}^\upk - \eta^\upk\nabla f(\widehat{\bmz}^\upk)\right);  \\
        \tau^\kpp & =\frac{1}{2}\left(1+\sqrt{1+4\left(\tau^\upk\right)^2}\right);\\
        \widehat{\bmz}^\kpp &= \bmz^\kpp + \frac{\tau^\upk-1}{\tau^\kpp} \left(\bmz^\kpp - \bmz^\upk\right).
    \end{aligned}\right.
\end{equation}

As proved in \cite{beck2009fast}, if $\bmz^*=\argmin_{\bmz} f(\bmz)+g(\bmz)$, and $\left\{ \bmz^\upk \right\}$ is generated by FISTA, then 
\begin{equation*}
    \left[ f\left(\bmz^\upk\right)+g\left(\bmz^\upk\right) \right] - \left[ f(\bmz^*)+g(\bmz^*) \right] = \bigO\left(\frac{1}{(k+1)^2}\right).
\end{equation*}
\subsection{FISTA for MFP}
\label{sec: fista for mfp}
To apply the above FISTA method to problem \cref{eqn: mfp opt prob}, let us write
\begin{equation}
    \min_{\rho,\bmm\in\calC(\rho_0,\rho_1)} \calE(\rho,\bmm):=\int_0^1\int_\Omega Y(\rho(t,\bmx), \bmm(t,\bmx), \bmx) \deri \bmx \deri t,
\label{eqn: mfp cts prob}
\end{equation}
where
\begin{equation}
Y(\beta_0,\bmb, \bmx) = L(\beta_0,\bmb) + \lambda_E F_E(\beta_0)+\lambda_Q \beta_0 Q(\bmx).
\end{equation}
For convenience, we write $\displaystyle F_E' = \frac{\mathrm{d}}{\mathrm{d} \beta_0}F_E, L_0 = \frac{\partial}{\partial \beta_0} Y$, $\displaystyle \nabla_{\bmb} Y = \left(\frac{\partial}{\partial \beta_d} Y \right)_{d=1,\cdots,D}$ and $\displaystyle \nabla_{\bmb}L = \left(\frac{\partial}{\partial \beta_d}L\right)_{d=1,\cdots,D}$.  
This yields
\begin{equation}
\left\{
\begin{aligned}
Y_0(\rho(t,\bmx), \bmm(t,\bmx), \bmx)
            &=L_0(\rho(t,\bmx),\bmm(t,\bmx))+\lambda_EF'_E(\rho(t,\bmx))+\lambda_QQ(\bmx),\\
\nabla_{\bmb}Y(\rho(t,\bmx), \bmm(t,\bmx), \bmx)
            &= \nabla_{\bmb} L(\rho(t,\bmx),\bmm(t,\bmx)), \quad d= 1,\cdots D
\end{aligned}\right.
\end{equation}

To apply FISTA to this problem, we need to compute  
the gradients $\delta_\rho\calE(\rho,\bmm)$, $\delta_{\bmm}\calE(\rho,\bmm)$ and the projection $\proj_{\calC(\rho_0,\rho_1)}(\rho,\bmm).$

\vspace{0.2cm}
\noindent\textbf{Gradient descent.} Let the boundary values $\rho(0,\cdot)=\rho_0$, $\rho(1,\cdot)=\rho_1$ and $\bmm(t,\bmx)\cdot\bmn=0$ for $\bmx\in\partial\Omega$ being fixed.
By variational calculus, we have
\begin{equation}
\begin{aligned}
\delta_\rho\calE(\rho,\bmm)(t,\bmx)&=
\begin{cases}
Y_0(\rho(t,\bmx), \bmm(t,\bmx), \bmx), &t\neq 0,1\\
0, &t=0,1
\end{cases}\\
\delta_{\bmm}\calE(\rho,\bmm)(t,\bmx)&=
\begin{cases}
\nabla_{\bmb}Y(\rho(t,\bmx), \bmm(t,\bmx), \bmx),&\bmx\not\in \partial\Omega\\
0,&\bmx\in\partial\Omega, \quad d = 1,\cdots, D
\end{cases}
\end{aligned}
\label{eqn: gd cts 1}
\end{equation}
Then with step-size $\eta^\upk$, the descent step can be written as
\begin{equation}
    \left(\rho^\kp,\bmm^\kp\right) = \left(\widehat{\rho}^\upk-\eta^\upk\delta_\rho\calE(\widehat{\rho}^\upk,\widehat{\bmm}^\upk), \widehat{\bmm}^\upk-\eta^\upk \delta_{\bmm} \calE(\widehat{\rho}^\upk,\widehat{\bmm}^\upk)\right),
\label{eqn: gd cts def}
\end{equation}

\noindent\noindent\textbf{Projection.} The projection step 
solves the following minimization problem
\begin{equation}
   \left(\rho^\kpp,\bmm^\kpp\right) =  \argmin_{\rho,\bmm\in\calC(\rho_0,\rho_1)} \half\left\|\rho-\rho^\kp\right\|_{L^2}^2 + \half\left\|\bmm-\bmm^\kp\right\|_{L^2}^2.
   \label{eqn: proj cts def}
\end{equation}
Since the boundary values are fixed and boundary conditions are always satisfied, we only need to introduce dual variable $\phi^\kpp(t,\bmx)$ for mass conservation equation $\partial_t\rho+ \divg_{\bmx} \bmm =0$. 
Consider a Lagrangian function
\begin{equation}
\begin{aligned}
    \calL(\rho,\bmm,\phi):&=\half\left\|\rho-\rho^\kp\right\|_{L^2}^2 
    + \half\left\|\bmm-\bmm^\kp\right\|_{L^2}^2 
    + \left\langle \phi,\partial_t\rho+ \divg_{\bmx} \bmm\right\rangle\\
    &=\half\left\|\rho-\rho^\kp\right\|_{L^2}^2 
    + \half\left\|\bmm-\bmm^\kp\right\|_{L^2}^2 
    - \left\langle \partial_t\phi,\rho\right\rangle 
    - \left\langle \nabla_{\bmx}\phi,\bmm \right\rangle \\
    &\quad + \left\langle \phi(1,\cdot),\rho_1 \right\rangle 
    - \left\langle \phi(0,\cdot),\rho_0 \right\rangle.
\end{aligned}
\end{equation}
$\left(\rho^\kpp,\bmm^\kpp,\phi^\kpp\right)$ is the saddle point of $\calL(\rho,\bmm,\phi)$ if and only if
\begin{equation}
\left\{\begin{aligned}
        &\delta_\rho\calL\left(\rho^\kpp,\bmm^\kpp,\phi^\kpp\right) = 0,\\
        &\delta_{\bmm}\calL\left(\rho^\kpp,\bmm^\kpp,\phi^\kpp\right) = 0,\\
        &\delta_\phi\calL\left(\rho^\kpp,\bmm^\kpp,\phi^\kpp\right) = 0.
\end{aligned}\right.
\end{equation}
This yields
\begin{equation}
\left\{\begin{aligned}
\rho^\kpp &= \rho^\kp + \partial_t\phi^\kpp, \\
\bmm^\kpp &= \bmm^\kp + \nabla_{\bmx}\phi^\kpp.
\end{aligned}\right.
\label{eqn: proj cts update}
\end{equation}
and 
\begin{equation}
    \partial_t\rho^\kpp + \divg_{\bmx} \bmm^\kpp = 0.
\label{eqn: proj cts cst}
\end{equation}
Combining \cref{eqn: proj cts update} and \cref{eqn: proj cts cst}, it is clear that the dual variable $\phi^\kpp$ solves the Poisson equation
\begin{equation}
\left\{\begin{aligned}
- \Delta_{t,\bmx}\phi^\kpp(t,\bmx) &=
 \partial_t\rho^\kp(t,\bmx) +  \divg_{\bmx}\bmm^\kp(t,\bmx), & 0\leq t\leq 1, \bmx\in\Omega\\
\partial_t\phi^\kpp (t,\bmx) &= 0, & t = 0,1, \bmx\in\Omega\\
\nabla_{\bmx}\phi^\kpp (t,\bmx)\cdot\bmn &= 0, & 0\leq t\leq 1, \bmx\in\partial\Omega,\\
\end{aligned}\right.
\label{eqn: proj cts phi}
\end{equation}
Therefore, we can obtain the projection \cref{eqn: proj cts def} in two steps: solving the Poisson equation \cref{eqn: proj cts phi} and update $\rho,\bmm$ by \cref{eqn: proj cts update}.

The FISTA algorithm for MFP problem \cref{eqn: mfp cts prob} is summarized in \cref{alg: fista cts}.

\begin{algorithm}[tb]
\caption{FISTA for MFP}
\begin{algorithmic}
\STATE{Parameters}
$\rho_0,\rho_1$
\STATE{Initialization} $\tau^{(1)}=1,$

\quad\quad\quad\quad\quad
$\rho^{(0)}(0,\cdot)=\widehat{\rho}^{(0)}(0,\cdot)=\rho_0,\ \rho^{(0)}(1,\cdot)=\widehat{\rho}^{(0)}(1,\cdot)=\rho_1,\ $

\quad\quad\quad\quad\quad
$\rho^{(0)}(t,\cdot)=\widehat{\rho}^{(0)}(t,\cdot)=1$ for $0<t<1$,

\quad\quad\quad\quad\quad
$\bmm^{(0)}(\cdot,\bmx)\cdot\bmn=\widehat{\bmm}^{(0)}(\cdot,\bmx)\cdot\bmn=0$ for $\bmx\in\partial\Omega$, 

\quad\quad\quad\quad\quad
$\bmm^{(0)}(\cdot,\bmx)=\widehat{\bmm}^{(0)}(\cdot,\bmx)=1$ for $\bmx\in\Omega\backslash\partial\Omega.$
\FOR{$k=0,1,2,\ldots$}
\STATE{\textbf{gradient descent}} 
\begin{equation}
\left\{\begin{aligned}
    \rho^\kp (t,\bmx) = \widehat{\rho}^\upk (t,\bmx) -  \eta^\upk Y_0\left(\widehat{\rho}^\upk(t,\bmx), \widehat{\bmm}^\upk(t,\bmx), \bmx\right),&\\
     0<t<1, & \bmx\in\Omega.\\
    \bmm^\kp (t,\bmx) = \widehat{\bmm}^\upk (t,\bmx) - \eta^\upk \nabla_{\bmm} Y\left(\widehat{\rho}^\upk(t,\bmx), \widehat{\bmm}^\upk(t,\bmx), \bmx\right),&\\
     0\leq t\leq 1, & \bmx\in\Omega\backslash \partial\Omega.
\end{aligned}\right.
\label{eqn: gd cts alg}
\end{equation}
\STATE{\textbf{projection}}\quad solve $\phi^\kpp$ for 
\begin{equation}
\left\{\begin{aligned}
- \Delta_{t,\bmx}\phi^\kpp(t,\bmx) &=
 \partial_t\rho^\kp(t,\bmx) +  \divg_{\bmx}\bmm^\kp(t,\bmx), & 0\leq t\leq 1, \bmx\in\Omega\\
\partial_t\phi^\kpp (t,\bmx) &= 0, & t = 0,1, \bmx\in\Omega\\
\nabla_{\bmx}\phi^\kpp (t,\bmx)\cdot\bmn &= 0, & 0\leq t\leq 1, \bmx\in\partial\Omega,\\
\end{aligned}\right.
\label{eqn: proj cts phi alg}
\end{equation}
and conduct
\begin{equation}
\left\{\begin{aligned}
\rho^\kpp &= \rho^\kp + \partial_t\phi^\kpp, \\
\bmm^\kpp &= \bmm^\kp + \nabla_{\bmx}\phi^\kpp.
\end{aligned}\right.
\label{eqn: proj cts update alg}
\end{equation}

\STATE{\textbf{update}} 
\begin{align}
&\tau^\kpp=\frac{1+\sqrt{1+4\left(\tau^\upk\right)^2}}{2},\nonumber\\
&\omega^\upk=\frac{\tau^\upk-1}{\tau^\kpp},\nonumber\\
&\left(\widehat{\rho}^\kpp, \widehat{\bmm}^\kpp \right) = \left(1+\omega^\upk\right) \left( \rho^\kpp,\bmm^\kpp \right) -\omega^\upk\left( \rho^\upk,\bmm^\upk \right).
\label{eqn: cts update alg}
\end{align}
\ENDFOR
\end{algorithmic}
\label{alg: fista cts}
\end{algorithm}

\begin{remark}
\label{rem: cts lap^-1}
To compute the projection, we need to solve a Poisson equation with Neumann boundary conditions \cref{eqn: proj cts phi}.
Since for any $\bmx\in\Omega$, $\rho^\kp(0,\bmx)=\rho_0(\bmx),\rho^\kp(1,\bmx)=\rho_1(\bmx)$ and for any $t\in[0,1],\bmx\in\partial\Omega$, $\bmm^\kp(t,\bmx)\cdot\bmn=0$,
we have
\begin{equation*}
\begin{aligned}
    &\int_{[0,1]\times\Omega}\partial_t\rho^\kp(t,\bmx) + \divg_{\bmx} \bmm^\kp(t,\bmx)\deri\bmx\deri t\\
   =&\int_{\Omega}\left( \rho^\kp(1,\bmx)-\rho^\kp(0,\bmx) \right)\deri\bmx
    +\int_0^1\int_{\partial\Omega} \bmm^\kp(t,\bmx)\cdot \deri \bm{s} \deri t\\
   =&\int_{\Omega}\left( \rho^\kp_1(\bmx)-\rho^\kp_0(\bmx) \right)\deri\bmx\\
   =&0,
\end{aligned}
\end{equation*}
This means \cref{eqn: proj cts phi} is solvable and the solution is unique up to a constant. In addition, the constant does not count in projection step because in \cref{eqn: proj cts update}, we only need $\partial_t\phi^\kpp$ and $\nabla_{\bmx}\phi^\kpp$.
Therefore the projection step is well-defined.
\end{remark}

\subsection{Discretization and Implementation}
\label{sec: disct and imple}

\begin{figure}[h]
\centering  
\includegraphics[width=.7\linewidth]{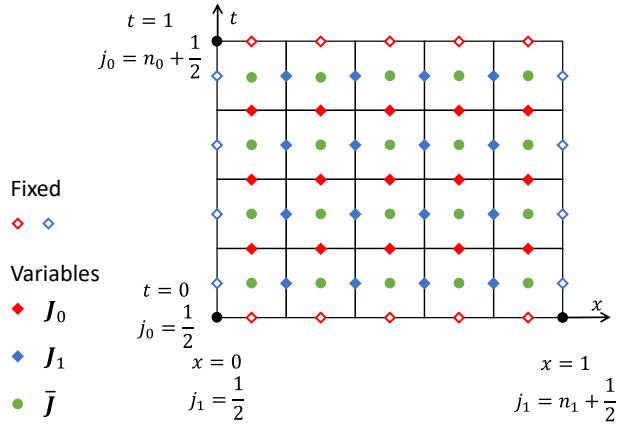}
\vspace{-1cm}
\caption{Illustration of staggered grids for the case $d=1$.}
\label{fig: grid mfp}
\end{figure}

For convenience, we here assume $\Omega = [0,1]^D$. Then the boundary conditions of $\bmm = (m_1,\cdots,m_D)$ is provided as:
\begin{equation*}
    m_d(t,\bmx)=0,\text{ if }x_d=0, \text{or}~ 1, \quad \text{for} ~d = 1,\cdots, D
\end{equation*}
Consider a uniform grid with $n_0$ segments on time interval $[0,1]$ and $n_d$ segments on the $d$-th space dimension. Namely, the mesh size on each dimension is $\Delta_d=\frac{1}{n_d}, d=0,\cdots,D$, and the staggered grid points are $t_j=(j-\half)\Delta_0,(x_d)_j=(j-\half)\Delta_d$.
We use a multi-dimensional index vector $\bmj=(j_0,j_1,\cdots,j_D)$ to indicate a grid point $(t_{j_0},\bmx_{\bmj})$, where $\bmx_{\bmj}:=\left((x_1)_{j_1},\cdots,(x_D)_{j_D}\right)$. We further write $u_{\bmj}:=u\left(t_{j_0},\bmx_{\bmj}\right)$ the value of function $u$ on the grid point and $U_{\bmj}$ the proposed numerical approximation of $u_{\bmj}$. Our discretization of $\rho$ and $\bmm$ defined on different staggered grids. For convenience, we list the following index sets:
\begin{equation*}
\begin{aligned}
    &J_d := \left\{\frac{3}{2},\frac{5}{2},\cdots,n_d-\frac{1}{2} \right\},\\
    &\Jbar_d := \{1,2,\cdots,n_d\},\\
    &\bmJ_d := \Jbar_0\times\Jbar_1\times\cdots\times \Jbar_{d-1}\times J_d\times\Jbar_{d+1}\times\cdots\Jbar_D,\\
    &\overline{\bmJ}:=\Jbar_0\times\Jbar_1\times\cdots\times\Jbar_D.\\
\end{aligned}
\end{equation*}
\cref{fig: grid mfp} illustrates an 1D example, where $n_0=4,n_1=5$ and grid points related to $\bmJ_0$, $\bmJ_1$ and $\overline{\bmJ}$ are annotated as red solid diamonds, blue solid diamonds and green solid dots, respectively. 

We use $\Rho, \bmM$ and $\Phibar$ to denote the discretization of $\rho, \bmm$ and $\phi$, respectively. They are defined as: 
\begin{equation*}
\begin{aligned}
    \Rho&:=\{\Rho_{\bmj}\}_{\bmj\in\bmJ_0}\in\bbV_{0}:=\bbR^{(n_0-1)\times n_1\times\cdots\times n_D},\\
    M_d &:= \{(M_d)_{\bmj}\}_{\bmj\in\bmJ_d}\in\bbV_{d}:=\bbR^{n_0\times n_1\times\cdots\times n_{d-1}\times(n_d-1)\times n_{d+1}\times\cdots\times n_D},\\
    \bmM &:= \{M_d\}_{d=1,2,\cdots,D}\in\bbV_{1}\times\cdots\bbV_{D},\\
     \Phibar &:=\{\Phibar_{\bmj}\}_{\bmj\in\overline{\bmJ}}\in\overline{\bbV}:=\bbR^{n_0\times n_1\times \cdots\times n_D}.
\end{aligned}
\end{equation*}
Moreover, we also define: 
\begin{equation*}
\begin{aligned}
    &\Rhobar:=\{\Rhobar_{\bmj}\}_{\bmj\in\overline{\bmJ}}\in\overline{\bbV}:=\bbR^{n_0\times n_1\times \cdots\times n_D},\\
    &\Mbar_d:=\{(\Mbar_d)_{\bmj}\}_{\bmj\in\overline{\bmJ}}\in\overline{\bbV},\\
    &\overline{\bmM}:=\{\Mbar_d\}_{d=1,2,\cdots,D}\in\overline{\bbV}^D.
\end{aligned}
\end{equation*}

Based on the above settings, we next discuss details of computing the objective value, implementing gradient descent and conducting the projection step.

\vspace{0.2cm}
\noindent\textbf{Objective value.}
To compute objective function, we need the value of $\rho(t,\bmx)$ and $\bmm(t,\bmx)$ on the same point $(t,\bmx)$. While $\Rho,M_d$ are defined on different grids, a natural idea is to transform them to the same central grid $\overline{\bmJ}$ first. For convenience, let $M_0\equiv\Rho$. We can define the average operators as:
\begin{equation*}
\begin{aligned}
    &\calA_d:\bbV_{d}\to\overline{\bbV}, \quad  M_d\mapsto\Mbar_d =\calA_d(M_d) , \quad \text{for}~d=0,1,\cdots,D,\\ 
    &\quad (\Mbar_d)_{\bmj} := \begin{cases}
    \displaystyle\half(M_d)_{\bmj+\frac{\bme_d}{2}},& j_d=1,\vspace{0.1cm}\\
    \displaystyle\half\left( (M_d)_{\bmj+\frac{\bme_d}{2}} + (M_d)_{\bmj-\frac{\bme_d}{2}} \right),& j_d=2,3,\cdots,n_d-1,\vspace{0.1cm}\\
    \displaystyle\half(M_d)_{\bmj-\frac{\bme_d}{2}},& j_d=n_d.\vspace{0.1cm}
    \end{cases}\quad \forall \bmj\in\overline{\bmJ} \\
    &\calA:\bbV_{1}\times\cdots\bbV_{D}\to\overline{\bbV}^D, \quad \bmM\mapsto\{\calA_d(M_d)\}_{d=1,\cdots,D}.\\
\end{aligned}
\end{equation*}
where $\bme_d\in\bbR^{D+1}$ has $1$ in the $(d+1)$-th entry and $0$ elsewhere, $ d=0,\cdots,D$. The boundary conditions of $\bmM$ are implicitly included in the average operator. We further define $\Rhobar_{\calA}\in\overline{\bbV}$  to indicate density boundary conditions,
\begin{equation*}
\begin{aligned}
    (\Rhobar_{\calA})_{\bmj}:=\begin{cases}
    \displaystyle \half\rho_0\left(\bmx_{\bmj}\right),& j_0=1,\vspace{0.1cm}\\
    \displaystyle 0,& j_0=2,3,\cdots,n_0-1,\vspace{0.1cm}\\
    \displaystyle \half\rho_1\left(\bmx_{\bmj}\right),& j_0=n_0.\vspace{0.1cm}
    \end{cases} \quad \quad \forall \bmj\in\overline{\bmJ}\\
\end{aligned}
\end{equation*}
As an example, in \cref{fig: grid mfp}, $\mathcal{A}_0$ maps the red solid dots to the green dots. 
 $\mathcal{A}_1$ maps the blue solid dots to the green dots, and 
the red hollow dots contribute to the non-zero entries of $\Rhobar_{\calA}$.

Now, we are ready to evaluate the objective function by averaging $\Rho$ and $M$ from their staggered grids to the central grid. Namely, We define $\Rhobar,\Mbar\in\overline{\bbV}$ as 
\begin{equation}
\Rhobar:=\calA_0(\Rho)+\Rhobar_{\calA},\quad
\overline{\bmM}:=\calA(\bmM)\label{eqn: average}
\end{equation}
then we approximate the objective value by $$\left(\prod_{d=0}^{D}\Delta_d\right)\calE(\Rho,\bmM),$$
where
\begin{equation}
    \calE(\Rho,\bmM):=\overline{\calE}(\Rhobar,\overline{\bmM}):=\sum_{\bmj\in\overline{\bmJ}}
Y\left( \Rhobar_{\bmj}, \overline{\bmM}_{\bmj},\bmx_{\bmj} \right) . 
\label{eqn: mfp obj func disct}
\end{equation}

\noindent\textbf{Gradient descent.}
To fulfil gradient descent, we first average $(\Rho,\bmM)$ from different grids $\bmJ_d$ to grid $\overline{\bmJ}$ by \cref{eqn: average} and compute gradient values pointwisely
\begin{equation}
\begin{aligned}
&\partial_{\Rhobar}\overline{\calE}(\Rhobar,\overline{\bmM}) 
:= \left\{Y_0\left( \Rhobar_{\bmj}, \overline{\bmM}_{\bmj},\bmx_{\bmj} \right)\right\}_{\bmj\in\overline{\bmJ}},\\
&\partial_{\Mbar_d}\overline{\calE}(\Rhobar,\overline{\bmM}) 
:= \left\{Y_d\left( \Rhobar_{\bmj}, \overline{\bmM}_{\bmj},\bmx_{\bmj} \right)\right\}_{\bmj\in\overline{\bmJ}},\quad d=1,\cdots,D,\\
&\partial_{\bm{\Mbar}}\overline{\calE}(\Rhobar,\overline{\bmM}) 
:= \left\{\partial_{\Mbar_d}\overline{\calE}(\Rhobar,\overline{\bmM})\right\}_{d=1,\cdots,D}.\\
\end{aligned}
\label{eqn: gd disct 1}
\end{equation}
Then we average gradient values back to different grids $\bmJ_d$.
Defining another sets of average operator as
\begin{equation*}
\begin{aligned}
    &\calA_d^*:\overline{\bbV}\to\bbV_{d}, \Mbar_d\mapsto M_d, \ 
    (M_d)_{\bmj} := \half\left( (\Mbar_d)_{\bmj+\frac{\bme_d}{2}} + (\Mbar_d)_{\bmj-\frac{\bme_d}{2}} \right),\\
    &\calA^*:\overline{\bbV}^D\to\bbV_{1}\times\cdots\bbV_{D}, \overline{\bmM}\mapsto\{\calA_d^*(\Mbar_d)\}_{d=1,\cdots,D},\\
\end{aligned}
\end{equation*}
we obtain desired gradient values:
\begin{equation}
\begin{aligned}
&\partial_{\Rho}\calE(\Rho,\bmM) = \calA_0^*\left( \partial_{\Rhobar}\overline{\calE}(\Rhobar,\overline{\bmM}) \right),\\
&\partial_{\bmM}\calE(\Rho,\bmM)= \calA^*\left( \partial_{\bm{\Mbar}}\overline{\calE}(\Rhobar,\overline{\bmM}) \right).
\end{aligned}
\label{eqn: gd disct 2}
\end{equation}
Combining \cref{eqn: gd disct 1},\cref{eqn: gd disct 2}, we can implement gradient descent step \cref{eqn: gd cts alg} on discrete meshes by:
\begin{equation}
    \left(\Rho^\kp,\bmM^\kp\right) =
    \left(\widehat{\Rho}^\upk,\widehat{\bmM}^\upk\right)
    -\eta^\upk\left(\partial_\Rho\calE(\widehat{\Rho}^\upk,\widehat{\bmM}^\upk), \partial_{\bmM} \calE(\widehat{\Rho}^\upk,\widehat{\bmM}^\upk)\right),
\label{eqn: gd disct update}
\end{equation}

\noindent\textbf{Projection. }
To compute projection, we use the finite difference method to discretize the corresponding differential operators in the PDE constraint. We first define discrete partial derivative:
\begin{equation*}
\begin{aligned}
    \calD_d:&\bbV_{d}\to\overline{\bbV}, \quad M_d\mapsto \calD_d(M_d), \quad \text{for}~ d = 0,\cdots,D\\ 
    &(\calD_d(M_d))_{\bmj}:=\begin{cases}
    \displaystyle \frac{1}{\Delta_d} (M_d)_{\bmj+\frac{\bme_d}{2}},
    & j_d=1,\vspace{0.1cm}\\
    \displaystyle \frac{1}{\Delta_d} \left( (M_d)_{\bmj+\frac{\bme_d}{2}} - (M_d)_{\bmj-\frac{\bme_d}{2}} \right), 
    & j_d=2,3,\cdots,n_d-1,\vspace{0.1cm}\\
    \displaystyle -\frac{1}{\Delta_d}(M_d)_{\bmj-\frac{\bme_d}{2}},
    & j_d=n_d,\vspace{0.1cm}
    \end{cases}\\
\end{aligned}
\end{equation*}
discrete divergence:
\begin{equation*}
\begin{aligned}
    \Divg:\bbV_{0}\times\bbV_{1}\times\cdots\bbV_{D}\to\overline{\bbV}, 
    \quad (\Rho,\bmM)\mapsto \calD_0(\Rho)+\sum_{d=1}^D \calD_d(M_d),\\
\end{aligned}
\end{equation*}
and the term $\Rhobar_{\calD}\in\overline{\bbV}$ to impose boundary conditions:
\begin{equation*}
\begin{aligned}
    (\Rhobar_{\calD})_{\bmj}:=\begin{cases}
    \displaystyle -\frac{1}{\Delta_0}\rho_0\left(\bmx_{\bmj}\right),& j_0=1,\vspace{0.1cm}\\
    \displaystyle 0,& j_0=2,3,\cdots,n_0-1,\vspace{0.1cm}\\
    \displaystyle \frac{1}{\Delta_0}\rho_1\left(\bmx_{\bmj}\right),& j_0=n_0.\vspace{0.1cm}
    \end{cases}, \\
\end{aligned}
\end{equation*}
Then the RHS of first equation in \cref{eqn: proj cts phi alg} can be approximated by $$\Divg\left(\Rho^\kp,\bmM^\kp\right)+\Rhobar_{\calD}.$$
We approximate $\partial_d$ with a central difference and $\partial_{dd}$ with a three-point stencil finite difference. By homogenenous Neumann boundary condition, we have discrete second-order derivative operators 
\begin{equation*}
\begin{aligned}
    \calD_{dd}: & \overline{\bbV}\to\overline{\bbV},  \quad\Phibar\mapsto \calD_{dd}(\Phibar), \\ 
    &\left(\calD_{dd}(\Phibar)\right)_{\bmj}:=\begin{cases}
    \displaystyle \frac{1}{\Delta_d^2}\left( -\Phibar_{\bmj}+\Phibar_{\bmj+\bme_d} \right),
    & j_d=1,\vspace{0.1cm}\\
    \displaystyle \frac{1}{\Delta_d^2}\left( \Phibar_{\bmj-\bme_d}-2\Phibar_{\bmj}+\Phibar_{\bmj+\bme_d} \right),
    & j_d=2,3,\cdots,n_d-1,\vspace{0.1cm}\\
    \displaystyle \frac{1}{\Delta_d^2}\left( \Phibar_{\bmj-\bme_d}-\Phibar_{\bmj}\right),
    & j_d=n_d,\vspace{0.1cm}
    \end{cases}\\
    \Lap:&\overline{\bbV}\to\overline{\bbV}, \quad \Phibar\mapsto \calD_{00}(\Phibar)+\sum_{d=1}^D \calD_{dd}(\Phibar).\\
    \end{aligned}
\end{equation*}
The Poisson equation \cref{eqn: proj cts phi alg} on grids is therefore
\begin{equation}
  -\Lap\left(\Phibar^\kpp\right) = \Divg\left(\Rho^\kp,\bmM^\kp\right)+\Rhobar_{\calD},
\label{eqn: proj disct phi}
\end{equation}

Defining another set of derivative operators
\begin{equation*}
\begin{aligned}
    &\calD_d^*:\overline{\bbV}\to\bbV_{d}, \Phibar\mapsto \calD_d^*(\Phibar), \ 
    (\calD_d^*(\Phibar))_{\bmj}:=
    \frac{1}{\Delta_d}\left( (\Phibar)_{\bmj+\frac{\bme_d}{2}} - \Phibar_{\bmj-\frac{\bme_d}{2}} \right)\\
    &\Grad:\overline{\bbV}\to\bbV_{0}\times\bbV_{1}\times\cdots\bbV_{D}, \Phibar\mapsto \left\{ \calD_d^*(\Phibar) \right\}_{d=0,1,\cdots,D},
\end{aligned}
\end{equation*}
we obtain the second step of projection, the discretization of \cref{eqn: proj cts update alg}:
\begin{equation}
\left( \Rho^\kpp, \bmM^\kpp \right) 
= \left( \Rho^\kp, \bmM^\kp \right)+\Grad\left(\Phibar^\kpp\right).
\label{eqn: proj disct update}
\end{equation}

Combining above ingredients, we have FISTA for MFP on discrete mesh summarized in \cref{alg: fista disct}. 
\begin{algorithm}
\caption{FISTA for MFP on discrete mesh}
\begin{algorithmic}
\STATE{Parameters}
$\rho_0,\rho_1$
\STATE{Initialization} $\tau^{(1)}=1,$
$\Rho^{(0)}=\widehat{\Rho}^{(0)}=\mathbf{1}$, and
$M_d^{(0)}=\widehat{M_d}^{(0)}=\mathbf{1}$.

\FOR{$k=0,1,2,\ldots$}
\STATE{\textbf{gradient descent}} 
\begin{equation}
\left\{\begin{aligned}
    &\Rho^\kp =
    \widehat{\Rho}^\upk
    -\eta^\upk\partial_\Rho\calE(\widehat{\Rho}^\upk,\widehat{\bmM}^\upk),\\
    &\bmM^\kp =
    \widehat{\bmM}^\upk
    -\eta^\upk \partial_{\bmM} \calE(\widehat{\Rho}^\upk,\widehat{\bmM}^\upk),
\end{aligned}\right.
\label{eqn: gd disct update alg}
\end{equation}

\STATE{\textbf{projection}}\quad solve $\Phibar^\kpp$ for 
\begin{equation}
  -\Lap\left(\Phibar^\kpp\right) = \Divg\left(\Rho^\kp,\bmM^\kp\right)+\Rhobar_{\calD},
\label{eqn: proj disct phi alg}
\end{equation}
and conduct
\begin{equation}
\left( \Rho^\kpp, \bmM^\kpp \right) 
= \left( \Rho^\kp, \bmM^\kp \right)+\Grad\left(\Phibar^\kpp\right).
\label{eqn: proj disct update alg}
\end{equation}

\STATE{\textbf{update}} 
\begin{align}
&\tau^\kpp=\frac{1+\sqrt{1+4(\tau^\upk)^2}}{2},\nonumber\\
&\omega^\upk=\frac{\tau^\upk-1}{\tau^\kpp},\nonumber\\
&\left(\widehat{\Rho}^\kpp, \widehat{\bmM}^\kpp \right) = \left(1+\omega^\upk\right) \left( \Rho^\kpp,\bmM^\kpp \right) -\omega^\upk\left( \Rho^\upk,\bmM^\upk \right). 
\label{eqn: disct update alg}
\end{align}
\ENDFOR
\end{algorithmic}
\label{alg: fista disct}
\end{algorithm}

\begin{remark}
\label{rem: disct operator consistency}
The discrete operators $\Divg,\Grad$ and $\Lap$ are consistent in the following sense. 
For space $\overline{\bbV}$ and $\bbV_{0}\times\bbV_{1}\times\cdots\times\bbV_{D}$, if we view the elements $\Phibar$ and $(\Rho,\bmM)$ as long vectors, we can define the inner product as
\begin{equation}
\begin{aligned}
    \left\langle\Phibar^1,\Phibar^2\right\rangle
    &:=\sum_{\bmj\in\overline{\bmJ}}\Phibar^1_{\bmj}\Phibar^2_{\bmj},\\
    \left\langle(\Rho^1,\bmM^1),(\Rho^2,\bmM^2)\right\rangle
    &:=\sum_{\bmj\in\bmJ_0}\Rho^1_{\bmj}\Rho^2_{\bmj}
    +\sum_{d=1}^D\sum_{\bmj\in\bmJ_d}(M_d^1)_{\bmj}(M_d^2)_{\bmj}.\\
\end{aligned}
\label{eqn: disct inner prod}
\end{equation}
and define the induced norm as $\|\cdot\|_F$. 
Then simple calculation shows that for any $\Phibar\in\overline{\bbV}$ and $(\Rho,\bmM)\in\bbV_{0}\times\bbV_{1}\times\cdots\times\bbV_{D}$, the following equation holds
\begin{equation}
\begin{aligned}
    &\Lap\left( \Phibar \right) = \Divg\circ\Grad\left( \Phibar \right),\\
    &\left\langle -\Grad(\Phibar), (\Rho,\bmM) \right\rangle = \left\langle \Phibar,\Divg(\Rho,\bmM) \right\rangle, 
\end{aligned}
\label{eqn: disct operators relation}
\end{equation}
These match the relations between $\divg_{t,\bmx},\nabla_{t,\bmx}$ and $\Delta_{t,\bmx}$ on continuous spaces.
\end{remark}

\begin{remark}
\label{rem: disct lap^-1}
Directly solving the large linear system  \cref{eqn: proj disct phi alg} could be very expensive. Thanks to the special structure of the operator $\Lap$, we can decompose it and solve \cref{eqn: proj disct phi alg} in a more efficient manner based on the cosine transformation.

Recall that $\Delta_d=\frac{1}{n_d}$. For $\bmi\in\overline{\bmJ}$, let $\lambda^{\bmi}\in\bbR,\Psi^{\bmi}\in\overline{\bbV}$ be
\begin{equation*}
\begin{aligned}
    &\lambda^{\bmi} = -4\sum_{d=0}^D n_d^2 \sin^2\left(\frac{(i_d-1)\pi}{2n_d}\right),\\
    &\Psi^{\bmi}_{\bmj} = \prod_{d=0}^D\sqrt{\frac{1+\delta_{1i_d}}{n_d}} \cos\left( \left(j_d-\half\right) \frac{(i_d-1)\pi}{n_d} \right),
\end{aligned}
\end{equation*}
then $\displaystyle\left\{ \Psi^{\bmi} \right\}_{\bmi\in\overline{\bmJ}}$ forms an orthogonal basis of $(\overline{\bbV},\|\cdot\|_F)$ and
\begin{equation*}
    \Lap\left( \Psi^{\bmi} \right) = \lambda^{\bmi}\Psi^{\bmi}.
\end{equation*}
Note that $\lambda^{\bmi}=0$ if and only if $\bmi=\bm{1}$, and $\Psi^{\bm{1}}$ has value $\prod_{d=0}^D\sqrt{\frac{1+\delta_{1i_d}}{n_d}}$ in all entries, which implies $\ker(\Lap)=\left\{ C\bm{1}:C\in\bbR \right\}$. This matches the fact that in continuous setting, $\Delta u=0$ if $u\equiv C$ on $\Omega=[0,1]^{D+1}$.
For any $\Phibar\in\overline{\bbV}$, we can define $\Lap^{-1}(\Phibar)$ by requiring
\begin{equation}
    \left\langle \Lap^{-1}(\Phibar),\Psi^{\bm{1}} \right\rangle=0.
\end{equation}
Therefore, we have 
\begin{equation}
    \Lap^{-1}(\Phibar):=\sum_{\bmi\in\overline{\bmJ}\backslash\{\bm{1}\}}\frac{1}{\lambda_{\bmi}} \left\langle \Phibar,\Psi^{\bmi} \right\rangle \Psi^{\bmi}.
\end{equation}
This leads to a discrete cosine transform method to solve \cref{eqn: proj disct phi alg}.

\end{remark}

\begin{figure}[htbp]
\centering
\begin{tikzpicture}
\node (dvp) [process] {discrete optimization problem \cref{eqn: mfp disct prob}};
\node (cvp) [process,below of = dvp, yshift=3cm] {continuous variational problem \cref{eqn: mfp cts prob}};
\node (dalg) [process,left of = dvp, xshift=7cm] {discrete algorithm: \cref{alg: fista disct}};
\node (calg) [process, below of = dalg, yshift=3cm] {continuous algorithm: \cref{alg: fista cts}};
\draw [->,thick,>=stealth] (calg) -- node[anchor=west] {discretize}(dalg);
\draw [->,dashed,thick,>=stealth] (cvp) -- node[anchor=west] {discretize}(dvp);
\draw [->,thick,>=stealth] (cvp) -- node[anchor=south] {optimize} (calg);
\draw [->,dashed,thick,>=stealth] (dvp) -- node[anchor=south] {optimize}(dalg);
\end{tikzpicture}
\label{fig: opt-dis = dis-opt}
\caption{Equivalent approaches to obtain the discrete \cref{alg: fista disct}}
\end{figure}
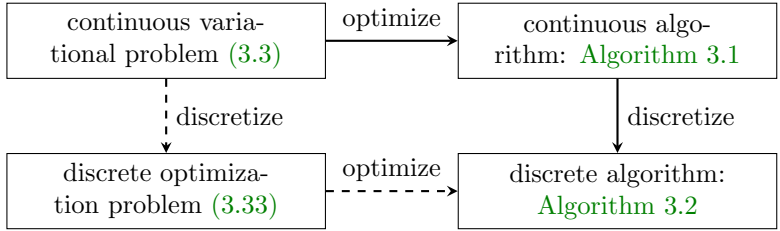

To derive the discrete \cref{alg: fista disct}, we optimize the continuous problem \cref{eqn: mfp cts prob} by \cref{alg: fista cts} then discretize the algorithm. This is a first-optimize-then-discretize approach. We can also consider a first-discretize-then-optimize approach. In fact, using our proposed discretization for MFP, the two approaches lead the same algorithm, as illustrated in \cref{fig: opt-dis = dis-opt}. This is mainly because of the consistent relation of discrete operators discussed in \cref{rem: disct operator consistency}.

Based on previous notations, we discretize the original problem \cref{eqn: mfp cts prob} to
\begin{equation}
    \min_{(\Rho,\bmM)\in\calC(\Rhobar_{\calD})}\calE(\Rho,\bmM):=\sum_{\bmj\in\overline{\bmJ}} 
Y\left( (\mathcal{A}_0(\Rho)+\Rhobar_{\calA})_{\bmj}, \mathcal{A}(\bmM)_{\bmj},\bmx_{\bmj} \right), 
\label{eqn: mfp disct prob}
\end{equation}
where the constraints are linear and the constraint set is convex:
\begin{equation}
    \calC(\Rhobar_{\calD}):=\left\{ \left(\Rho,\bmM\right): \calD_0(\Rho)+\Rhobar_{\calD}+\Divg(\bmM)=\bmO \right\}.
\label{eqn: mfp disct constraint}
\end{equation}
To optimize the problem with FISTA, we first compute gradient.
For any $\Rho,\bmM$, we define the corresponding values on $\overline{\bmJ}$ by $\Rhobar:=\mathcal{A}_0(\Rho)+\Rhobar_{\calA}, \overline{\bmM}:=\mathcal{A}(\bmM)$ be , then 
\begin{equation}
    \calE(\Rho,\bmM)=\sum_{\bmj\in\overline{\bmJ}} 
Y\left( \Rhobar_{\bmj}, \overline{\bmM}_{\bmj},\bmx_{\bmj} \right).
\end{equation}
We will have \cref{eqn: gd disct 1} by taking partial derivatives w.r.t $\Rhobar,\overline{\bmM}$, and then \cref{eqn: gd disct 2} by chain rule. Therefore the gradient descent step is exactly \cref{eqn: gd disct update alg}.
For projection, based on the inner product defined as \cref{eqn: disct inner prod} and induced norm, we can formulate the Lagrangian as 
\begin{equation}
    \calL(\Rho,M,\Phibar):=\half \left\|(\Rho,M)-\left(\Rho^\kp,M^\kp\right)\right\|_F^2 + \left\langle \Phibar^\kpp,\Divg(\Rho,\bmM)+\Rhobar_D\right\rangle.
\end{equation}
Because of the consistency of the discrete operators \cref{eqn: disct operators relation}, we know that \cref{eqn: proj disct phi alg},\cref{eqn: proj disct update alg} computes the projection to $\calC(\Rhobar_{\calD})$. 
Therefore the FISTA algorithm to the discrete MFP problem \cref{eqn: mfp disct prob} is exactly \cref{alg: fista disct}.

\section{Convergence}
\label{sec: convergence}
One major difference between \cref{alg: fista cts} and \cref{alg: fista disct} is that the former one is for the continuous setup while the latter one is for a given discretized mesh grid although both algorithms provide convergence sequences according to the FISTA theory. It is natural to ask if the discretized solution can converge to the continuous solution when the mesh grid size $\Delta_d$ goes to zero. Specifically, 
with a given step-size sequence $\left\{\eta^\upk\right\}_{k}$, let the sequences 
$\left\{ \left(\widehat{\rho}^\upk,\widehat{\bmm}^\upk\right) \right\}_{k}$,
$\left\{ \left(\rho^\upk,\bmm^\upk\right),\left(\rho^\kp,\bmm^\kp\right) \right\}_{k}$ 
be obtained from \cref{alg: fista cts}, 
and $\left\{ \left(\widehat{\Rho}^\upk,\widehat{\bmM}^\upk\right) \right\}_{k}$, 
$\left\{ \left(\Rho^\upk,\bmM^\upk\right),\left(\Rho^\kp,\bmM^\kp\right) \right\}_{k}$ from \cref{alg: fista disct}.
If $\left(\rho^\upk,\bmm^\upk\right)\to(\rho^*,\bmm^*)$ and $\left(\Rho^\upk,\bmM^\upk\right)\to\left(\Rho^*,\bmM^*\right)$ as $k\to\infty$, we would like to check whether $ (\Rho^*,\bmM^*)$ converge to $ (\rho^*,\bmm^*) $ as the mesh grid size converge to zero. In this section, we theoretically analyze  and provide a positive answer to this question. To the best of our knowledge, this is for the first time to examine the discretization error between optimization path of the continuous variational MFP and its discretized optimization counterpart.

We first introduce some notations. 
With given step-size sequence $\{\eta^\upk\}_k$, 
let $\left\{ \left(\widehat{\rho}^\upk,\widehat{\bmm}^\upk\right) \right\}_{k}$, 
$\left\{ \left(\rho^\upk,\bmm^\upk\right),\left(\rho^\kp,\bmm^\kp\right) \right\}_{k}$ be obtained from \cref{alg: fista cts}.
With the same step-size sequence and initialization
$\Rho^{(0)}=\widehat{\Rho}^{(0)}=\rho^{(0)}_{\bmJ_0}, \quad \bmM^{(0)}=\widehat{\bmM}^{(0)}=\bmm^{(0)}_{\bmJ}$,
let $\left\{ \left(\widehat{\Rho}^\upk,\widehat{\bmM}^\upk\right) \right\}_{k}$, $\left\{ \left(\Rho^\upk,\bmM^\upk\right),\left(\Rho^\kp,\bmM^\kp\right) \right\}_{k}$ be obtained from \cref{alg: fista disct}.
For any index set $\bmJ_d,\overline{\bmJ}$, we write the continuous functions $\rho$ and $\bmm$ evaluating on corresponding discrete grids as
\begin{equation*}
\begin{aligned}
    &\rho_{\bmJ_0}:=\left\{ \rho_{\bmj} \right\}_{\bmj\in\bmJ_0},
    \quad (m_d)_{\bmJ_d}:=\left\{ (m_d)_{\bmj} \right\}_{\bmj\in\bmJ_d},
    \quad \bmm_{\bmJ}:=\left\{ (m_d)_{\bmJ_d} \right\}_{d=1,2,\cdots,D},\\
    &\rho_{\overline{\bmJ}}:=\left\{ \rho_{\bmj} \right\}_{\bmj\in\overline{\bmJ}},
    \quad (m_d)_{\overline{\bmJ}}:=\left\{ (m_d)_{\bmj} \right\}_{\bmj\in\overline{\bmJ}},
    \quad \bmm_{\overline{\bmJ}}:=\left\{ (m_d)_{\overline{\bmJ}} \right\}_{d=1,2,\cdots,D}.
\end{aligned}
\end{equation*}
Let $m_0\equiv\rho , M_0\equiv\Rho$. For any $ k=0,\half,1,1+\half,\cdots $, we define the error on grid points $\bmJ_d$ by 
\begin{equation*}
\begin{aligned}
     &E_d^\upk:=M_d^\upk - \left(m_d^\upk\right)_{\bmJ_d}, \quad d = 0,\cdots,D\\
    &\bmE^\upk_{\bmm}:=\left\{ E_d^\upk \right\}_{d=1,\cdots,D},\\
    &\bmE^\upk:=\left\{ E_0^\upk, \bmE^\upk_{\bmm} \right\}.
\end{aligned}
\end{equation*}
Similarly, for $ k=0,1,2,\cdots $, we define
\begin{equation*}
\begin{aligned}
     &\widehat{E}_d^\upk:=\widehat{M}_d^\upk - \left(\widehat{m_d}^\upk\right)_{\bmJ_d},\quad d = 0,\cdots,D\\
     &\widehat{E}_{\bmm}^\upk:=\left\{ \widehat{E}_d^\upk \right\}_{d=1,\cdots,D},\\
    &\widehat{\bmE}^\upk:=\left\{ \widehat{E}_0^\upk, \widehat{\bmE}^\upk_{\bmm}  \right\}\\
     &E_\phi^\upk:=\Phibar^\upk - \phi^\upk_{\overline{\bmJ}},\\
\end{aligned}
\end{equation*}
Recall that in \cref{rem: disct lap^-1}, we introduce induced norm $\|\cdot\|_F$ on space $\bbV_{0}\times\bbV_{1}\times\cdots\times\bbV_{D}$ and $\overline{\bbV}$ as
\begin{equation*}
\begin{aligned}
    &\left\| \bmE \right\|_F:=\left( \sum_{d=0}^D\sum_{\bmj\in\bmJ_d}(E_d)_{\bmj}^2 \right)^\half,
    &\left\| E_\phi \right\|_F:=\left(\sum_{\bmj\in\overline{\bmJ}}(E_{\phi})^2_{\bmj} \right)^\half.\\
\end{aligned}
\end{equation*}
We here define 2-norm $\|\cdot\|_2$ as 
\begin{equation*}
\begin{aligned}
    &\|\cdot\|_2:=\left(\prod_{d=0}^D\Delta_d\right)^{\half}
    \|\cdot\|_F.\\
\end{aligned}
\end{equation*}

Next we propose several assumptions before stating the main theorem.
\begin{assumption}
Let $\rho_0,\rho_1$ be given initial and terminal densities. With above notations, we assume the following conditions hold for any $k=0,1,\cdots$,
\begin{enumerate}
    \item $\rho_0,\rho_1,\rho^\upk,\bmm^\upk,\rho^\kp,\bmm^\kp,\widehat{\rho}^\upk,\widehat{\bmm}^\upk,$ are $C^2$,
    \item There exist $\underline{\rho} 
    \leq \overline{\rho}$, such that $\widehat{\rho}^\upk(t,\bmx), \widehat{\Rho}^\upk_{\bmj} \in[\underline{\rho},\overline{\rho}]$,
    \item $\widehat{\bmm}^\upk(t,\bmx), \widehat{M}^\upk_{\bmj}\in\Omega_{\bmm}\subset\bbR^D$,
    \item $Y_d$'s are $C_{Y}$-Lipschitz continuous on $[\underline{\rho},\overline{\rho}]\times\Omega_{\bmm}\times[0,1]^D$, 
    i.e. for $d=0,\cdots,D$ and any $(\beta_0^1,\bmb^1,\bmx^1),(\beta_0^2,\bmb^2,\bmx^2)\in[\underline{\rho},\overline{\rho}]\times\Omega_{\bmm}\times[0,1]^D$,
    \begin{equation}
    \begin{aligned}
    \left| Y_d(\beta_0^1,\bmb^1,\bmx^1)-Y_d(\beta_0^2,\bmb^2,\bmx^2) \right|\leq C_Y\left\| (\beta_0^1,\bmb^1,\bmx^1)-(\beta_0^2,\bmb^2,\bmx^2) \right\|_1.
    \end{aligned}
    \end{equation}
\end{enumerate}
\label{assumption}
\end{assumption}

\begin{remark}
\cref{assumption} is accessible for very general cases. 
In fact, when $\rho_0,\rho_1,\rho^{(0)},\bmm^{(0)}$ are $C^2$ and $Y$ is $C^1$, one can show that assumption 1 holds by induction on $k$. 
And assumptions 2 and 3 are true as long as $\left\{\rho^\upk,\bmm^\upk\right\}$ and $\left\{\Rho^\upk,\bmM^\upk\right\}$ converges. 
With a typical choice $Y(\beta_0,\bmb,\bmx)=L(\beta_0,\bmb)$ where $L$ is defined in \cref{eqn: dynamical term L}, we retrieve the optimal transport problem. Both \cref{eqn: mfp cts prob} and \cref{eqn: mfp disct prob} have unique minimizers $\{\rho^*,m^*\}$ and $\{\Rho^*,M^*\}$ and both algorithms converges. 
If in addition $\rho_0,\rho_1$ are $C^2$ and $\min_{\bmx}\{\rho_0(\bmx),\rho_1(\bmx)\}\geq \underline{\rho}>0$, then \cref{assumption} hold with continuous initialization and carefully chosen step-sizes.
\end{remark}

We now state our main theorem which characterizes the error bound with respect the grid size.
\begin{theorem}
If \cref{assumption} hold for $k=0,1,\cdots$,
then
\begin{equation}
\begin{aligned}
    \left\|\bmE^\upk  \right\|_2 \leq C\left( \sum_{d=0}^D\Delta_d \right)
    =\bigO\left( \sum_{d=0}^D\Delta_d \right).\\
\end{aligned}
\label{eqn: thm1-1}
\end{equation}
Here $C$ is a constant depending on dimension $D$, Lipschitz constant $C_Y$, stepsizes $\{\eta^{(s)}\}_{s=1,\cdots,k}$ and sequences $\{\widehat{\rho}^{(s)}, \widehat{\bmm}^{(s)}\}_{s=1,\cdots,k}$ but it is independent of $\{\Delta_d\}_{d=0,\cdots,D}$.
\label{thm: main}
\end{theorem}

Note that the above theorem analyzes error bounds at each iteration along optimization paths from the continuous setup and its discretized counterpart. Consequently, we can have the following convergence analysis if both sequences from the continuous and discretized optimization converge (i.e. choice of the step size satisfies convergence conditions used in FISTA~\cite{beck2009fast}). 

\begin{corollary} Suppose that $\left\{\left(\Rho^\upk,\bmM^\upk\right)\right\}_k$ and $\left\{\left(\rho^\upk,\bmm^\upk\right)\right\}_k$ satisfy all conditions in \cref{thm: main}.
If in addition, there exist $\left(\Rho^*,\bmM^*\right)$, $\left(\rho^*,\bmm^*\right)$ such that 
$\rho^*\in C^1,m_d^*\in C^1$
and 
\begin{equation}
\begin{aligned}
    &\lim_{k\to\infty}\left\|\left(\Rho^\upk,\bmM^\upk\right)-\left(\Rho^*,\bmM^*\right)\right\|_2 = 0,\\
    &\lim_{k\to\infty}\left\|\left(\rho^\upk,\bmm^\upk\right)-\left(\rho^*,\bmm^*\right)\right\|_{L_2} = 0,
\end{aligned}
\end{equation}
where $\|\cdot\|_{L_2}$ denotes the standard $L_2$-norm in the function space. Let $\displaystyle \Delta = \max_{d=0,\cdots,D}{\Delta_d}$, then
\begin{equation}
\begin{aligned}
    \lim_{\Delta \to 0}\left\|\bmE^*  \right\|_2:= \lim_{\Delta \to 0}
    \left\| (\Rho^*,\bmM^*)-\left(\rho^*_{\bmJ_0},\bmm^*_{\bmJ}\right) \right\|_2 = 0.
\end{aligned}
\label{eqn: thm1-2}
\end{equation}
\label{thm: coro}
\end{corollary}

\begin{proof}

By triangular inequality,
\begin{equation*}
\begin{aligned}
    &\left\|(\Rho^*,\bmM^*)-(\rho_{\bmJ_0}^*,\bmm_{\bmJ}^*)\right\|_2\\
    \leq& \left\|\left(\Rho^\upk,\bmM^\upk\right)-(\Rho^*,\bmM^*)\right\|_2
    +\left\|\bmE^\upk\right\|_2
    +\left\|\left(\rho_{\bmJ_0}^\upk,\bmm_{\bmJ}^\upk\right)-(\rho_{\bmJ_0}^*,\bmm_{\bmJ}^*)\right\|_2\\
\end{aligned}
\end{equation*}
For any $\epsilon>0$, there exists $k_\epsilon$ such that
\begin{equation}
\begin{aligned}
    \left\|\left(\Rho^{(k_\epsilon)},\bmM^{(k_\epsilon)}\right)-\left(\Rho^*,\bmM^*\right)\right\|_2 &\leq\frac{\epsilon}{4},\\
    \left\|\left(\rho^{(k_\epsilon)},\bmm^{(k_\epsilon)}\right)-\left(\rho^*,\bmm^*\right)\right\|_{L_2}&\leq\frac{\epsilon}{4}
\end{aligned}
\label{eqn: coro pf tri}
\end{equation}
By numerical integration, there exists a constant $C_1$ depending on $d,\rho^{(k_\epsilon)},\bmm^{(k_\epsilon)},$ $\rho^*,\bmm^*$ and independent of $\Delta_d$ such that
\begin{equation}
\begin{aligned}
    &\left\|\left(\rho_{\bmJ_0}^{(k_\epsilon)},\bmm_{\bmJ}^{(k_\epsilon)}\right)-\left(\rho_{\bmJ_0}^*,\bmm_{\bmJ}^*\right)\right\|_2^2\\
    \leq& \int_0^1\int_\Omega\left\|\left(\rho^{(k_\epsilon)},\bmm^{(k_\epsilon)}\right)-\left(\rho^*,\bmm^*\right)\right\|_2^2\deri\bmx\deri t 
    +C_1\sum_{d=0}^D\Delta_d,
\end{aligned}
\label{eqn: coro pf converge1}
\end{equation}
By \cref{thm: main}, there exist a constant $C_2$ independent of $\Delta_d$ such that 
$$\left\|\bmE^{(k_\epsilon)}\right\|_2\leq C_2\sum_{d=0}^D\Delta_d.$$
Let $\displaystyle \delta = \frac{\epsilon}{(D+1)(|C_1|+|C_2|)}$.Then for any $\Delta_d$ satisfying $\max_{d=0,\cdots,D}\Delta_d\leq \delta$, we have 
\begin{equation}
\begin{aligned}
    \left|C_1\sum_{d=0}^D\Delta_d\right|
    +\left|C_2\sum_{d=0}^D\Delta_d\right|\leq\frac{\epsilon}{2},
\end{aligned}
\label{eqn: coro pf converge2}
\end{equation}
Combining \cref{eqn: coro pf tri}, \cref{eqn: coro pf converge1} and \cref{eqn: coro pf converge2}, we conclude that for any $\epsilon>0$, there exist $\delta$ with all $\{\Delta_d\}_d$ satisfying $\Delta \leq \delta$ such that $\|\bmE^*\|_2\leq\epsilon.$ 
\qquad
\end{proof}

To prove \cref{thm: main}, we need to establish three lemmas to analyse the error introduced in each main steps of the algorithm. After that, the proof of \cref{thm: main} can be obtained by induction.

\begin{lemma}
If \cref{assumption} hold for $k=0,1,\cdots$,
then 
\begin{equation}
    \left\| \bmE^\kp \right\|_2
\leq C(D,C_Y,\eta^\upk)\left\| \widehat{\bmE}^\upk \right\|_2
+ C\left(\widehat{\rho}^\upk,\widehat{\bmm}^\upk\right)\left(\sum_{d=0}^D\Delta_d\right)
    + \bigO\left(\sum_{d=0}^D\Delta_d^2\right).
\label{eqn: err gd}
\end{equation}

\label{lem: err gd}
\end{lemma}

\begin{proof}

By definition of $E_d^\kp$, we substitute discrete variables in \cref{eqn: gd disct update alg} by the sum of error and continuous variables. This leads to
\begin{equation*}
\begin{aligned}
    \left(E_d^\kp\right)_{\bmj} + \left(m_d^\kp\right)_{\bmj}
    &=\left(\widehat{E}_d^\upk\right)_{\bmj} + \left(\widehat{m}_d^\upk\right)_{\bmj}-\eta^\upk\left(\partial_{M_d}\calE(\widehat{\Rho}^\upk,\widehat{\bmM}^\upk)\right)_{\bmj}
\end{aligned}    
\end{equation*}
From \cref{eqn: gd cts alg}, we have
\begin{equation*}
\begin{aligned}
    \left(m_d^\kp\right)_{\bmj}
    &=\left(\widehat{m}_d^\upk\right)_{\bmj}-\eta^\upk Y_d\left(\widehat{\rho}^\upk_{\bmj},
              \widehat{\bmm}^\upk_{\bmj},
              \bmx_{\bmj}\right).
\end{aligned}    
\end{equation*}
Combining above gives us
\begin{equation*}
\begin{aligned}
    \left(E_d^\kp\right)_{\bmj}
    =\left(\widehat{E}_d^\upk\right)_{\bmj}
    &-\eta^\upk \left[\left(\partial_{M_d}\calE(\widehat{\Rho}^\upk,\widehat{\bmM}^\upk)\right)_{\bmj} 
    - Y_d\left(\widehat{\rho}^\upk_{\bmj},
              \widehat{\bmm}^\upk_{\bmj},
              \bmx_{\bmj}\right)\right] \\
\end{aligned}    
\end{equation*}
Therefore we have the norm estimation 
\begin{equation}
\begin{aligned}
    \left\|\bmE^\kp\right\|_F
    \leq&\left\|\widehat{\bmE}^\upk\right\|_F\\
    &+\eta^\upk\left[\sum_{d=0}^D\sum_{\bmj\in\bmJ_d}
    \left| \left(\partial_{M_d}\calE(\widehat{\Rho}^\upk,\widehat{\bmM}^\upk)\right)_{\bmj} - Y_d\left(\widehat{\rho}^\upk_{\bmj},
              \widehat{\bmm}^\upk_{\bmj},
              \bmx_{\bmj}\right) \right|^2\right]^\half
\end{aligned} 
\label{eqn: err gd1}
\end{equation}
For any $\bmj \in \bmJ_d$, the definition of $\left(\partial_{M_d}\calE(\widehat{\Rho}^\upk,\widehat{\bmM}^\upk)\right)_{\bmj}$ in \cref{eqn: gd disct 2} yields, 
\begin{equation}
\begin{aligned}
    &\left| \left(\partial_{M_d}\calE\left(\widehat{\Rho}^\upk,\widehat{\bmM}^\upk\right)\right)_{\bmj}
    -Y_d\left(\widehat{\rho}^\upk_{\bmj},\widehat{\bmm}^\upk_{\bmj},\bmx_{\bmj}\right) \right| \\
    \leq&\half \left|\left( \partial_{\Mbar_d}\overline{\calE}\left(\overline{\widehat{\Rho}}^\upk,\overline{\widehat{\bmM}}^\upk\right) \right)_{\bmj+\frac{\bme_d}{2}}
    - Y_d\left(\widehat{\rho}^\upk_{\bmj},\widehat{\bmm}^\upk_{\bmj},\bmx_{\bmj}\right)\right|\\
    &+ \half\left| \left(\partial_{\Mbar_d}\overline{\calE}\left(\overline{\widehat{\Rho}}^\upk,\overline{\widehat{\bmM}}^\upk\right)  \right)_{\bmj-\frac{\bme_d}{2}}
    - Y_d\left(\widehat{\rho}^\upk_{\bmj},\widehat{\bmm}^\upk_{\bmj},\bmx_{\bmj}\right) \right|\\
    =&\half \left| Y_d\left( \overline{\widehat{\Rho}}^\upk_{\bmj+\frac{\bme_d}{2}}, \overline{\widehat{\bmM}}^\upk_{\bmj+\frac{\bme_d}{2}}, \bmx_{\bmj+\frac{\bme_d}{2}}\right) 
    - Y_d\left(\widehat{\rho}^\upk_{\bmj},\widehat{\bmm}^\upk_{\bmj},\bmx_{\bmj}\right)\right|\\
    &+ \half \left| Y_d\left( \overline{\widehat{\Rho}}^\upk_{\bmj-\frac{\bme_d}{2}}, \overline{\widehat{\bmM}}^\upk_{\bmj-\frac{\bme_d}{2}}, \bmx_{\bmj-\frac{\bme_d}{2}}\right)
    - Y_d\left(\widehat{\rho}^\upk_{\bmj},\widehat{\bmm}^\upk_{\bmj},\bmx_{\bmj}\right)\right| \\
    \leq&\frac{C_Y}{2} \left\| \left( \overline{\widehat{\Rho}}^\upk_{\bmj+\frac{\bme_d}{2}}, \overline{\widehat{\bmM}}^\upk_{\bmj+\frac{\bme_d}{2}}, \bmx_{\bmj+\frac{\bme_d}{2}}\right) 
    - \left(\widehat{\rho}^\upk_{\bmj},\widehat{\bmm}^\upk_{\bmj},\bmx_{\bmj}\right)\right\|_1\\
    &+ \frac{C_Y}{2} \left\| \left( \overline{\widehat{\Rho}}^\upk_{\bmj-\frac{\bme_d}{2}}, \overline{\widehat{\bmM}}^\upk_{\bmj-\frac{\bme_d}{2}}, \bmx_{\bmj-\frac{\bme_d}{2}}\right)
    - \left(\widehat{\rho}^\upk_{\bmj},\widehat{\bmm}^\upk_{\bmj},\bmx_{\bmj}\right)\right\|_1\\
\end{aligned}
\label{eqn: err grad1}
\end{equation}
Note that $\overline{\widehat{\Rho}}^\upk_{\bmj},\left(\overline{\widehat{M}}_d^\upk\right)_{\bmj}$ can be written as the sum of errors and continuous values:
\begin{equation*}
\begin{aligned}
    \overline{\widehat{\Rho}}^\upk_{\bmj}
    &=\left(\mathcal{A}_0\left( \widehat{E}_0^\upk+\widehat{\rho}^\upk_{\bmJ_0} \right)+\Rhobar_{\calA}\right)_{\bmj}\\
    &=\left(\mathcal{A}_0\left( \widehat{E}_0^\upk \right)\right)_{\bmj}+\left(\mathcal{A}_0\left( \widehat{\rho}^\upk_{\bmJ_0} \right)+\Rhobar_{\calA}\right)_{\bmj}\\
    &=\left(\mathcal{A}_0\left( \widehat{E}_0^\upk \right)\right)_{\bmj}+\widehat{\rho}^\upk_{\bmj}+\bigO(\Delta_0^2)\\
    \left(\overline{\widehat{M}}_d^\upk\right)_{\bmj}
    &=\left(\mathcal{A}_d\left( \widehat{E}_d^\upk \right)\right)_{\bmj}+\left(\widehat{m}_d^\upk\right)_{\bmj}+\bigO(\Delta_d^2),
\end{aligned}    
\end{equation*}
where the last equality in the above two equations are obtained from using Taylor expansion to $\widehat{\rho}^\upk$ and $\widehat{m}_d^\upk$. We further have:
\begin{equation}
\begin{aligned}
    & \left\| \left( \overline{\widehat{\Rho}}^\upk_{\bmj\pm\frac{\bme_d}{2}}, \overline{\widehat{\bmM}}^\upk_{\bmj\pm\frac{\bme_d}{2}}, \bmx_{\bmj\pm\frac{\bme_d}{2}}\right) 
    - \left(\widehat{\rho}^\upk_{\bmj},\widehat{\bmm}^\upk_{\bmj},\bmx_{\bmj}\right)\right\|_1\\
    \leq& \left\| \left( \left(\calA_0\left( \widehat{E}_0^\upk \right)\right)_{\bmj\pm\frac{\bme_d}{2}}, \left(\calA\left(\widehat{\bmE}^\upk_{\bmm}\right)\right)_{\bmj\pm\frac{\bme_d}{2}},\bmO \right) \right\|_1 \\
    &+\left\| \left(\widehat{\rho}^\upk_{\bmj\pm\frac{\bme_d}{2}},\widehat{\bmm}^\upk_{\bmj\pm\frac{\bme_d}{2}},\bmx_{\bmj\pm\frac{\bme_d}{2}}\right)
    - \left(\widehat{\rho}^\upk_{\bmj},\widehat{\bmm}^\upk_{\bmj},\bmx_{\bmj}\right)\right\|_1
    +\bigO\left(\Delta_d^2\right)\\
    \leq&\half\sum_{d'=0}^D \left|\left(\widehat{E}_{d'}^\upk\right)_{\bmj} +\left(\widehat{E}_{d'}^\upk\right)_{\bmj+\bme_d} \right| 
    +\half\sum_{d'=0}^D \left|\left(\widehat{E}_{d'}^\upk\right)_{\bmj} +\left(\widehat{E}_{d'}^\upk\right)_{\bmj-\bme_d} \right| \\ 
    &+C\left(\widehat{\rho}^\upk,\widehat{\bmm}^\upk\right)\Delta_d+\bigO\left(\Delta_d^2\right),\\
\end{aligned}
\label{eqn: err grad2}
\end{equation}
where $ C\left(\widehat{\rho}^\upk,\widehat{\bmm}^\upk\right) =\max\left\{\frac{\partial}{\partial x_d}\widehat{m}_{d'}^\upk(t,\bmx):d'=0,\cdots,D\right\}.$
Combining \cref{eqn: err grad1} and \cref{eqn: err grad2} provides:
\begin{equation*}
\begin{aligned}
    &\left| \left(\partial_{M_d}\calE\left(\widehat{\Rho}^\upk,\widehat{\bmM}^\upk\right)\right)_{\bmj}
    -Y_d\left(\widehat{\rho}^\upk_{\bmj},\widehat{\bmm}^\upk_{\bmj},\bmx_{\bmj}\right) \right| \\
    \leq&\frac{C_Y}{4}\sum_{d'=0}^D\left|\left(\widehat{E}_{d'}^\upk\right)_{\bmj} +\left(\widehat{E}_{d'}^\upk\right)_{\bmj+\bme_d} \right|
    +\frac{C_Y}{4}\sum_{d'=0}^D\left|\left(\widehat{E}_{d'}^\upk\right)_{\bmj} +\left(\widehat{E}_{d'}^\upk\right)_{\bmj-\bme_d} \right|\\
    &+C\left(\widehat{\rho}^\upk,\widehat{\bmm}^\upk\right)\Delta_d+\bigO\left(\Delta_d^2\right).
\end{aligned}
\end{equation*}
and applying the triangle inequality yields: 
\begin{equation*}
\begin{aligned}
    &\left[\sum_{d=0}^D\sum_{\bmj\in\bmJ_d}
    \left| \left(\partial_{M_d}\calE(\widehat{\Rho}^\upk,\widehat{\bmM}^\upk)\right)_{\bmj} 
    - Y_d\left(\widehat{\rho}^\upk_{\bmj},
              \widehat{\bmm}^\upk_{\bmj},
              \bmx_{\bmj}\right) \right|^2\right]^\half\\
    \leq&\left[\sum_{d=0}^D\sum_{\bmj\in\bmJ_d}\sum_{d'=0}^D \frac{C_Y^2}{8}
    \left( \left(\widehat{E}_{d'}^\upk\right)_{\bmj-\bme_d}^2
    +2\left(\widehat{E}_{d'}^\upk\right)_{\bmj}^2
    +\left(\widehat{E}_{d'}^\upk\right)_{\bmj+\bme_d}^2\right)\right]^\half\\
    &+ \left[\sum_{d=0}^D\sum_{\bmj\in\bmJ_d} C^2\left(\widehat{\rho}^\upk,\widehat{\bmm}^\upk\right)\Delta_d^2\right]^\half 
    + \left(\prod_{d=0}^D n_d\right)^\half\bigO\left(\sum_{d=0}^D\Delta_d^2\right)
\end{aligned}    
\end{equation*}
\begin{equation*}
\begin{aligned}
    \leq&C(D,C_Y)\left\|\widehat{\bmE}^\upk\right\|_F\\ 
    &+ \left(\prod_{d=0}^D n_d\right)^\half C\left(\widehat{\rho}^\upk,\widehat{\bmm}^\upk\right)\left(\sum_{d=0}^D \Delta_d\right)
    + \left(\prod_{d=0}^D n_d\right)^\half\bigO\left(\sum_{d=0}^D\Delta_d^2\right),\\
\end{aligned}    
\end{equation*}
Together with estimation \cref{eqn: err gd1}, we have
\begin{equation}
\begin{aligned}
    \left\|\bmE^\kp\right\|_2
    \leq&\left(1+C(D,C_Y,\eta^\upk)\right)
    \left\|\widehat{\bmE}^\upk\right\|_2\\
    &+ C\left(\widehat{\rho}^\upk,\widehat{\bmm}^\upk\right)\left(\sum_{d=0}^D\Delta_d\right)
    + \bigO\left(\sum_{d=0}^D\Delta_d^2\right).
\end{aligned}  
\label{eqn: err gd2}
\end{equation}
Therefore we prove the lemma.
\qquad
\end{proof}

Next we examine the error introduced in projection step. 
\begin{lemma}
Suppose that $\rho_0,\rho_1,\rho^\kp,\bmm^\kp$ are $C^2$,
then $$\left\| \bmE^\kpp \right\|_2\leq 2\left\| \bmE^\kp \right\|_2+ \bigO\left(\sum_{d=0}^D\Delta_d^2\right). $$
\label{lem: err proj}
\end{lemma}

\begin{proof}
By definition of error terms and \cref{eqn: proj disct phi alg}, we have 
\begin{equation}
    -\Lap \left(E_\phi^\kpp+\phi_{\overline{\bmJ}}^\kpp\right)
    =\Divg\left(\bmE^\kp+\left(\rho^\kp_{\bmJ_0},\bmm^\kp_{\bmJ}\right)\right)+\Rhobar_{\calD}.
\label{eqn: err proj1-1}
\end{equation}
Since $\left(\phi^\kpp,\rho^\kp,\bmm^\kp\right)$ satisfies \cref{eqn: proj cts phi alg}, and $\rho^\kp,\bmm^\kp$ are $C^2$, by Taylor expansions, we have
\begin{equation}
    -\Lap \left(\phi_{\overline{\bmJ}}^\kpp\right)
    = \Divg\left(\rho^\kp_{\bmJ_0},\bmm^\kp_{\bmJ}\right)+\Rhobar_D+\bm{\bigO}_1\left(\sum_{d=0}^D\Delta_d^2\right).
\label{eqn: err proj1-2}
\end{equation}
Here $\displaystyle\bm{\bigO}_1\left(\sum_{d=0}^D\Delta_d^2\right)\in\overline{\bbV}$ indicates its entry-wise contribution is  $\displaystyle\bigO\left(\sum_{d=0}^D\Delta_d^2\right)$.
Combining \cref{eqn: err proj1-1} and \cref{eqn: err proj1-2} gives us
\begin{equation}
    -\Lap \left(E_\phi^\kpp\right) = \Divg\left(\bmE^\kp\right)+\bm{\bigO}_1\left(\sum_{d=0}^D\Delta_d^2\right).
\label{eqn: err proj1}
\end{equation}
Similarly, the second step on discrete mesh \cref{eqn: proj disct update alg} gives
\begin{equation}
\begin{aligned}
    \bmE^\kpp + \left(\rho^\kpp_{\bmJ_0},\bmm^\kpp_{\bmJ}\right) 
    =& \bmE^\kp + \Grad\left( E_\phi^\kpp \right)\\
    &+\left(\rho^\kp_{\bmJ_0},\bmm^\kp_{\bmJ}\right) 
     + \Grad\left(\phi_{\overline{\bmJ}}^\kpp \right),
\end{aligned}
\label{eqn: err proj2-1}
\end{equation}
and on continuous setting \cref{eqn: proj cts update alg} gives
\begin{equation}
\begin{aligned}
    \left(\rho^\kpp_{\bmJ_0},\bmm^\kpp_{\bmJ}\right) 
    = \left(\rho^\kp_{\bmJ_0},\bmm^\kp_{\bmJ}\right) 
    + \Grad\left( \phi_{\overline{\bmJ}}^\kpp \right) + \bm{\bigO}_2\left(\sum_{d=0}^D\Delta_d^2\right).
\end{aligned}
\label{eqn: err proj2-2}
\end{equation}
Thus we have:
\begin{equation}
\begin{aligned}
    \bmE^\kpp = \bmE^\kp + \Grad\left( E_\phi^\kpp \right) + \bm{\bigO}_2\left(\sum_{d=0}^D\Delta_d^2\right),
\end{aligned}
\label{eqn: err proj2}
\end{equation}
where $ \displaystyle\bm{\bigO}_2\left(\sum_{d=0}^D\Delta_d^2\right)\in\bbV_{0}\times\bbV_{1}\times\cdots\times\bbV_{D} $ indicates its entry-wise contribution is $\displaystyle\bigO\left(\sum_{d=0}^D\Delta_d^2\right)$.

Combining \cref{eqn: err proj1} and \cref{eqn: err proj2}, we obtain
\begin{equation}
\begin{aligned}
    \bmE^\kpp =& \left(\Id-\Grad\circ\Lap^{-1}\circ\Divg\right)\bmE^\kp \\
    &- \Grad\circ\Lap^{-1}\bm{\bigO}_1\left(\sum_{d=0}^D\Delta_d^2\right)
    + \bm{\bigO}_2\left(\sum_{d=0}^D\Delta_d^2\right).
\end{aligned}
\end{equation}

\textit{Claim: }
$\displaystyle\| \Grad\circ\Lap^{-1}\circ\Divg \|_2\leq 1,$ 
$\displaystyle\|\Grad\circ\Lap^{-1}\|_2\leq \frac{1}{4}.$

Therefore
\begin{equation*}
\begin{aligned}
    \left\|\bmE^\kpp\right\|_2
    &\leq2\left\|\bmE^\kp\right\|_2
    +\frac{1}{4}\left\| \displaystyle\bm{\bigO}_1\left(\sum_{d=0}^D\Delta_d^2\right) \right\|_2 
    +\left\| \displaystyle\bm{\bigO}_2\left(\sum_{d=0}^D\Delta_d^2\right) \right\|_2\\
    &\leq2\left\|\bmE^\kp\right\|_2 + \bigO\left(\sum_{d=0}^D\Delta_d^2\right).
\end{aligned}
\end{equation*}

\textit{Proof of claim:}
%
Recall that $\Lap\left(\Psi^{\bmi}\right)=\lambda^{\bmi}\Psi^{\bmi},\bmi\in\overline{\bmJ}$, where
\begin{equation*}
\begin{aligned}
    &\lambda^{\bmi} = -4\sum_{d=0}^D n_d^2 \sin^2\left(\frac{(i_d-1)\pi}{2n_d}\right),\\
    &\Psi^{\bmi}_{\bmj} = \prod_{d=0}^D\sqrt{\frac{1+\delta_{1i_d}}{n_d}} \cos\left( \left(j_d-\half\right) \frac{(i_d-1)\pi}{n_d} \right),
\end{aligned}
\end{equation*}
and $\left\{\Psi^{\bmi} \right\}_{\bmi\in\overline{\bmJ}}$ forms a basis of $(\overline{\bbV},\|\cdot\|_F).$

For $d=0,1,\cdots,D$, and $\bmi\in\overline{\bmJ},i_d\neq1$, let $\sigma^{d,\bmi}\in\bbR$ and $\bm{\Psi}^{d,\bmi}\in\bbV_0\times\bbV_1\times\cdots\times\bbV_d$ be:
\begin{equation*}
\begin{aligned}
    &\sigma^{d,\bmi}=-2n_d\sin\left( \frac{(i_d-1)\pi}{2n_d} \right),\\
    &\bm{\Psi}^{d,\bmi} = \left\{ \Psi^{d,\bmi}_{d'} \right\}_{d'=0,1,\cdots,D},\\
    &\quad\text{where }\Psi^{d,\bmi}_{d} = \frac{1}{\sigma^{d,\bmi}}\calD_d^*\left(\Psi^{\bmi}\right),
    \quad\Psi^{d,\bmi}_{d'} = \bmO, d'\neq d,
\end{aligned}
\end{equation*}
then 
\begin{equation*}
\begin{aligned}
    \langle\bm{\Psi}^{d,\bmi},\bm{\Psi}^{d',\bmi'} \rangle &= 0,  \text{ if } d\neq d',\\
    \langle\bm{\Psi}^{d,\bmi},\bm{\Psi}^{d,\bmi'} \rangle &= \frac{1}{\sigma^{d,\bmi}\sigma^{d,\bmi'}} \langle\calD_d^*\left(\bm{\Psi}^{\bmi}\right),\calD_d^*(\bm{\Psi}^{\bmi'}) \rangle
    =\frac{\sigma^{d,\bmi'}}{\sigma^{d,\bmi}}\langle\Psi^{\bmi},\Psi^{\bmi'}\rangle
    =\begin{cases}
    0, &\bmi\neq\bmi'\\
    1, &\bmi=\bmi'
    \end{cases}
\end{aligned}
\end{equation*}
i.e. $\left\{ \bm{\Psi}^{d,\bmi} \right\}$ forms an orthonormal basis of $(\bbV_{0}\times\bbV_{1}\times\cdots\times\bbV_{D},\|\cdot\|_F).$

Since $\displaystyle\|\cdot\|_2 =\left(\prod_{d=0}^D\Delta_d\right)^\half \|\cdot\|_F$, we next compute the $\left\|\Grad\circ\Lap^{-1}\circ\Divg \right\|_2$, $\left\|\Grad\circ\Lap^{-1}\right\|_2$ with basis of $(\overline{\bbV},\|\cdot\|_F)$ and $(\bbV_{0}\times\bbV_{1}\times\cdots\times\bbV_{D},\|\cdot\|_F).$
For any basis $\Psi^{\bmi}\in\overline{\bbV}$, 
\begin{equation*}
\begin{aligned}
    &\Grad\circ\Lap^{-1}\left( \Psi^{\bm{1}} \right) = \bmO,\\
    &\Grad\circ\Lap^{-1}\left( \Psi^{\bmi} \right) = 
    \sum_{\substack{d=0,\\i_d\neq1} }^D \frac{\sigma^{d,\bmi}}{\lambda^{\bmi}}\Psi^{d,\bmi},\quad\bmi\neq\bm{1}
\end{aligned}
\end{equation*}
thus when $n_d>1$ for $d=0,\cdots,D$, we have
\begin{equation*}
\begin{aligned}
    &\left\|\Grad\circ\Lap^{-1}\right\|_2
\leq\max_{\bmi\in\overline{\bmJ}\backslash\{\bm{1}\}} 
    \left( \frac{1}{\left(\lambda^{\bmi}\right)^2}\sum_{\substack{d=0,\\i_d\neq1} }^D
    \left(\sigma^{d,\bmi}\right)^2 \right)^{\half}
=&\max_{\bmi\in\overline{\bmJ}\backslash\{\bm{1}\}}
    \frac{1}{\left|\lambda^{\bmi}\right|}
\leq&\frac{1}{4}.
\end{aligned}
\end{equation*}
And for any basis $\Psi^{d,\bmi}\in\bbV_{0}\times\bbV_{1}\times\cdots\times\bbV_{D}$, 
\begin{equation*}
\begin{aligned}
    &\Grad\circ\Lap^{-1}\circ\Divg\left(\Psi^{d,\bmi}\right)\\
    =& \Grad\circ\Lap^{-1}\left( \frac{1}{\sigma^{d,\bmi}}\calD_d\circ \calD_d^*\left(\Psi^{\bmi}\right) \right)
    = \Grad\circ\Lap^{-1}\left( -\sigma^{d,\bmi}\Psi^{\bmi} \right)\\
    =& \Grad\left( -\frac{\sigma^{d,\bmi}}{\lambda^{\bmi}}\Psi^{\bmi} \right)
    =- \sum_{\substack{d'=0,\\i_{d'}\neq 1} }^D \frac{\sigma^{d,\bmi}\sigma^{d',\bmi}}{\lambda^{\bmi}}\Psi^{d',\bmi}
\end{aligned}
\end{equation*}
therefore
\begin{equation*}
\begin{aligned}
    \left\|\Grad\circ\Lap^{-1}\circ\Divg\right\|_2
\leq &\max_{d=0,1,\cdots,D}\max_{\bmi\in\overline{\bmJ},i_d\neq1}
    \left( \left(\frac{\sigma^{d,\bmi}}{\lambda^{\bmi}} \right)^2
    \sum_{\substack{d'=0,\\i_{d'}\neq 1} }^D
    \left(\sigma^{d',\bmi}\right)^2 \right)^{\half}\\
= &\max_{d=0,1,\cdots,D}\max_{\bmi\in\overline{\bmJ},i_d\neq1}
    \left( \frac{\left(\sigma^{d,\bmi}\right)^2}{\left|\lambda^{\bmi}\right|} \right)^{\half}
\leq 1.
\end{aligned}
\end{equation*}
The claim  and thus the lemma are proved.
\qquad
\end{proof}

The last step is to estimate the error introduced in linear interpolation step \cref{eqn: cts update alg},\cref{eqn: disct update alg}.

\begin{lemma} 
\begin{equation*}
\begin{aligned}
    \left\|\widehat{\bmE}^\kpp\right\|_2
    &\leq \left|1+\omega^\upk\right|\left\|\bmE^\kpp\right\|_2 + \left|\omega^\upk\right|\left\|\bmE^\upk\right\|_2.\\
\end{aligned}
\end{equation*}

\label{lem: err update}
\end{lemma}

\begin{proof}
By definition of error terms and \cref{eqn: disct update alg}
\begin{equation*}
\begin{aligned}
    &\widehat{E}_d^\kpp+\left(\widehat{m}_d^\kpp\right)_{\bmJ_d}\\
    =& \left(1+\omega^\upk\right) \left(E_d^\kpp+\left(m_d^\kpp\right)_{\bmJ_d} \right) - \omega^\upk \left(E_d^\upk+\left(m_d^\upk\right)_{\bmJ_d} \right),
\end{aligned}
\end{equation*}
and by \cref{eqn: cts update alg}
\begin{equation*}
    \left(\widehat{m}_d^\kpp\right)_{\bmJ_d}
    = \left(1+\omega^\upk\right)\left(m_d^\kpp\right)_{\bmJ_d} - \omega^\upk\left(m_d^\upk\right)_{\bmJ_d}.
\end{equation*}
Therefore we have
\begin{equation*}
\begin{aligned}
    \widehat{\bmE}^\kpp
    &= \left(1+\omega^\upk\right)\bmE^\kpp - \omega^\upk \bmE^\upk.
\end{aligned}
\end{equation*}
By triangular inequality, the lemma is proved.
\qquad
\end{proof}

With \cref{lem: err gd}-\cref{lem: err update}, we can show \cref{thm: main} by induction.

\noindent\textbf{Proof of \cref{thm: main}: } 
\begin{proof} We first restate results from \cref{lem: err gd}-\cref{lem: err update}:
\begin{align*}
    &\left\| \bmE^\kpp \right\|_2\leq 2\left\| \bmE^\kp \right\|_2+ \bigO\left(\sum_{d=0}^D\Delta_d^2\right),\\
    &\left\| \bmE^\kp \right\|_2
\leq C(D,C_Y,\eta^\upk)\left\| \widehat{\bmE}^\upk \right\|_2
+ C\left(\widehat{\rho}^\upk,\widehat{\bmm}^\upk\right)\left(\sum_{d=0}^D\Delta_d\right)
    + \bigO\left(\sum_{d=0}^D\Delta_d^2\right),\\
    &\left\| \widehat{\bmE}^\upk \right\|_2 \leq \left|1+\omega^{(k-1)}\right|\left\| \bmE^\upk \right\|_2 + \left|\omega^{(k-1)}\right|\left\| \bmE^{(k-1)} \right\|_2.
\end{align*}
From these, we obtain
\begin{align}
    &\left\| \bmE^{(1)} \right\|_2
    \leq C\left\| \widehat{\bmE}^{(0)} \right\|_2
    +C\sum_{d=0}^D\Delta_d
    + \bigO\left(\sum_{d=0}^D\Delta_d^2\right),\\
    &\left\| \bmE^\kpp \right\|_2
    \leq C\left\| \bmE^\upk \right\|_2 
    + C\left\| \bmE^{(k-1)} \right\|_2
    + C\sum_{d=0}^D\Delta_d
    + \bigO\left(\sum_{d=0}^D\Delta_d^2\right),\quad k\geq 1
    \label{eqn:error recursive}
\end{align}
where $C$ depends on $D,C_Y,\eta^\upk,\rho^\kp,\bmm^\kp,\widehat{\rho}^\upk,\widehat{\bmm}^\upk$.

The initialization gives us
\begin{equation*}
\begin{aligned}
    \left\|\bmE^{(0)}  \right\|_2= 0, \left\|\widehat{\bmE}^{(0)}\right\|_2=0,
\end{aligned}
\end{equation*}
Then based on \cref{eqn:error recursive}, it is straightforward to show \cref{eqn: thm1-1} by  applying induction on $k$ directly.  
\qquad
\end{proof}

\section{Generalization and Acceleration}
\label{sec: acc and general}
In this section, we generalize the proposed algorithm to solve potential MFG problems. Moreover, we also discuss how to use multilevel and multigrid strategies to speed up our algorithm. 

\subsection{Generalization to Potential MFG}
\label{sec: generalization mfg}

\begin{figure}[htbp]
\centering
\vspace{-1cm}
\includegraphics[width=10cm]{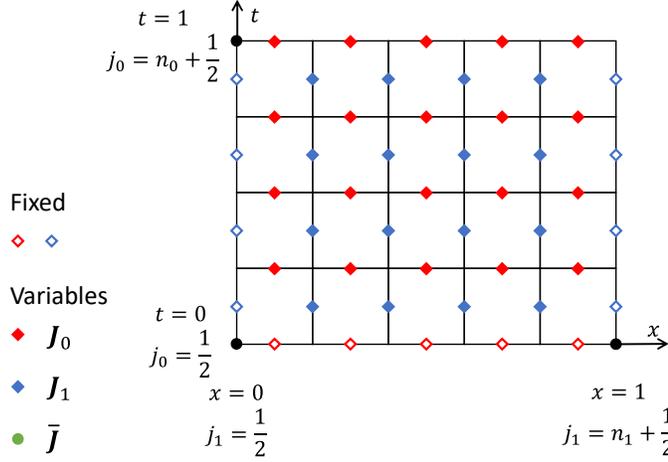}
\vspace{-1cm}
\caption{Illustration of discretization (MFG).}
\label{fig: grid mfg}
\end{figure}

To apply FISTA to MFG, we follow a first-discretize-then-optimize approach. One crucial difference between MFG and MFP is whether $\rho(1,\cdot)$ is provided explicitly. For MFG, we consider a discretization in \cref{fig: grid mfg} and modify our previous notations related to $\rho$. 

The index set and discrete variable are now
\begin{equation*}
\begin{aligned}
    &J_0 := \left\{\frac{3}{2},\frac{5}{2},\cdots,n_0+\frac{1}{2} \right\},\quad
    \Rho:=\{\Rho_{\bmj}\}_{\bmj\in\bmJ_0}\in\bbV_{0}:=\bbR^{n_0\times n_1\times\cdots\times n_D},\\
\end{aligned}
\end{equation*}
and the discrete operators are 
\begin{equation*}
\begin{aligned}
    \mathcal{A}_0:\bbV_{0}\to\overline{\bbV}, &\Rho\mapsto\Rhobar, \\ 
    &\Rhobar_{\bmj} := \begin{cases}
    \displaystyle\half\Rho_{\bmj+\frac{\bme_0}{2}},& j_0=1,\vspace{0.1cm}\\
    \displaystyle\half\left( \Rho_{\bmj+\frac{\bme_0}{2}} + \Rho_{\bmj-\frac{\bme_0}{2}} \right),& j_0=2,3,\cdots,n_0,\vspace{0.1cm}
    \end{cases}\\
    \calD_0:\bbV_{0}\to\overline{\bbV}, &\Rho\mapsto \calD_0(\Rho), \\ 
    &(\calD_0(\Rho))_{\bmj}:=\begin{cases}
    \displaystyle\frac{1}{\Delta_0}\Rho_{\bmj+\frac{\bme_0}{2}},& j_0=1,\vspace{0.1cm}\\
    \displaystyle\frac{1}{\Delta_0}\left( \Rho_{\bmj+\frac{\bme_0}{2}} - \Rho_{\bmj-\frac{\bme_0}{2}} \right),& j_0=2,3,\cdots,n_0.\vspace{0.1cm}
    \end{cases}\\
\end{aligned}
\end{equation*}
Since the boundary condition is only at $t=0$, we modify $\Rhobar_{\calA},\Rhobar_{\calD}\in\overline{\bbV}$ to
\begin{equation*}
\begin{aligned}
    &(\Rhobar_{\calA})_{\bmj}:=\begin{cases}
    \displaystyle \half\rho_0\left(\bmx_{\bmj}\right),& j_0=1,\vspace{0.1cm}\\
    \displaystyle 0,& j_0=2,3,\cdots,n_0,\vspace{0.1cm}
    \end{cases} \\ 
    &(\Rhobar_{\calD})_{\bmj}:=\begin{cases}
    \displaystyle -\frac{1}{\Delta_0}\rho_0\left(\bmx_{\bmj}\right),& j_0=1,\vspace{0.1cm}\\
    \displaystyle 0,& j_0=2,3,\cdots,n_0.\vspace{0.1cm}\\
    \end{cases} \\
\end{aligned}
\end{equation*}
Take model \cref{eqn: mfp eg mfg} in \cref{sec: review} as an example, the discrete problem can be formulated as
\begin{equation}
\begin{aligned}
    \min_{(\Rho,\bmM)\in\calC_{MFG}(\Rhobar_{\calD})}\calE_{MFG}(\Rho,\bmM):=&\Delta_0\sum_{d=0}^D \sum_{j_d=1}^{n_d} 
J_{MFG}\left( (\mathcal{A}_0(\Rho)+\Rhobar_{\calA})_{\bmj}, \mathcal{A}(\bmM)_{\bmj},\bmx_{\bmj} \right)  \\
    & + \lambda_G\sum_{\substack{\bmj\in\bmJ_0,\\ j_0=n_0+\half}} \Rho_{\bmj}G(\bmx_{\bmj})
\end{aligned}
\label{eqn: mfg disct prob}
\end{equation}
where 
\begin{equation}
    J_{MFG}(\beta_0,\bmb,\bmx):=L(\beta_0,\bmb)+\lambda_E\beta_0\log(\beta_0)+\lambda_Q \beta_0 Q(\bmx),
\end{equation}
\begin{equation}
    \calC_{MFG}(\Rhobar_{\calD}):=\left\{ \left(\Rho,\bmM\right): \calD_0(\Rho)+\Rhobar_{\calD}+\Divg(\bmM)=\bmO \right\}.
\label{eqn: mfg disct constraint}
\end{equation}
Since this is an optimization problem with linear constraints, we apply FISTA to it as detailed in \cref{alg: fista mfg}. In the algorithm, $\calA_0^\top, \calA^\top, \calD_0^\top,\Divg^\top$ are conjugate operators of $\calA_0, \calA, \calD_0,\Divg$ in norm $\|\cdot\|_F$. Similar to what we discussed before, one can have:
\begin{equation*}
\left\{\begin{aligned}
\partial_{\Rhobar}\calE_{MFG}(\Rho,\bmM) 
&:= \left\{\Delta_0(J_{MFG})_{\beta_0}\left( (\mathcal{A}_0(\Rho)+\Rhobar_{\calA})_{\bmj}, \mathcal{A}(\bmM)_{\bmj},\bmx_{\bmj} \right)\right\}_{\bmj\in\overline{\bmJ}},\\
\partial_{\overline{\bmM}}\calE_{MFG}(\Rho,\bmM) 
&:= \left\{\Delta_0(J_{MFG})_{\bmb}\left( (\mathcal{A}_0(\Rho)+\Rhobar_{\calA})_{\bmj}, \mathcal{A}(\bmM)_{\bmj},\bmx_{\bmj} \right)\right\}_{\bmj\in\overline{\bmJ}},
\end{aligned}\right.
\end{equation*}
and
\begin{equation*}
\left\{\begin{aligned}
\left(\partial_{\Rho}\calE_{MFG}(\Rho,\bmM)\right)_{\bmj} 
&:= \left(\mathcal{A}_0^\top\left( \partial_{\Rhobar}\calE_{MFG}(\Rho,\bmM)\right)\right)_{\bmj},
&\quad j_0\neq n_0+\half\\
\left(\partial_{\Rho}\calE_{MFG}(\Rho,\bmM)\right)_{\bmj} 
&:= \left(\mathcal{A}_0^\top\left( \partial_{\Rhobar}\calE_{MFG}(\Rho,\bmM)\right)\right)_{\bmj}
    +\lambda_G G(\bmx_{\bmj}),&\\
& &\quad j_0=n_0+\half\\
\partial_{\bmM}\calE_{MFG}(\Rho,\bmM) 
&:= \mathcal{A}^\top\left( \partial_{\overline{\bmM}}\calE_{MFG}(\Rho,\bmM) \right).&
\end{aligned}\right.
\end{equation*}

\begin{algorithm}
\caption{FISTA for MFG}
\begin{algorithmic}
\STATE{Parameters}
$\rho_0,\rho_1$
\STATE{Initialization} $\tau^{(1)}=1$, 
$\Rho^{(0)}=\widehat{\Rho}^{(0)}=\mathbf{1}$, and 
$M_d^{(0)}=\widehat{M_d}^{(0)}=\mathbf{1}$.
\FOR{$k=0,1,2,\ldots$}
\STATE{\textbf{gradient descent:}} 
\begin{equation*}
\left\{\begin{aligned}
    &\Rho^\kp = \Rho^\upk-\eta^\upk\partial_\Rho\calE_{MFG}(\Rho^\upk,\bmM^\upk),\\
    &\bmM^\kp = \bmM^\upk-\eta^\upk \partial_{\bmM} \calE_{MFG}(\Rho^\upk,\bmM^\upk)
\end{aligned}\right.
\end{equation*}
\STATE{\textbf{projection:}} solve $\Phibar^\kpp$ for 
\begin{equation}
\begin{aligned}
\left(\calD_0\calD_0^\top+\Divg\Divg^\top\right)\Phibar^\kpp = \calD_0\left(\Rho^\kp\right)+\Rhobar_{\calD}+\Divg\left(\bmM^\kp\right),
\end{aligned}
\label{eqn: proj phi mfg}
\end{equation}
and project $\left(\Rho^\kp,\bmM^\kp\right)$ to $\calC_{MFG}(\Rhobar_{\calD})$ by
\begin{equation*}
\left\{\begin{aligned}
    &\Rho^\kpp=\Rho^\kp-\calD_0^\top\left(\Phibar^\kpp\right),\\
    &\bmM^\kpp=\bmM^\kp-\Divg^\top\left(\Phibar^\kpp\right).
\end{aligned}\right.
\label{eqn: proj update mfg}
\end{equation*}
\STATE{\textbf{update}} 
\begin{equation*}
\begin{aligned}
&\tau^\kpp=\frac{1+\sqrt{1+4\left(\tau^\upk\right)^2}}{2},\\
&\omega^\upk=\frac{\tau^\upk-1}{\tau^\kpp},\\
&\left(\widehat{\Rho}^\kpp, \widehat{\bmM}^\kpp \right) = \left(1+\omega^\upk\right) \left( \Rho^\kpp,\bmM^\kpp \right) -\omega^\upk\left( \Rho^\upk,\bmM^\upk \right).
\end{aligned}
\end{equation*}
\ENDFOR
\end{algorithmic}
\label{alg: fista mfg}
\end{algorithm}

\begin{remark}
\label{rem: disct lap^-1 mfg}
In \cref{alg: fista mfg}, we also need to solve a discrete Poisson equation \cref{eqn: proj phi mfg} and the approach is similar as presented in \cref{rem: disct lap^-1}.
If we write
$\Lap=\calD_0\calD_0^\top+\Divg\Divg^\top$,
where $\calD_0$ is the modified definition in this section. Then the linear decomposition of $\Lap$ is
\begin{equation*}
    \Lap(\Phibar)=\sum_{\bmi\in\overline{\bmJ}}\lambda^{\bmi}\left\langle \Phibar,\Psi^{\bmi} \right\rangle\Psi^{\bmi},
\end{equation*}
where $\lambda^{\bmi}\in\bbR,\Psi^{\bmi}\in\overline{\bbV}$ are
\begin{equation*}
\begin{aligned}
    \lambda^{\bmi} &= -4n_0^2\sin^2\left(\frac{(i_0-1/2)\pi}{2(n_0+1/2)}\right) 
     -4\sum_{d=1}^D n_d^2 \sin^2\left(\frac{(i_d-1)\pi}{2n_d}\right),\\
    \Psi^{\bmi}_{\bmj} &=  \sqrt{\frac{4}{2n_0+1}}\cos\left( \left(j_0-\half\right) \frac{(i_0-1/2)\pi}{n_0+1/2} \right) \\
    &\quad\times\prod_{d=1}^D\sqrt{\frac{1+\delta_{1i_d}}{n_d}} \cos\left( \left(j_d-\half\right) \frac{(i_d-1)\pi}{n_d} \right),
\end{aligned}
\end{equation*}
Note that $\lambda_{\bmi}<0$ for any $\bmi\in\overline{\bmJ}$, we can define $\Lap^{-1}(\Phibar)$ by 
\begin{equation}
    \Lap^{-1}(\Phibar):=\sum_{\bmi\in\overline{\bmJ}}\frac{1}{\lambda^{\bmi}} \left\langle \Phibar,\Psi^{\bmi} \right\rangle \Psi^{\bmi}.
\end{equation}

\end{remark}
\subsection{Multilevel and Multigrid FISTA}
\label{sec: ml and mg}

Inspired by \cite{borzi2009multigrid,liu2021multilevel}, we borrow ideas from multigrid and multilevel methods in numerical PDEs to our variational problem. These methods can reduce computational cost on the finest level and thus accelerate the proposed algorithm. The implementation details are presented in this section. 

For notation simplicity, we assume $h=\Delta_0=\Delta_1=\cdots=\Delta_D$ in this section. 
Let $\prescript{h}{}{\Omega}$ be a grid with $h=\Delta_d$, $\prescript{h}{}{\bmJ_d}$ be the certain $\bmJ_d$ on the grid. Then index $\prescript{h}{}{\bmj}\in\prescript{h}{}{\bmJ_d}$ stands for the point $h\bmj.$ If there is no ambiguity, we can omit the prescript of $\bmj$. For example, we define $\prescript{h}{}{u}_{\bmj}=u\left(h(\bmj-\frac{\bm{1}}{2})\right)$ for any function $u$ and approximate the value by $\prescript{h}{}{U}_{\bmj}.$

\begin{figure}[ht]
\centering
\vspace{-0.8cm}
\subfigure{
\includegraphics[width=6cm]{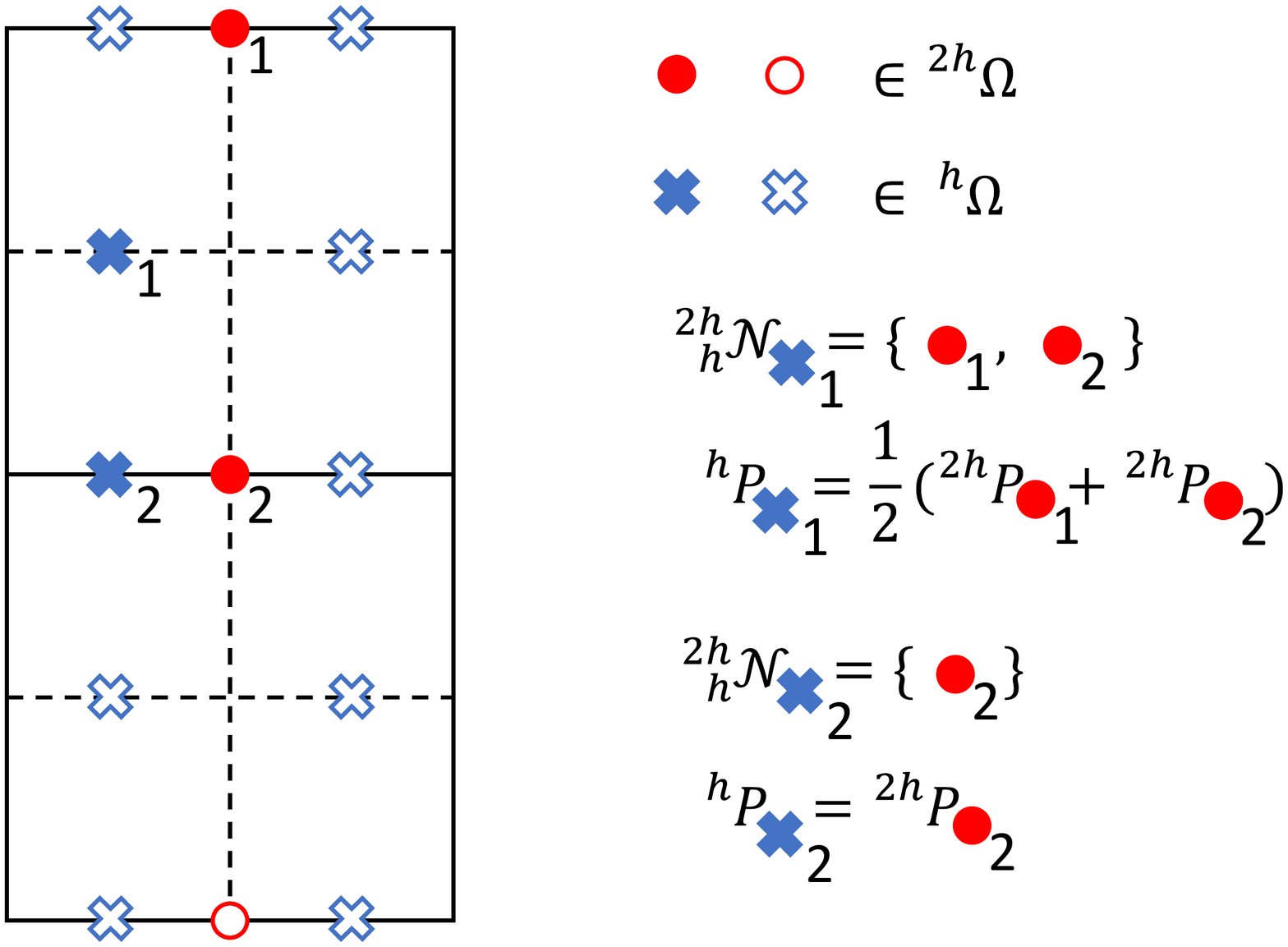}
}
\subfigure{
\includegraphics[width=6cm]{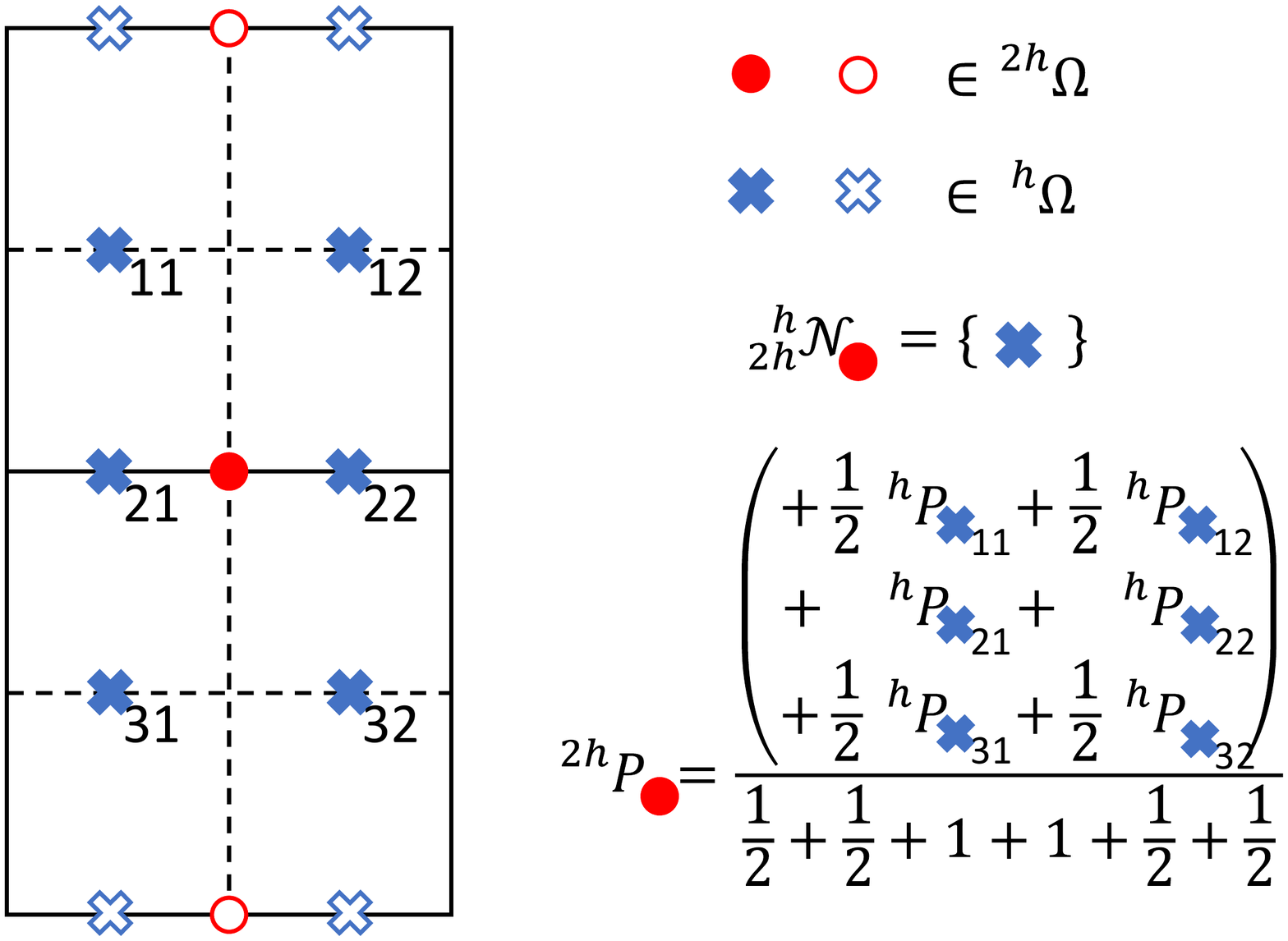}
}
\vspace{-0.8cm}
\caption{Illustration of Prolongation (left) and Restriction (right) for 1D case.}
\label{fig: res and pro}
\end{figure}

Consider $L$ levels of grids $\prescript{h_1}{}{\Omega},\ldots,\prescript{h_L}{}{\Omega}$ where the finest level is $\prescript{h_1}{}{\Omega}$, and $h_l:=2^{l-1}h_1$. 
We first define how to prolongate values on a coarser grid into a finer grid. 
Assume that $\prescript{h}{}{\bmj}\in\prescript{h}{}{\bmJ_d}$ stands for point $h\left(\bmj-\frac{\bm{1}}{2}\right)$ on the finer grid $\prescript{h}{}{\Omega}$, we define its neighbourhood on the coarser grid $\prescript{2h}{}{\Omega}$ as
\begin{equation}
\begin{aligned}
    \prescript{2h}{h}{\mathcal{N}}_{\bmj}:=
    \Bigg\{ \prescript{2h}{}{\bmi}\in\prescript{2h}{}{\bmJ_d}:
    &\left\|2h\left(\prescript{2h}{}{\bmi}-\frac{\bm{1}}{2}\right)-h\left(\prescript{h}{}{\bmj}-\frac{\bm{1}}{2}\right) \right\|_2  \\
    &= \min_{\prescript{2h}{}{\bmi'}\in\prescript{2h}{}{\bmJ_d}}\left\|2h\left(\prescript{2h}{}{\bmi'}-\frac{\bm{1}}{2}\right)-h\left(\prescript{h}{}{\bmj}-\frac{\bm{1}}{2}\right) \right\|_2\Bigg\}.
\end{aligned}
\end{equation}
Then with boundary values
\begin{equation*}
\begin{aligned}
    &\prescript{2h}{}{\Rho_{\bmi}}=\prescript{2h}{}{(\rho_0)}_{\bmi},\quad i_0=\half,\\
    &\prescript{2h}{}{\Rho_{\bmi}}=\prescript{2h}{}{(\rho_1)}_{\bmi},\quad i_0=\frac{1}{2h}+\half,\\
    &\prescript{2h}{}{(M_d)_{\bmi}}=0,\quad i_d=\half,\frac{1}{2h}+\half,\\
\end{aligned}
\end{equation*}
we define the prolongation $(\prescript{h}{}{\Rho},\prescript{h}{}{\bmM})=\Pro(\prescript{2h}{}{\Rho},\prescript{2h}{}{\bmM})$ by averaging values in neighbourhoods:
\begin{equation}
\left\{\begin{aligned}
    &\prescript{h}{}{\Rho}_{\bmj}:= \frac{1}{\left| \prescript{2h}{h}{\mathcal{N}}_{\bmj} \right|} \sum_{\prescript{2h}{}{\bmi}\in\prescript{2h}{h}{\mathcal{N}}_{\bmj}} \prescript{2h}{}{\Rho}_{\bmi}, &\forall\prescript{h}{}{\bmj}\in\prescript{h}{}{\bmJ_0},\\
    &\prescript{h}{}{(M_d)}_{\bmj}:= \frac{1}{\left| \prescript{2h}{h}{\mathcal{N}}_{\bmj} \right|} \sum_{\prescript{2h}{}{\bmi}\in\prescript{2h}{h}{\mathcal{N}}_{\bmj}} \prescript{2h}{}{(M_d)}_{\bmi}&\forall\prescript{h}{}{\bmj}\in\prescript{h}{}{\bmJ_d}.
\end{aligned}\right.
\label{eqn: prolongation}
\end{equation}
An example of prolongation in 1D is shown in the left panel of \cref{fig: res and pro}.

From a finer grid to a coarser grid, the neighbourhood is defined inversely. Suppose $\prescript{2h}{}{\bmi}\in\prescript{2h}{}{\bmJ_d}$, its neighbourhood is the set of all $\prescript{h}{}{\bmj}\in\prescript{h}{}{\bmJ_d}$ whose neighbourhood includes $\prescript{2h}{}{\bmi}$: 
\begin{equation}
    \prescript{h}{2h}{\mathcal{N}}_{\bmi}:=\left\{ \prescript{h}{}{\bmj}\in\prescript{h}{}{\bmJ_d}:\prescript{2h}{}{\bmi}\in\prescript{2h}{h}{\mathcal{N}}_{\bmj}\right\}.
\end{equation}
and the restriction from finer level to coarser level $(\prescript{2h}{}{\Rho},\prescript{2h}{}{M})=\Res(\prescript{h}{}{\Rho},\prescript{h}{}{M})$ is defined by a weighted average over neighbourhoods: 
\begin{equation}
\left\{\begin{aligned}
    &\prescript{2h}{}{\Rho}_{\bmi}:=\left.
    \sum_{\bmj\in\prescript{h}{2h}{\mathcal{N}}_{\bmi}}\frac{1}{\left| \prescript{2h}{h}{\mathcal{N}}_{\bmj} \right|}  \prescript{h}{}{\Rho}_{\bmj}
    \middle/
    \sum_{\bmj\in\prescript{h}{2h}{\mathcal{N}}_{\bmi}}\frac{1}{\left| \prescript{2h}{h}{\mathcal{N}}_{\bmj} \right|} \right.,&\forall\prescript{2h}{}{\bmi}\in\prescript{2h}{}{\bmJ_0},\\
    &\prescript{2h}{}{(M_d)}_{\bmi}:=\left. 
    \sum_{\bmj\in\prescript{h}{2h}{\mathcal{N}}_{\bmi}}\frac{1}{\left| \prescript{2h}{h}{\mathcal{N}}_{\bmj} \right|}  \prescript{h}{}{(M_d)}_{\bmj}
    \middle/
    \sum_{\bmj\in\prescript{h}{2h}{\mathcal{N}}_{\bmi}}\frac{1}{\left| \prescript{2h}{h}{\mathcal{N}}_{\bmj} \right|}\right., &\forall\prescript{2h}{}{\bmi}\in\prescript{2h}{}{\bmJ_d}.
\end{aligned}\right.
\label{eqn: restriction}
\end{equation}
An example of restriction is shown in the right panel of \cref{fig: res and pro}.

We describe our multigrid FISTA in \cref{alg: mg FISTA}, in which $\Solve_K(\cdot)$ means run \cref{alg: fista disct} for $K$ iterations and $\Solve(\cdot)$ means run the algorithm till convergence. The first two inputs of $\Solve(\cdot)$ are initial and terminal densities $\rho_0,\rho_1$, and the last two inputs are initialization $\Rho^{(0)}=\widehat{\Rho}^{(0)},\ \bmM^{(0)}=\widehat{\bmM}^{(0)}$. 

\begin{algorithm}
\caption{Multigrid FISTA for MFP}
\begin{algorithmic}
\STATE{Parameters}
$L, h_l=2^{l-1}h, K, \prescript{h_l}{}{\rho_0},\prescript{h_l}{}{\rho_1}   (l=1,\ldots,L)$

\STATE{Initialization} 
$\prescript{h_1}{}{\Rho^{(0)}}=\mathbf{1},\ \prescript{h_1}{}{\bmM^{(0)}}=\mathbf{1}$
\STATE

\STATE{\textbf{pre-smoothing}}
$$\left( \prescript{h_1}{}{\Rho},\prescript{h_1}{}{\bmM} \right) = \Solve_{K}\left( 
\rho_0,\rho_1,
\prescript{h_1}{}{\Rho^{(0)}},\prescript{h_1}{}{\bmM^{(0)}}
\right)$$

\FOR{$l=2,3,\ldots,L$}
\STATE 
$$\left( \prescript{h_l}{}{\Rho},\prescript{h_l}{}{\bmM} \right) = \Solve_{K}\left( 
\rho_0,\rho_1,
\Res\left(\prescript{h_{l-1}}{}{\Rho},\prescript{h_{l-1}}{}{\bmM}\right)
\right)$$
\ENDFOR

\STATE

\STATE{\textbf{correction and post-smoothing}}
$$\left( \prescript{h_L}{}{\Rho},\prescript{h_L}{}{\bmM} \right) = \Solve\left( 
\rho_0,\rho_1,
\prescript{h_L}{}{\Rho},\prescript{h_L}{}{\bmM}
\right)$$

\FOR{$l=L-1,L-2,\ldots,1$}
\STATE
\begin{equation*}
\begin{aligned}
    \left( \prescript{h_l}{}{\Rho},\prescript{h_l}{}{\bmM} \right) = \left( \prescript{h_l}{}{\Rho},\prescript{h_l}{}{\bmM} \right)
    &+ \Solve\left( \rho_0,\rho_1,
    \Pro\left(\prescript{h_{l+1}}{}{\Rho},\prescript{h_{l+1}}{}{\bmM}\right)\right) \\
    &- \Pro\left(\prescript{h_{l+1}}{}{\Rho},\prescript{h_{l+1}}{}{\bmM}\right)
\end{aligned}
\end{equation*}
\ENDFOR

\end{algorithmic}
\label{alg: mg FISTA}
\end{algorithm}

Note that to keep the cost of \cref{alg: mg FISTA} low, we need to choose a $K$ not very large. Motivated by \cite{liu2021multilevel}, we can remove the pre-smoothing steps by setting $K=0$ and this leads to our \cref{alg: ml FISTA}: multilevel FISTA.

\begin{algorithm}
\caption{Multilevel FISTA for MFP}
\begin{algorithmic}
\STATE{Parameters}
$L, h_l=2^{l-1}h, \prescript{h_l}{}{\rho_0},\prescript{h_l}{}{\rho_1}   (l=1,\ldots,L)$

\STATE{Initialization} 
$\prescript{h_L}{}{\Rho^{(0)}}=\mathbf{1},\ \prescript{h_L}{}{\bmM^{(0)}}=\mathbf{1}$

\STATE $\left( \prescript{h_L}{}{\Rho},\prescript{h_L}{}{\bmM} \right) = \Solve\left( 
\rho_0,\rho_1,
\prescript{h_L}{}{\Rho^{(0)}},\prescript{h_L}{}{\bmM^{(0)}}
\right)$

\FOR{$l=L-1,L-2,\ldots,1$}
\STATE 
$\left( \prescript{h_l}{}{\Rho},\prescript{h_l}{}{\bmM} \right) = \Solve\left( 
\rho_0,\rho_1,
\Pro\left(\prescript{h_{l+1}}{}{\Rho},\prescript{h_{l+1}}{}{\bmM}\right)
\right)$
\ENDFOR

\end{algorithmic}
\label{alg: ml FISTA}
\end{algorithm}

\section{Numerical Experiments}
\label{sec: num eg}

In this section, we conduct comprehensive experiments to show the efficiency and effectiveness of the proposed numerical algorithms. We first numerical verify the convergence of rate of the algorithm related to the mesh size. After that, our computation efficiency tests demonstrate that the proposed \cref{alg: fista disct} has comparable efficiency with the state-of-the-art methods. Interestingly, the proposed multilevel method performs around 10 times faster than existing methods. We further illustrate the flexibility of our algorithms on different MFP problems. In all the numerical experiments, we use the dynamic cost function $L$ defined in \cref{eqn: dynamical term L}. All of our numerical experiments are implemented in Matlab on a PC with an Intel(R) i7-8550U 1.80GHz CPU  and 16 GB memory. 

\subsection{Convergence Rate}
\label{sec: num eg convergence}
To numerically verify the theoretical convergence analysis discussed in \cref{sec: convergence},  we first apply the proposed numerical algorithm to a simple 1D OT example with exact solution as follows. 

Let $\Omega=[0,1]$, $\rho_0(x) = x+\half, \rho_1(x) = 1.$ Then we can have the following theoretical solution of the OT between $\rho_0$ and $\rho_1$.
\begin{equation}
\rho^*(t,x) =
\left\{\begin{aligned}
    \displaystyle &x+\half,&\quad t=0,\\
    \displaystyle &\frac{ \sqrt{2tx+\left(\frac{t}{2}-1\right)^2} + t-1 }{ t \sqrt{2tx+\left(\frac{t}{2}-1\right)^2} },&\quad 0<t\leq 1.
\end{aligned}\right.
\end{equation}
\begin{equation}
 m^*(t,x) = 
 \left\{\begin{aligned}
    \displaystyle &\frac{1}{4}x(x-1)(2x+1), &\quad t=0,\\
    \displaystyle &\frac{x}{t^2} - \frac{3-t}{2t^3}\sqrt{2tx+\left(\frac{t}{2}-1\right)^2} \\
    & - \frac{(t-1)(t^2-4)}{8t^3}\frac{1}{\sqrt{2tx+\left(\frac{t}{2}-1\right)^2}} - \frac{3t-4}{2t^3}, &\quad 0<t\leq 1.
\end{aligned}\right.
\end{equation}
We also know $\displaystyle W_2^2(\rho_0,\rho_1) = \frac{1}{120}$.

\begin{table}[htbp]
\centering
\caption{Convergence rate of \cref{alg: fista disct} applied to 1D OT problem ($k=50000$).}
\begin{tabular}{ll|ll|ll|ll}
\toprule
$\Delta_0$  & $\Delta_1$     & $\left\| \bmE^{(k,*)} \right\|_2$  & order  & $\left\| \bmE^{(k,*)} \right\|_\infty$ & order  & $W_2^2$ error     & order  \\ \midrule
1/16  & 1/64    &  3.19E-04 &      & 2.88E-03 &      & 4.88E-06 &      \\
1/32  & 1/128   &  1.08E-04 & 1.56 & 1.47E-03 & 0.97 & 1.22E-06 & 2.00 \\
1/64  & 1/256   &  3.76E-05 & 1.53 & 7.44E-04 & 0.98 & 3.05E-07 & 2.00 \\
1/128 & 1/512   &  1.37E-05 & 1.46 & 3.62E-04 & 1.04 & 7.63E-08 & 2.00 \\
\bottomrule
\end{tabular}
\label{tab: convergence rate}
\end{table}

Note that it would be quite difficult to check $\bmE^\upk$ as we do not have evolution path, $\rho^\upk$ and $\bmm^\upk$, in the continuous \cref{alg: fista cts}. Instead, we compute the following values:
\begin{equation*}
\begin{aligned}
    &\left\| \bmE^{(k,*)} \right\|_2:=
    \sqrt{\Delta_0\Delta_1}\left[\sum_{\bmj\in\bmJ_0} \left|\Rho^\upk_{\bmj}-\rho^*_{\bmj}\right|^2+ \sum_{\bmj\in\bmJ_1}\left|M^\upk_{\bmj}-m^*_{\bmj} \right|^2  \right]^{\half},\\
    &\left\| \bmE^{(k,*)} \right\|_\infty:=
    \max\left\{\max_{\bmj\in\bmJ_0} \left|\Rho^\upk_{\bmj}-\rho^*_{\bmj}\right|, \max_{\bmj\in\bmJ_1}\left|M^\upk_{\bmj}-m^*_{\bmj} \right|  \right\},\\
    &\text{$W_2^2$ error: } \left| \Delta_0\Delta_1\calE\left(\Rho^\upk,\bmM^\upk\right)-W_2^2(\rho_0,\rho_1) \right|.
\end{aligned}
\end{equation*}
Here $\left\| \bmE^{(k,*)} \right\|_2$ is related to $\left\| \bmE^\upk \right\|_2$ by:
\begin{equation*}
\begin{aligned}
    &\left\| \bmE^{(k,*)} \right\|_2 \leq \left\| \bmE^\upk \right\|_2 
    +\left\| \left(\rho^\upk_{\bmJ_0},\bmm^\upk_{\bmJ}\right)-\left(\rho^*_{\bmJ_0},\bmm^*_{\bmJ}\right) \right\|_2\\
\end{aligned}
\end{equation*}
For given $\Delta_0,\Delta_1$, we can choose very large $k$ such that
\begin{equation*}
\begin{aligned}
    &\left\| \bmE^{(k,*)} \right\|_2 \leq \left\| \bmE^\upk \right\|_2 
    +\epsilon^\upk\\
\end{aligned}
\end{equation*}
and $\epsilon^\upk \ll \Delta_0+\Delta_1$. 
Fixing $k$, according to our theoretical analysis, we expect to observe at least
\begin{equation*}
\begin{aligned}
    &\left\| \bmE^{(k,*)} \right\|_2 =\bigO\left(\Delta_0+\Delta_1\right).
\end{aligned}
\end{equation*}
and
\begin{equation*}
\begin{aligned}
    &\left\| \bmE^{(k,*)} \right\|_\infty \leq \left\| \bmE^{(k,*)} \right\|_F = (\Delta_0\Delta_1)^{-\half}\left\| \bmE^{(k,*)} \right\|_2 =\bigO\left(1\right). 
\end{aligned}
\end{equation*}
Numerical results are shown in \cref{tab: convergence rate} where we observe
\begin{equation*}
\begin{aligned}
    \left\| \bmE^{(k,*)} \right\|_2 =\bigO\left(\Delta_0^{1.5}+\Delta_1^{1.5}\right),
\quad\left\| \bmE^{(k,*)} \right\|_\infty =\bigO\left(\Delta_0+\Delta_1\right).
\end{aligned}
\end{equation*}
This indicates that convergence rate of our numerical experiments perform better than theoretical prediction. This is not surprise as the way of our theoretical analysis may not be sharp.

\subsection{Computation Efficiency}
\label{sec: num eg speed comp}
\begin{table}[htbp]
\centering
\caption{Time comparison of OT in 1D for different grid sizes ($n_0=64, tol=10^{-4}$, F=FISTA, A=ALG, G=G-prox) with the best performance highlighted in red. }
\begin{tabular}{c|lll|lll|lll}
\toprule
          & \multicolumn{3}{c|}{Iter}                        & \multicolumn{3}{c|}{Time (s)}                 & \multicolumn{3}{c}{Time(s)/Iter}                      \\ \cline{2-10} 
  $\nx$       & F & A                        & G                     & F                        & A   & G & F                           & A      & G   \\ \midrule
256  & 611 & 435 & \color[HTML]{FE0000}426 & \color[HTML]{FE0000}1.74  & 2.99  & 2.93  & \color[HTML]{FE0000}2.85E-03 & 6.86E-03 & 6.88E-03 \\
512  & 611 & 435 & \color[HTML]{FE0000}429 & \color[HTML]{FE0000}3.06  & 5.91  & 6.92  & \color[HTML]{FE0000}5.00E-03 & 1.36E-02 & 1.61E-02 \\
1024 & 611 & 435 & \color[HTML]{FE0000}430 & \color[HTML]{FE0000}7.60  & 12.93  & 12.85  & \color[HTML]{FE0000}1.24E-02 & 2.97E-02 & 2.99E-02 \\
2048 & 611 & 435 & \color[HTML]{FE0000}431 & \color[HTML]{FE0000}24.84 & 36.15 & 32.72 & \color[HTML]{FE0000}4.07E-02 & 8.31E-02 & 7.59E-02 \\
4096 & 611 & 435 & \color[HTML]{FE0000}431 & \color[HTML]{FE0000}51.79 & 69.99 & 68.09 & \color[HTML]{FE0000}8.48E-02 & 1.61E-01 & 1.58E-01
\\ \bottomrule
\end{tabular}
\label{tab: 1dot gaussian diffsize}
\end{table}

\begin{table}[htbp]
\centering
\caption{Time comparison of OT in 2D for different grid sizes ($n_0=64$, F=FISTA, A=ALG, G=G-prox) with the best performance highlighted in red.}
\resizebox{5.1in}{!}{
\begin{tabular}{c|lll|lll|lll}
\toprule
& \multicolumn{3}{c|}{Iter}  & \multicolumn{3}{c|}{Time (s)}    & \multicolumn{3}{c}{Time(s)/Iter}  \\ \cline{2-10} 
$n_1,n_2$                                & F                      & A                       & G                    & F                         & A                           & G                        & F                           & A                             & G                         \\ \midrule
128 & 116 & \color[HTML]{FE0000}64 & 66 & 46.20   & \color[HTML]{FE0000}46.07  & 46.49  & \color[HTML]{FE0000}3.98E-01 & 7.20E-01 & 7.92E-01 \\
256 & 116  & \color[HTML]{FE0000}64 & 66 & 212.31 & 201.52 & \color[HTML]{FE0000}190.31 & \color[HTML]{FE0000}1.83E+00 & 3.15E+00 & 3.27E+00 \\
512 & 116  & \color[HTML]{FE0000}64 & 66 & 810.86 & 761.65 & \color[HTML]{FE0000}752.59  & \color[HTML]{FE0000}6.99E+00 & 1.19E+01 & 1.14E+01 
\\\bottomrule
\end{tabular}}
\label{tab: 2dot gaussian diffsize}
\end{table}

In this part, we would like to demonstrate the efficiency of our algorithms by comparing with state-of-the-art methods for dynamic OT problems. We apply our algorithms to OT problems with $\rho_0,\rho_1$ being Gaussian distribution densities, and compare the results and computation time with those using ALG(augmented Lagrangian)~\cite{benamou1999numerical,benamou2000computational} and G-prox~\cite{jacobs2019solving}.
For all approaches, the stopping criteria are 
\begin{equation*}
    \left\|\left(\Rho^\kpp,\bmM^\kpp\right) - \left(\Rho^\upk,\bmM^\upk\right)\right\|_2\leq tol.
\end{equation*}

In \cref{tab: 1dot gaussian diffsize} and \cref{tab: 2dot gaussian diffsize}, we report computation time and number of iterations for each algorithms on different grid sizes in 1D and 2D. 
From the tables, the proposed \cref{alg: fista disct} outperforms ALG and G-prox in 1D and achieves similar efficiency in 2D.
Interestingly, CPU time per iteration in our algorithm is the least among these three algorithms.
This is because, at each iteration, solving a Poisson equation is required for all three algorithms while our method does not need to solve $\prod_{d=0}^D n_d$ cubic equations required in ALG and G-prox.
Therefore our method needs less time in 1D experiment although it needs more iterations to achieve the given stopping criteria. 
While, this computation save is marginal comparing with the cost of solving Poisson equation in 2D. Thus,  our method spend comparable time instead of less time in this 2D experiment. 

\begin{table}[htbp]
\centering
\caption{Efficiency and accuracy comparisons of OT in 1D ($\nt=64, \nx=256$) with the best performance highlighted in red. }
\begin{tabular}{l|lllll}
\toprule
            & \makecell[l]{Num\\Iter} & Time (s)                & \makecell[l]{Stationarity\\ Residue} &\makecell[l]{Feasibility \\Residue} &\makecell[l]{Mass \\Residue}\\ \midrule
FISTA          & 611  & 1.723 & \color[HTML]{FE0000}3.27E-05 &\color[HTML]{FE0000}2.28E-13   & 1.33E-15         \\
ALG            & 435  & 2.840 &9.43E-05& 2.41E-04 &1.64E-08\\
G-prox         & 426  & 2.761 &1.93E-04& 1.88E-04 &2.96E-08\\
MLFISTA        & 882 & \color[HTML]{FE0000}0.422 & 7.97E-05  &\color[HTML]{FE0000}2.28E-13& \color[HTML]{FE0000}1.11E-15         \\
MGFISTA $(K=5)$  & 1448 & 1.195 & 4.79E-05 &2.33E-13&1.77E-15          \\
MGFISTA $(K=10)$ & 1517 & 1.341 & 3.95E-05 &\color[HTML]{FE0000}2.28E-13&2.22E-15         
\\ \bottomrule
\end{tabular}
\label{tab: 1dot gaussian samesize}
\end{table}

\begin{table}[htbp]
\centering
\caption{Efficiency and accuracy comparison of OT in 2D ($\nt=64, \nx=n_2=256$) with the best performance highlighted in red. }
\begin{tabular}{l|lllll}
\toprule
            & \makecell[l]{Num\\Iter} & Time (s)                & \makecell[l]{Stationarity\\ Residue} &\makecell[l]{Feasibility \\Residue} &\makecell[l]{Mass \\Residue}\\ \midrule
FISTA          & 116  & 232.560 & 9.22E-04&6.01E-13  &\color[HTML]{FE0000}1.42E-14           \\
ALG            & 64  & 211.043 &\color[HTML]{FE0000}8.75E-04& 4.99E-03 & 3.10E-03\\
G-prox         & 66  & 208.696 &9.29E-04& 6.88E-03 & 3.10E-03\\
MLFISTA        & 162 & \color[HTML]{FE0000}12.853  & 3.43E-03&\color[HTML]{FE0000}2.95E-13  &2.07E-14            \\
MGFISTA $(K=5)$  & 315 & 134.226 & 1.07E-03 &5.99E-13& 1.67E-14            \\
MGFISTA $(K=10)$ & 315 & 170.580 & 9.86E-04 &6.00E-13& 2.02E-14          
\\ \bottomrule
\end{tabular}
\label{tab: 2dot gaussian samesize}
\end{table}

\begin{figure}[htbp]
\centering
\includegraphics[width=4cm]{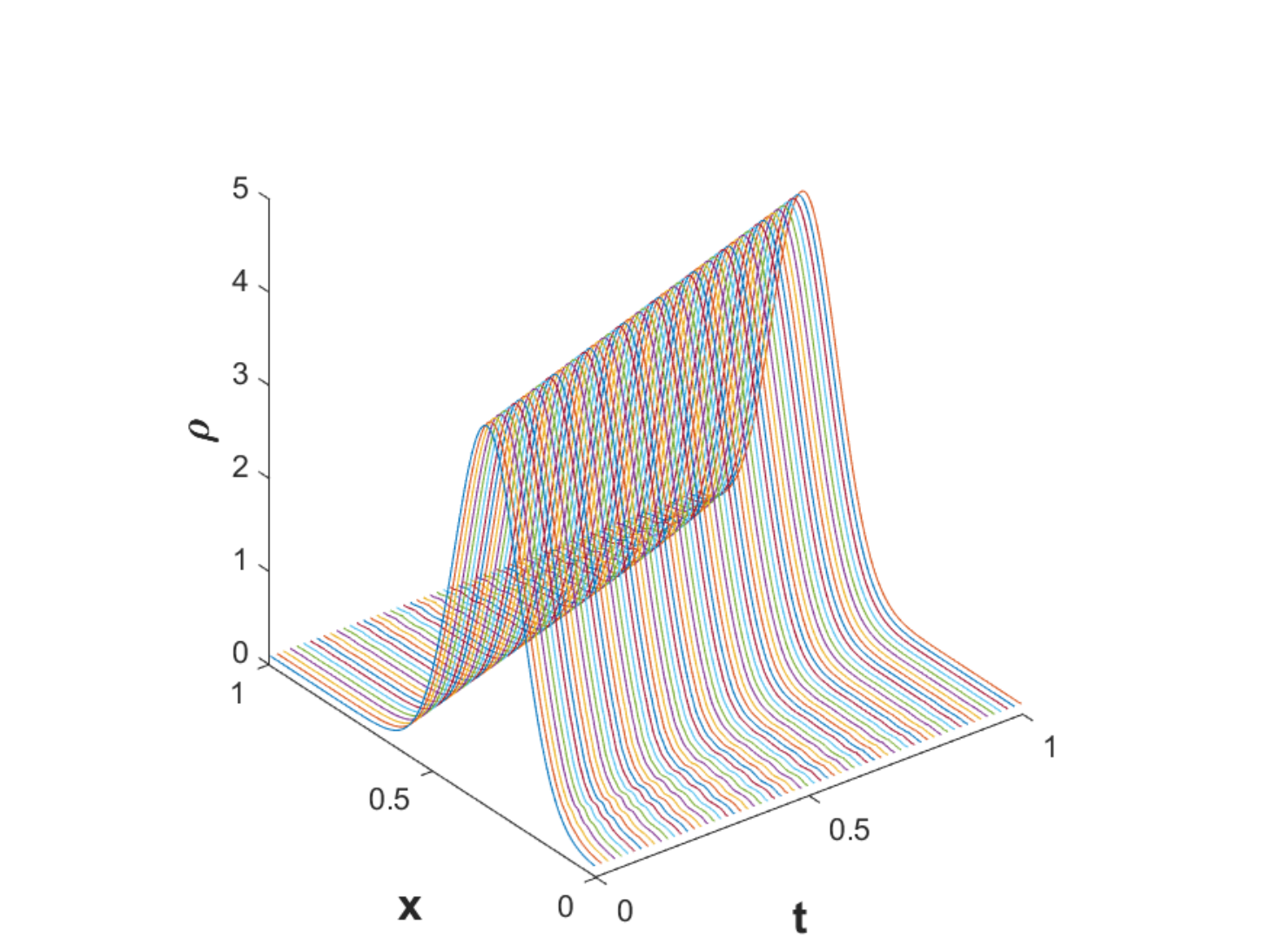}
\includegraphics[width=4cm]{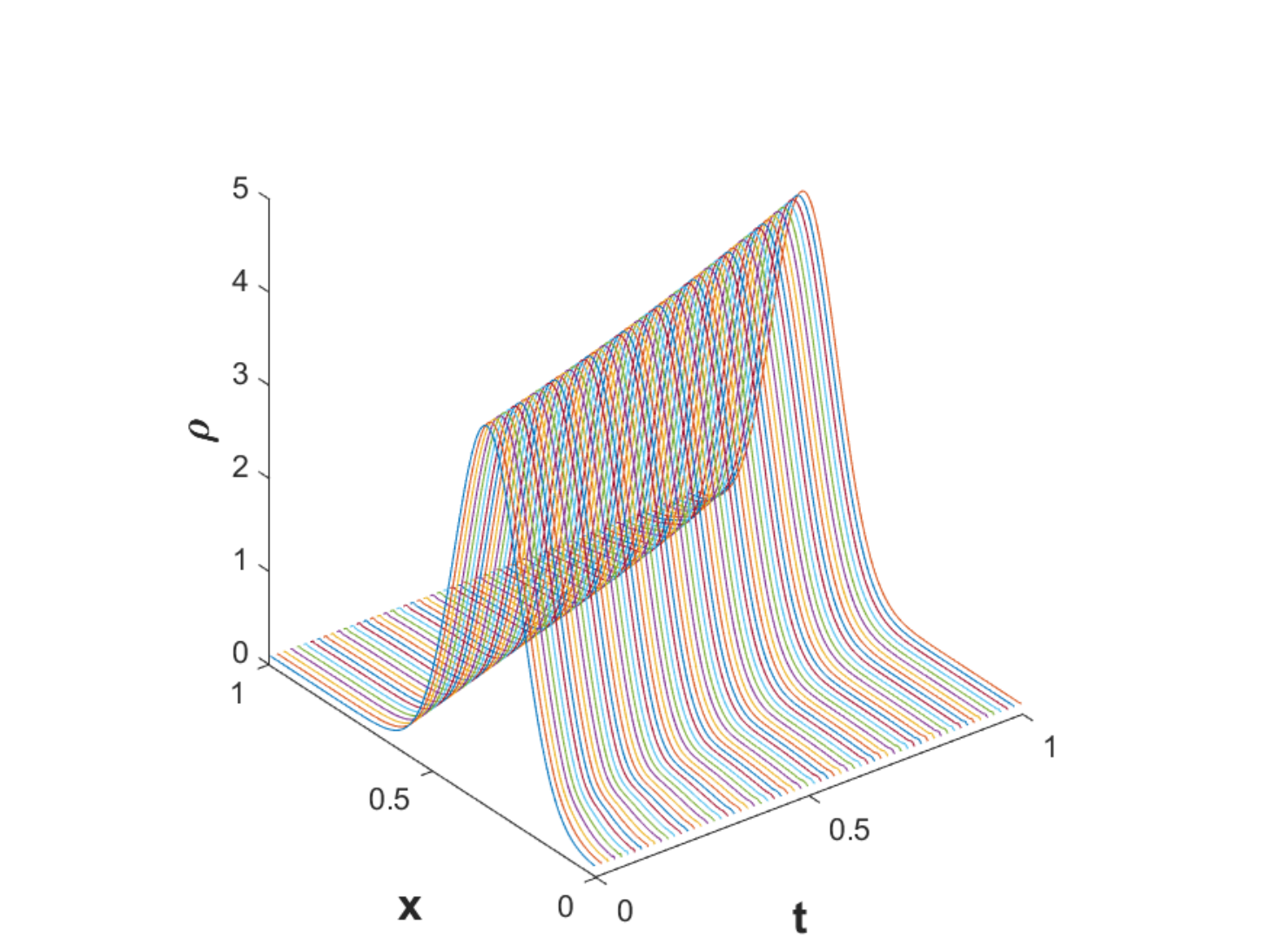}
\includegraphics[width=4cm]{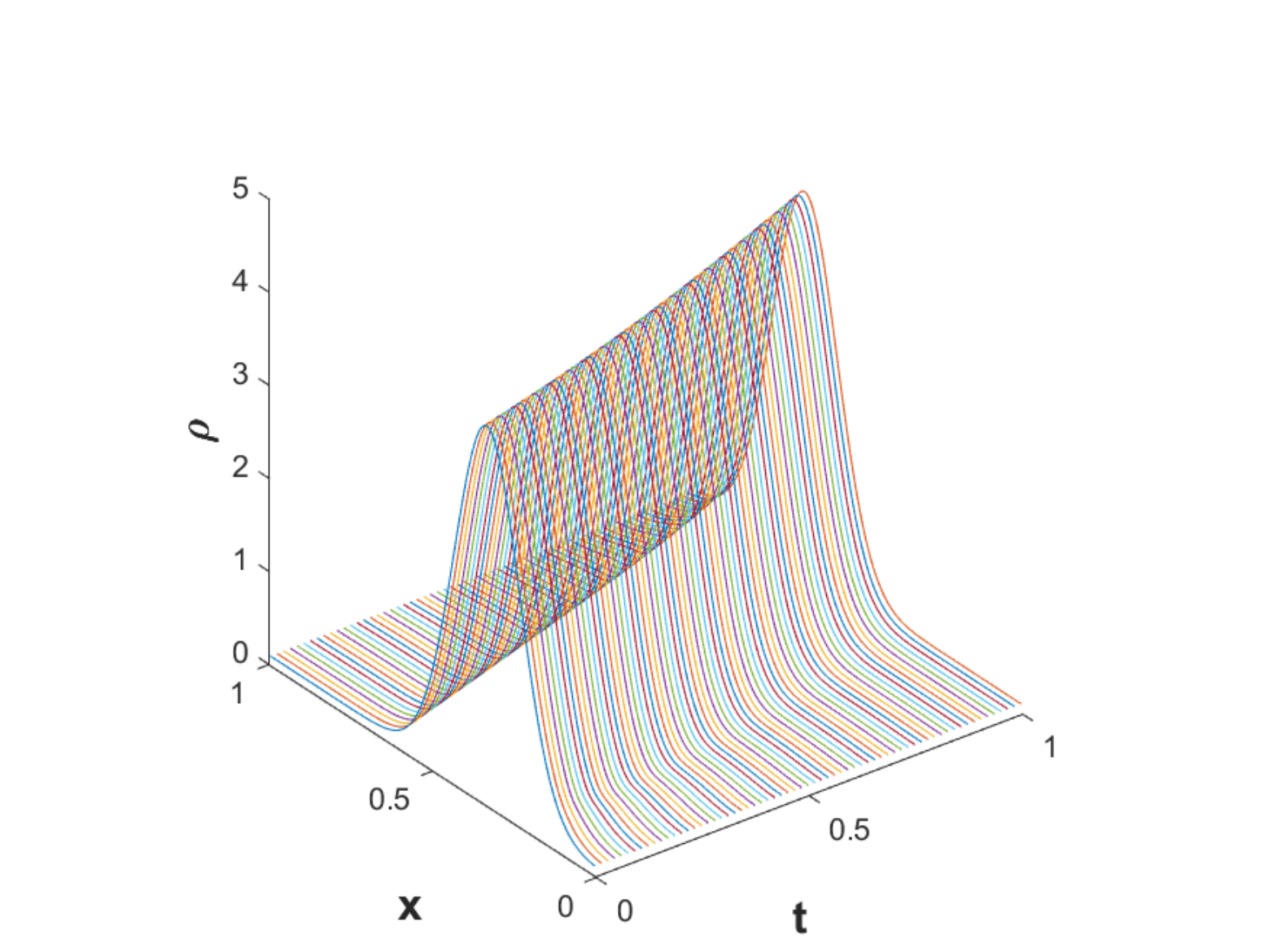}\\
\includegraphics[width=4cm]{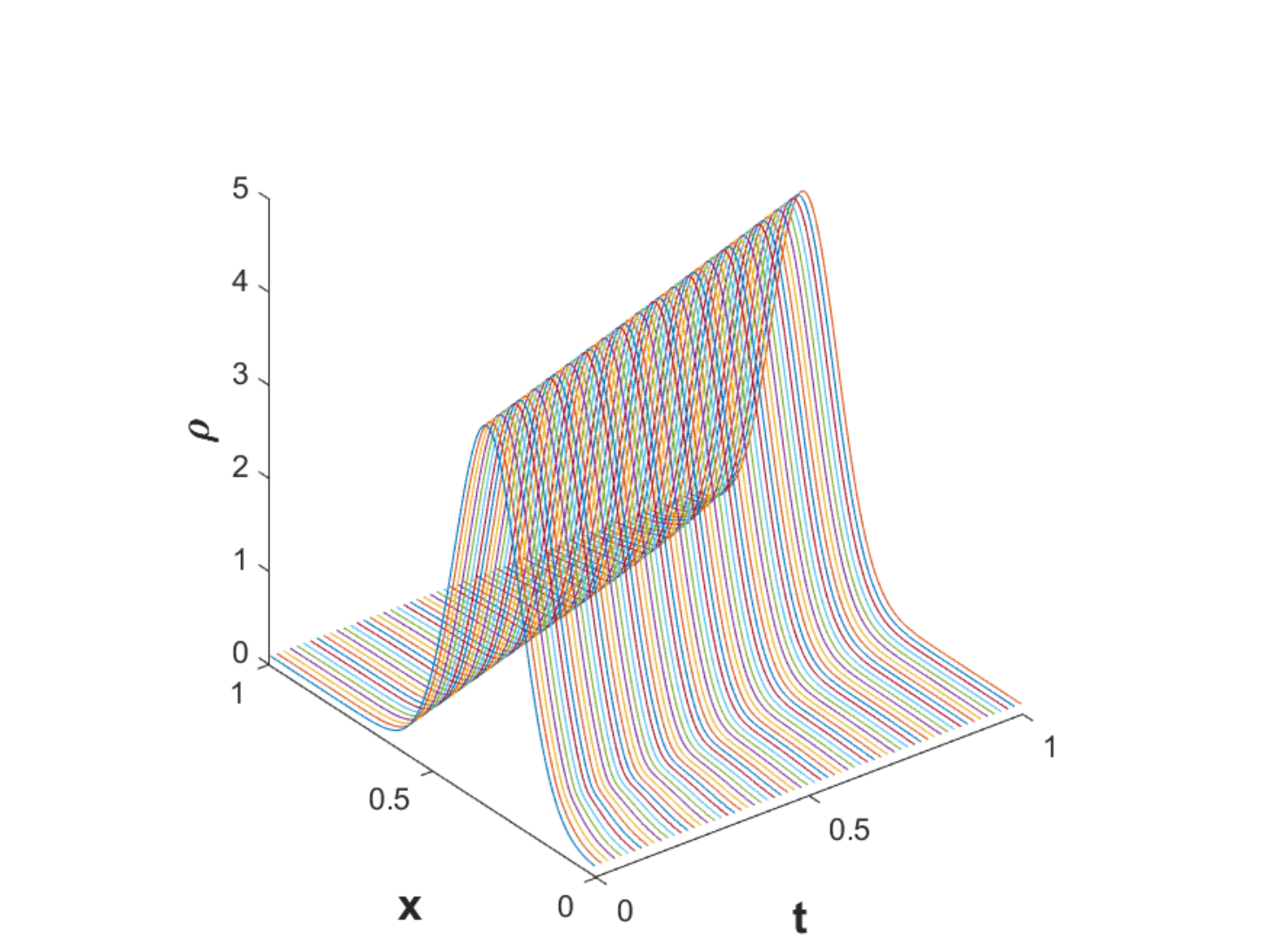}
\includegraphics[width=4cm]{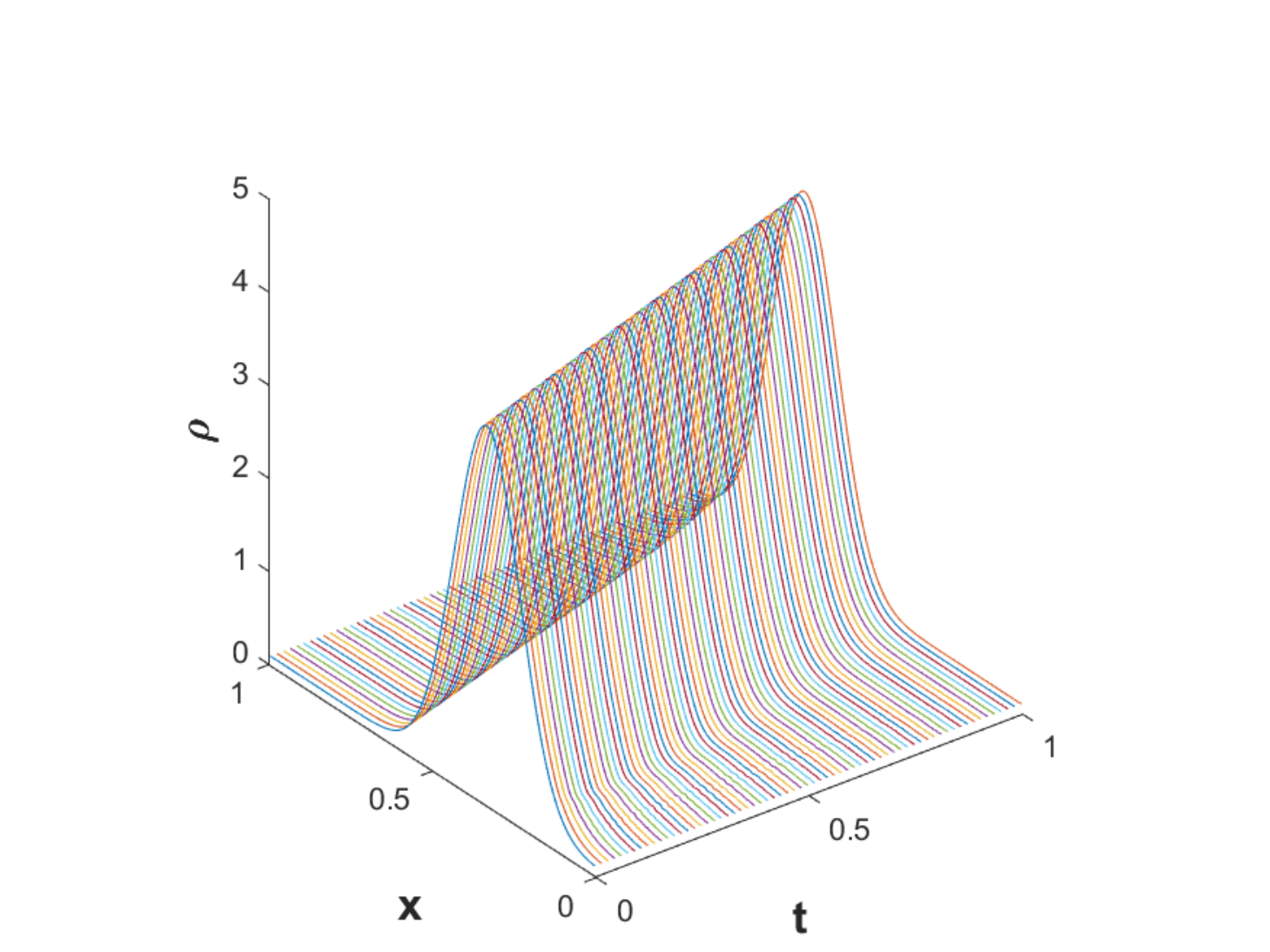}
\includegraphics[width=4cm]{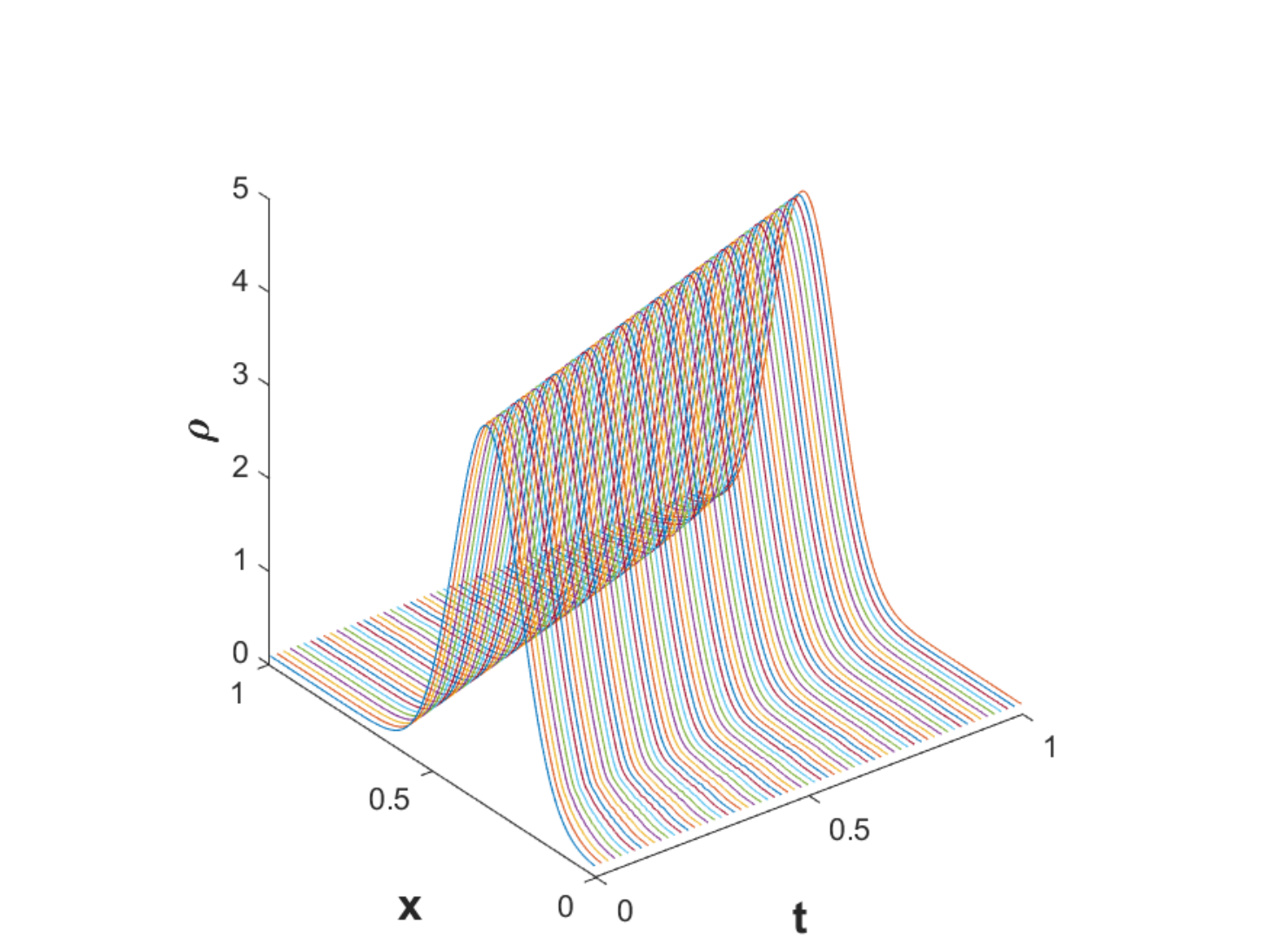}
\caption{Qualitative comparisons of $\rho(t,\cdot)$ in 1D. 
Row 1 from left to right: FISTA, ALG, G-Prox.  
Row 2 from left to right: MLFISTA, MGFISTA$(K=5)$, MGFISTA$(K=10)$.}
\label{fig: 1dot gaussian comp plot3}
\end{figure}

Moreover, as shown in \cref{tab: 1dot gaussian samesize} and \cref{tab: 2dot gaussian samesize}, we further accelerate the proposed algorithm by at most 10 times with the help of multilevel and multigrid strategies. We also compute the residue of being a stationary point, residue of feasibility constraint \cref{eqn: mfp opt constraint}, and residue of mass conservation to check the accuracy of the solutions. From the residue comparisons listed in the tables, it is clear to see that all of our algorithms provide solutions with far more better mass preservation property than results from ALG and G-prox methods due to the nature of the projection step in our method. Qualitatively, \cref{fig: 1dot gaussian comp plot3} also shows that all 6 algorithms in our experiments provide satisfactory results in accuracy.

\subsection{MFP with Obstacles}
\label{sec: num comp mfp}

Most numerical examples of MFP in literature consider $\Omega$ to be a regular region, i.e. $\Omega = [0,1]\times[0,1]$. However, in real application, problems defined in irregular regions might make the  implementation very complicated. One potential way of handling irregular domain is to set $Q$ to be an indicator function of obstacles which leads to solutions staying in the irregular domain. In a different example, \cite{papadakis2014optimal} provides an interesting optimal transport example where the region is a maze with many ``walls''. Here we consider several illustrative cases where there are one or two pieces of obstacles in our square domain and show that our algorithm can deal with this case without modification of implementation. More detailed study along this direction will be explored in our future work. 

\begin{figure}[htbp]
\centering
\subfigure[Example 1]{
\label{fig: num maze illustration 1}
\includegraphics[width=3cm]{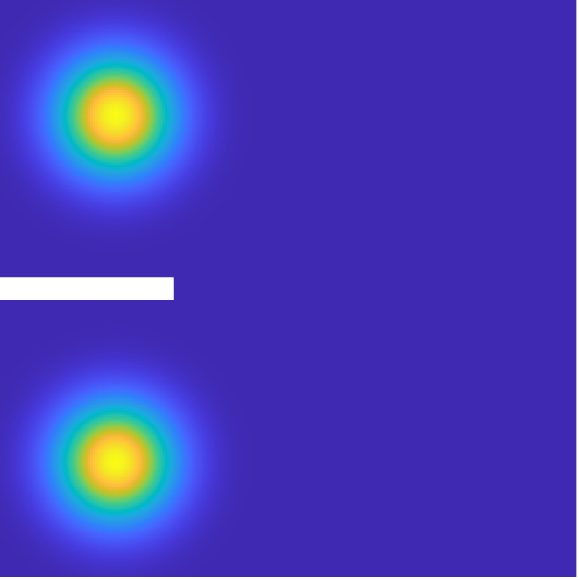}
}
\subfigure[Example 2]{
\label{fig: num maze illustration 2}
\includegraphics[width=3cm]{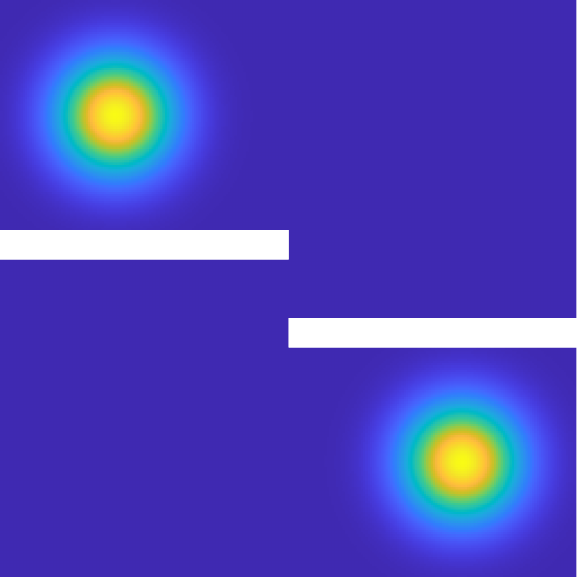}
}
\subfigure[Example 3]{
\label{fig: num maze illustration 3}
\includegraphics[width=3cm]{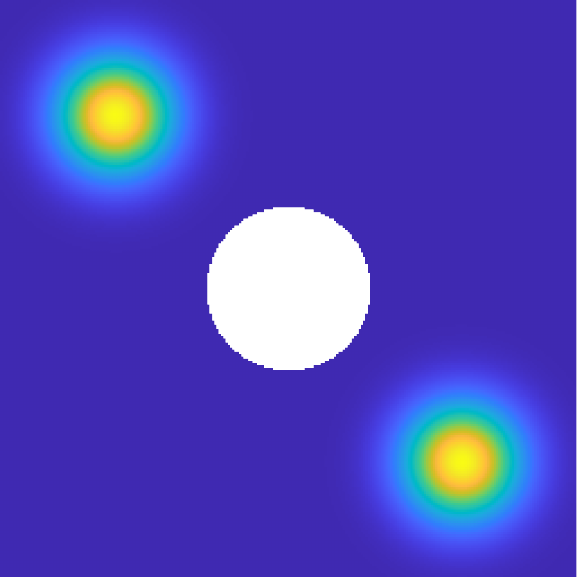}
} \\
\centering  
\subfigure[Example 1]{
\label{fig: num maze snapshot 1}
\includegraphics[width=2.4cm]{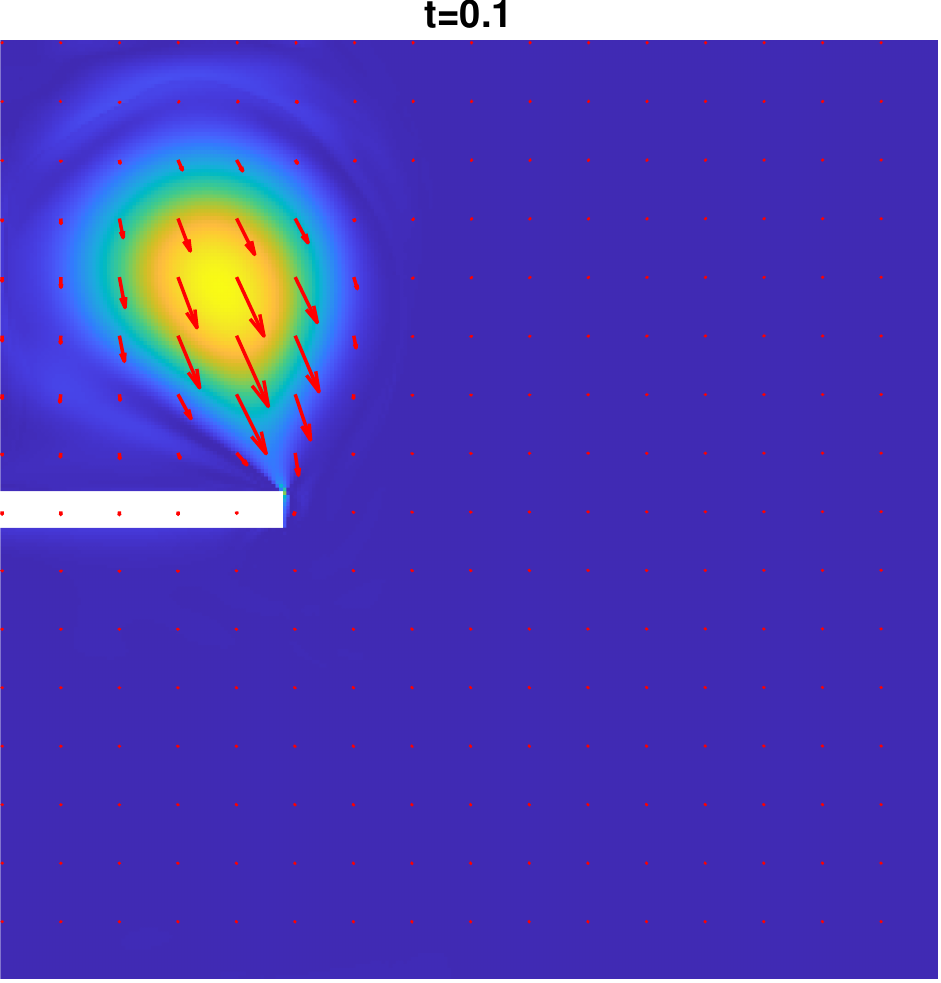}
\includegraphics[width=2.4cm]{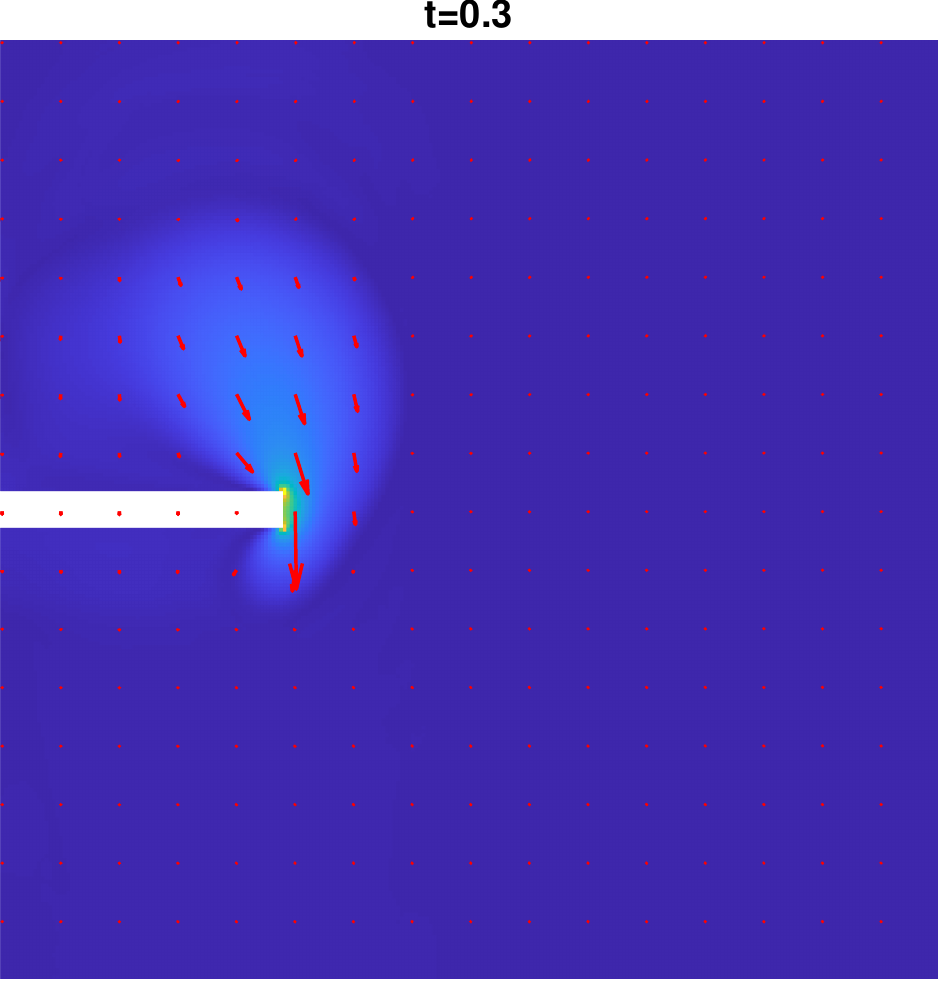}
\includegraphics[width=2.4cm]{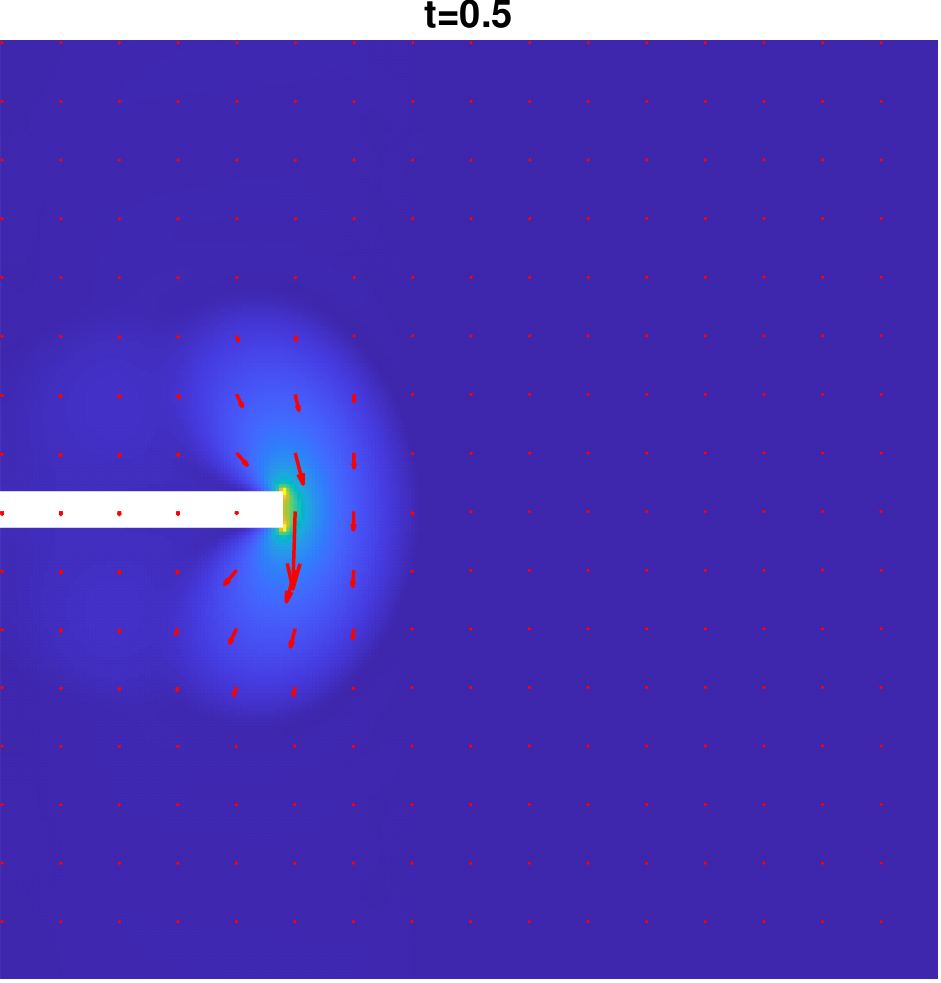}
\includegraphics[width=2.4cm]{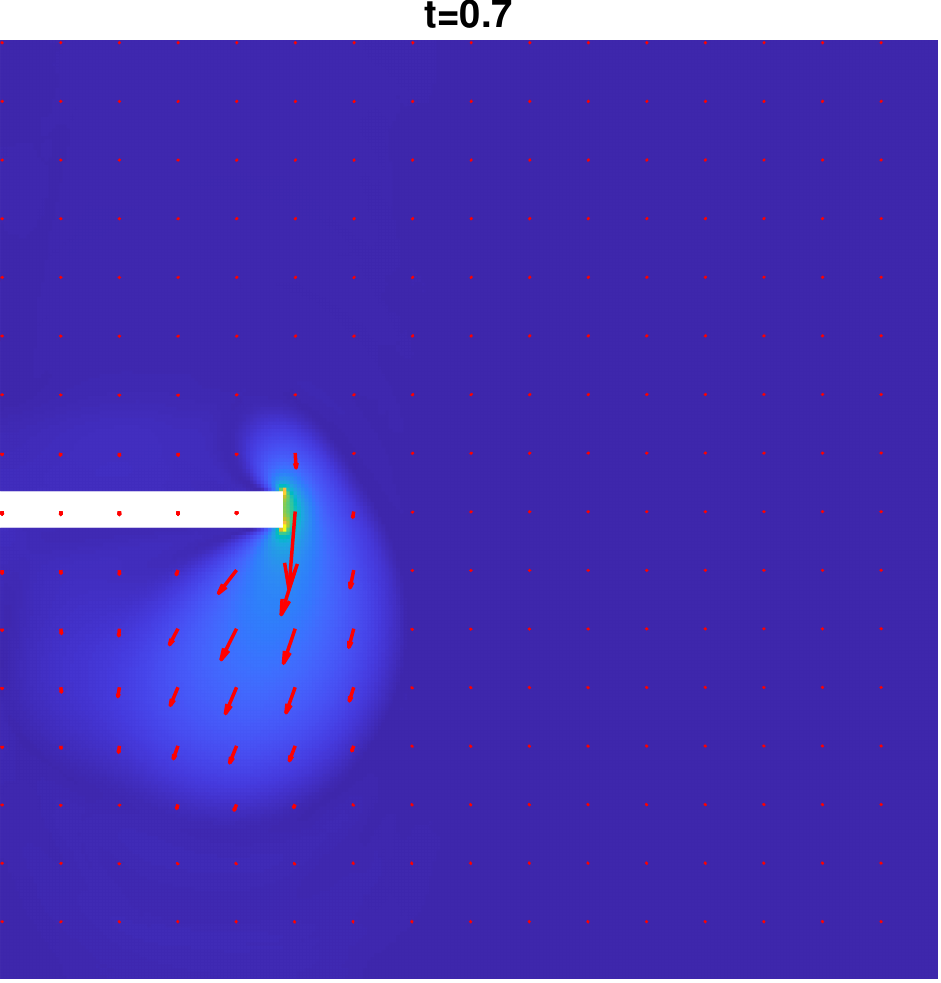}
\includegraphics[width=2.4cm]{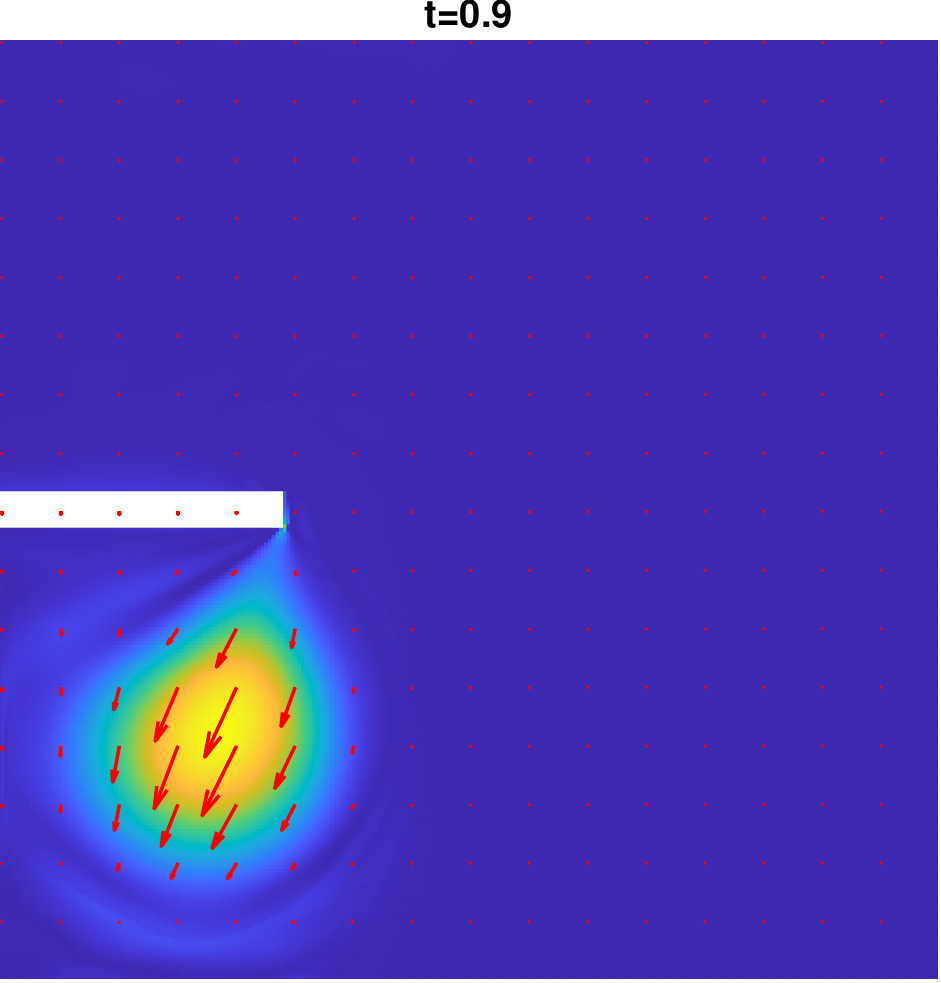}
}
\subfigure[Example 2]{
\label{fig: num maze snapshot 2}
\includegraphics[width=2.4cm]{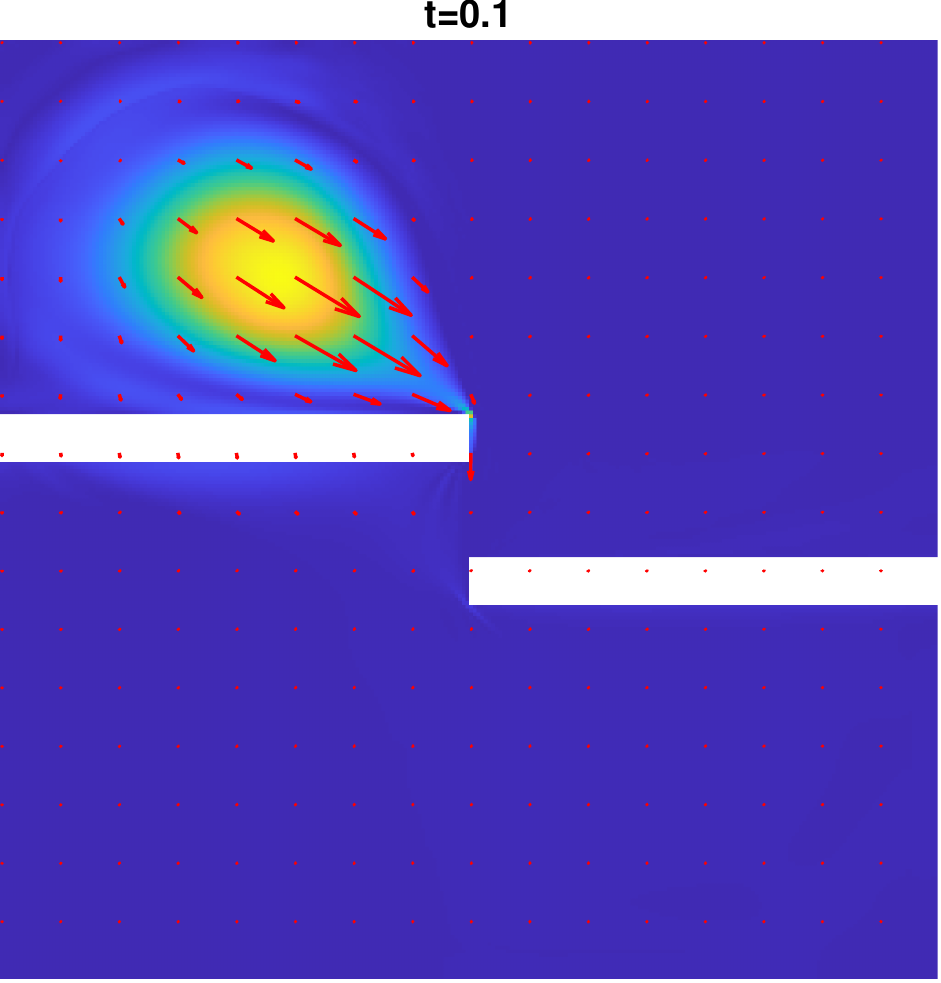}
\includegraphics[width=2.4cm]{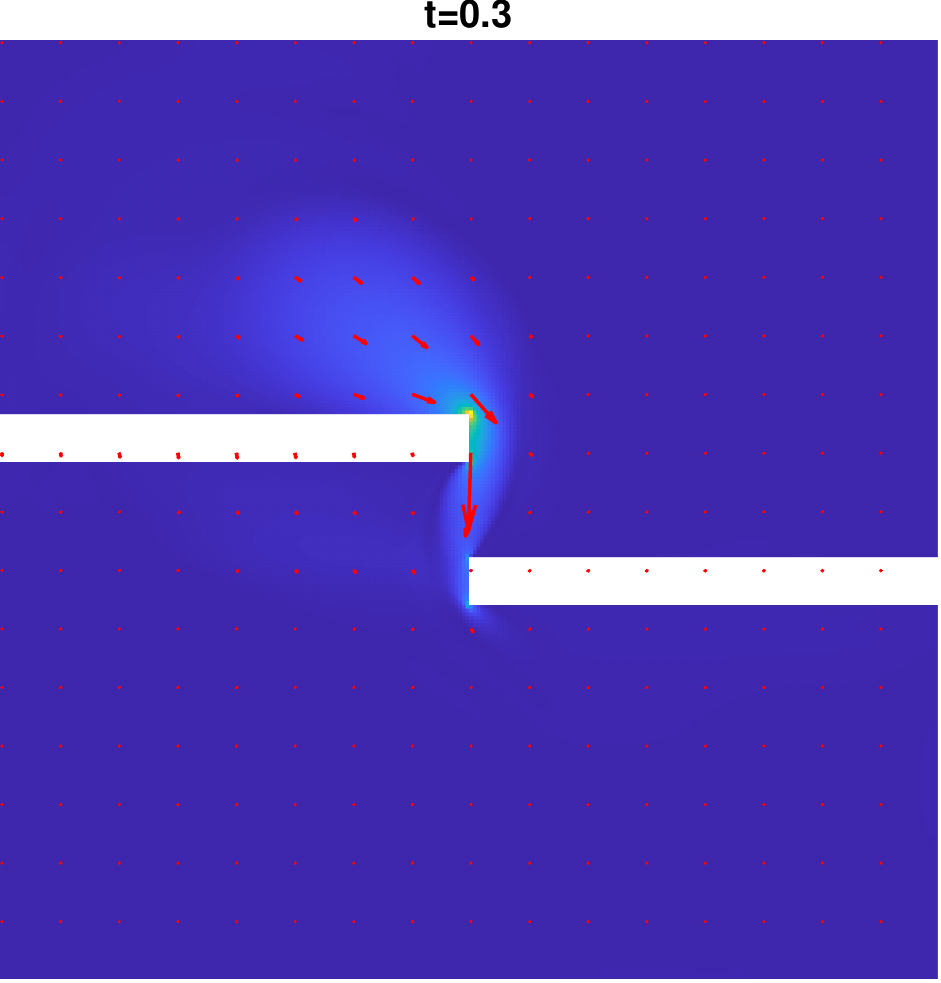}
\includegraphics[width=2.4cm]{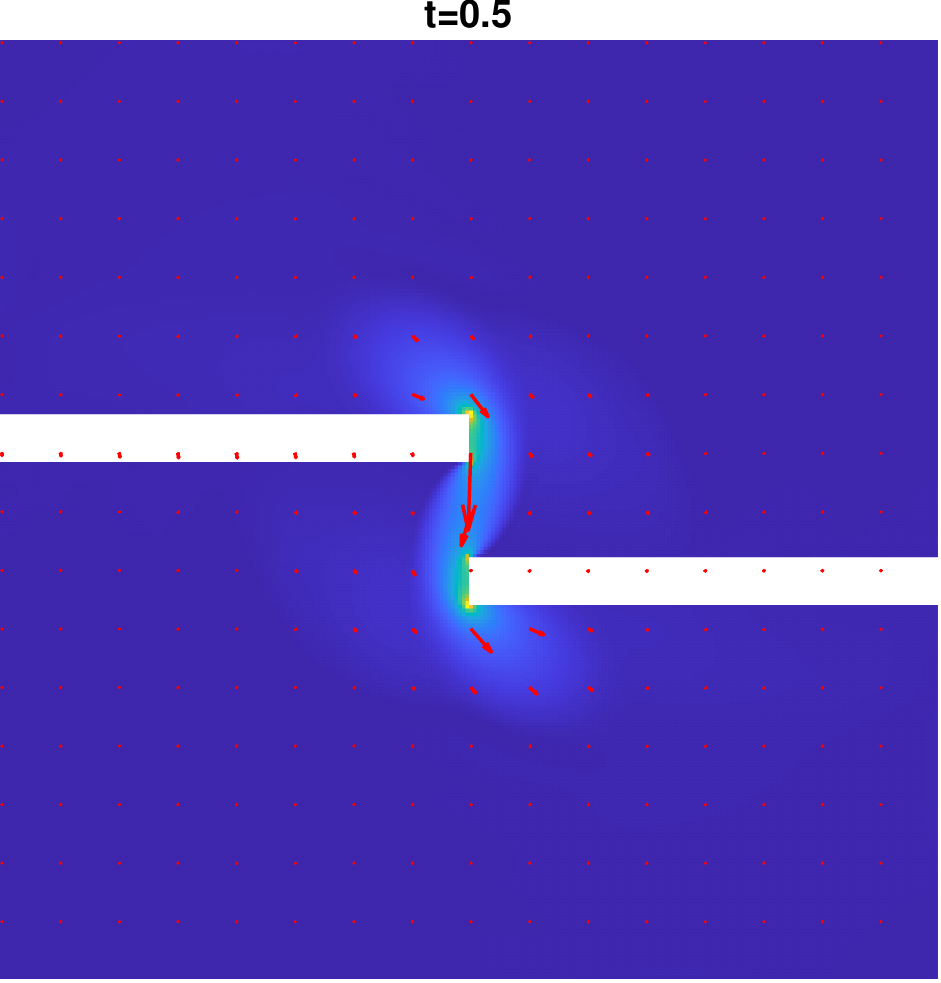}
\includegraphics[width=2.4cm]{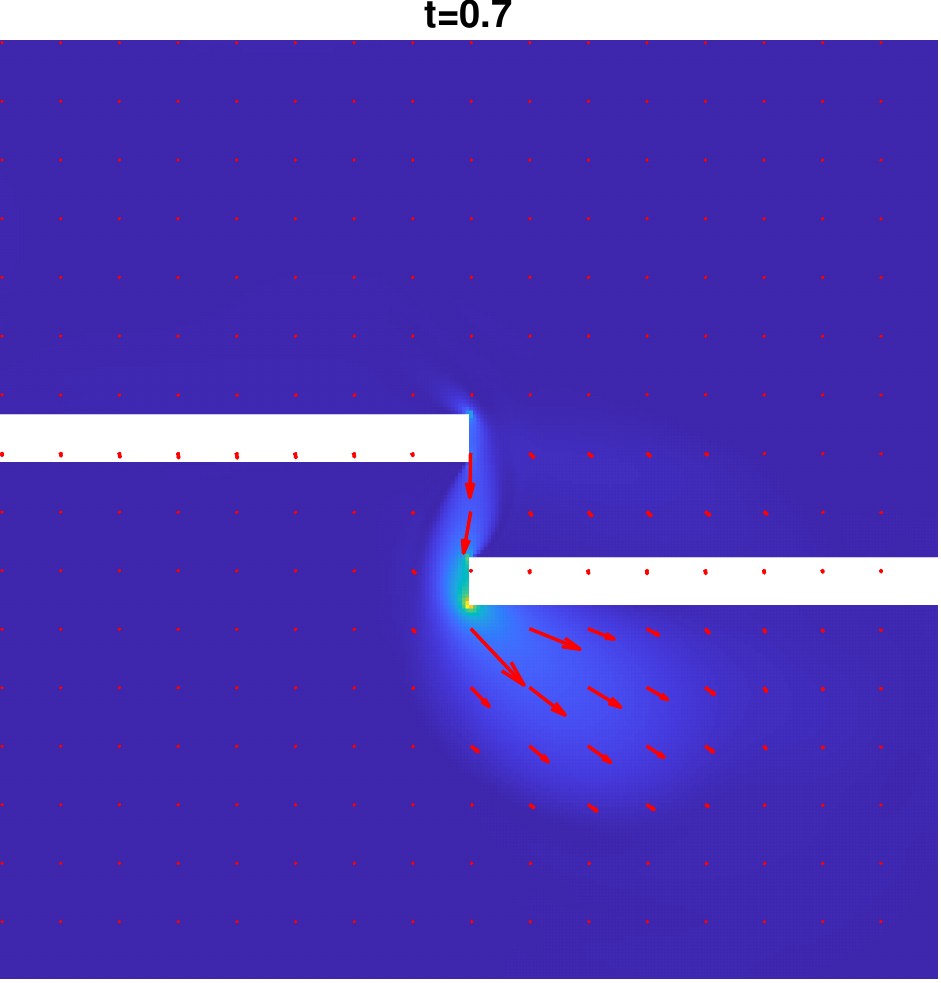}
\includegraphics[width=2.4cm]{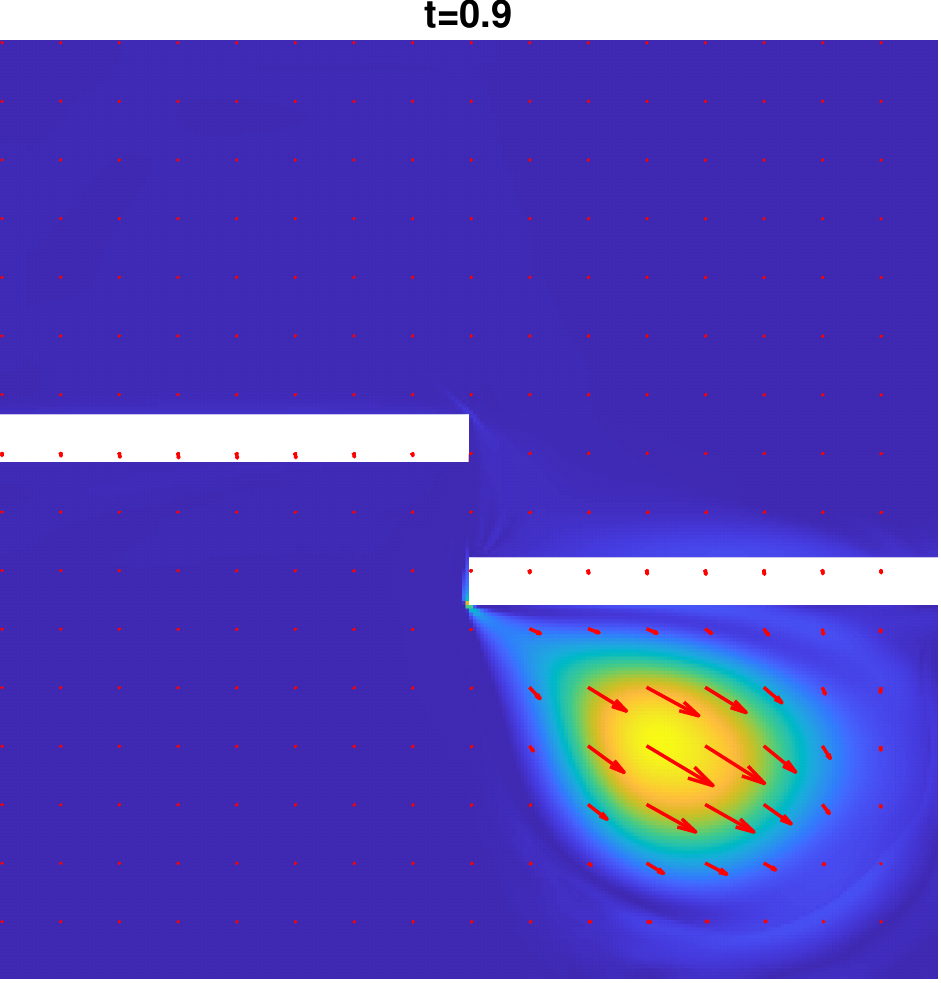}
}
\subfigure[Example 3]{
\label{fig: num maze snapshot 3}
\includegraphics[width=2.4cm]{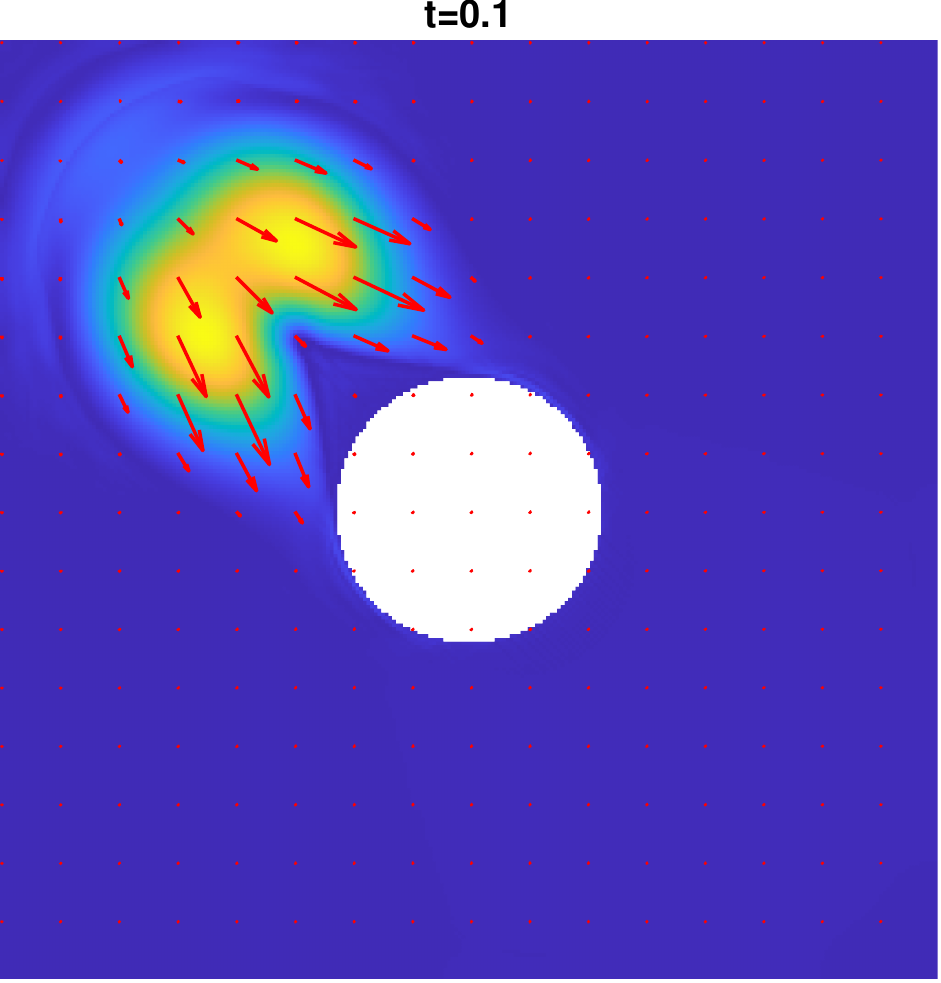}
\includegraphics[width=2.4cm]{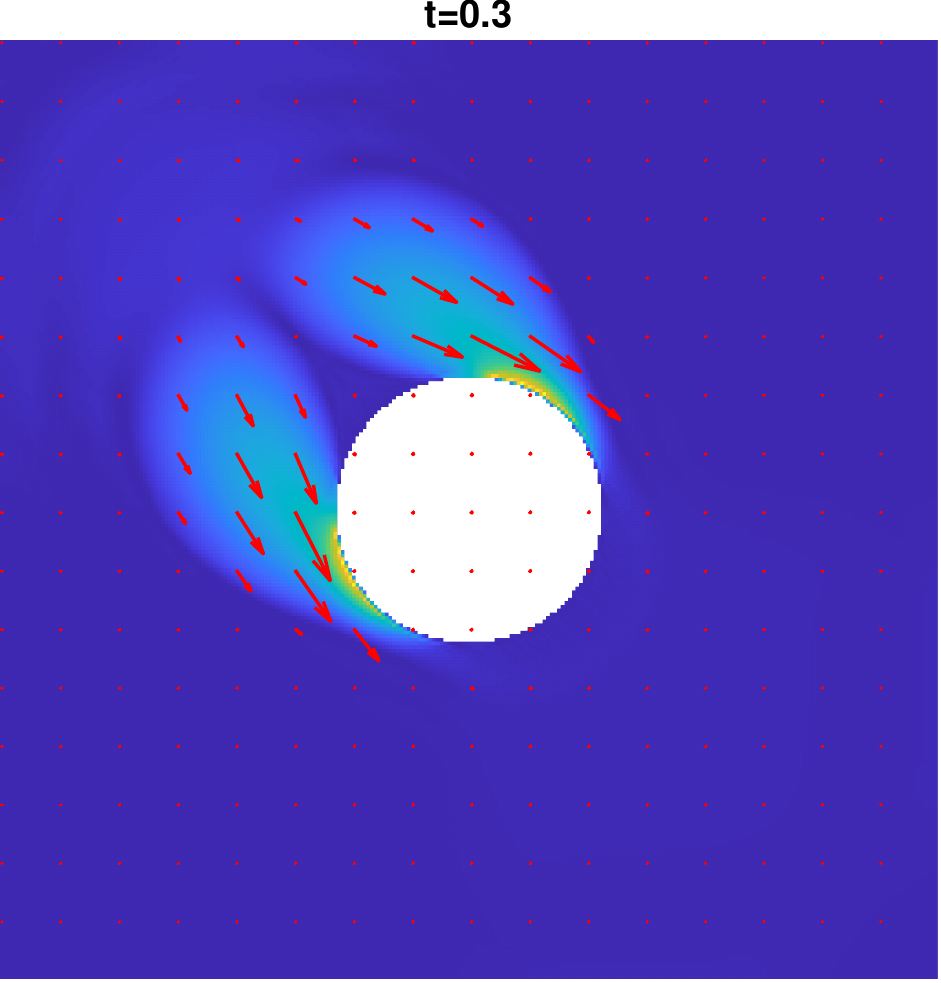}
\includegraphics[width=2.4cm]{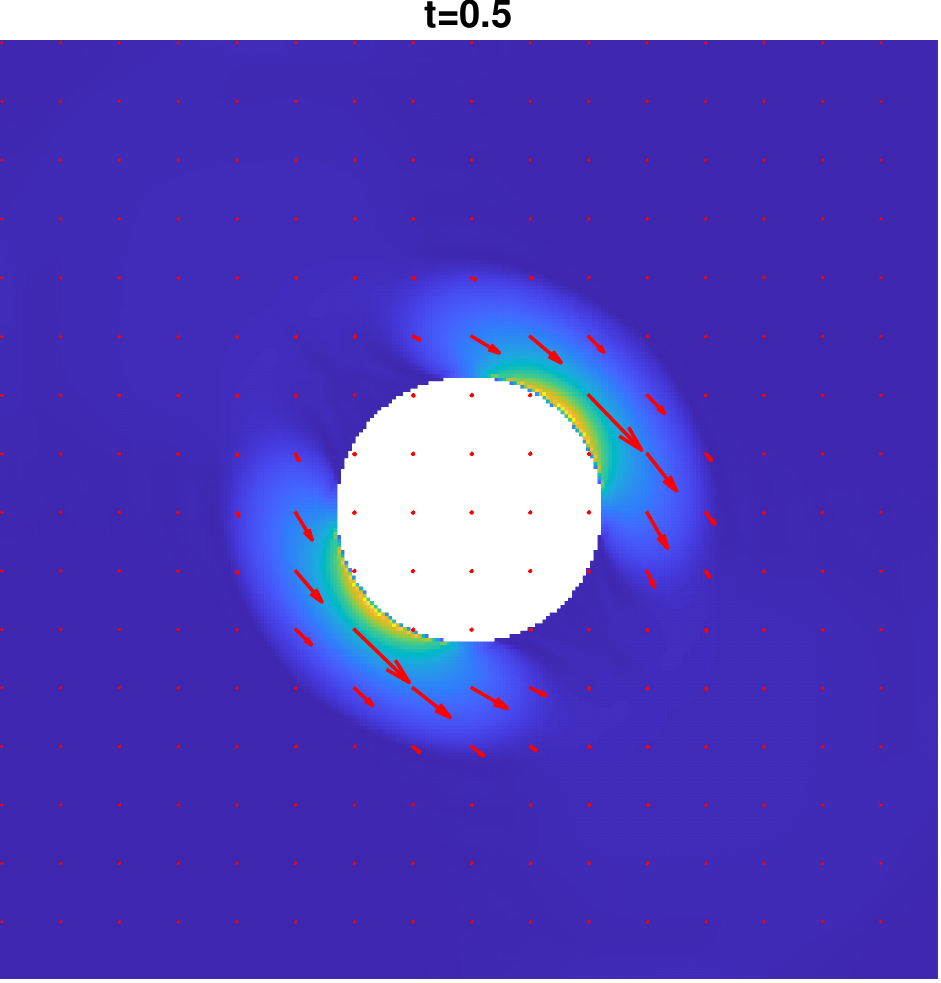}
\includegraphics[width=2.4cm]{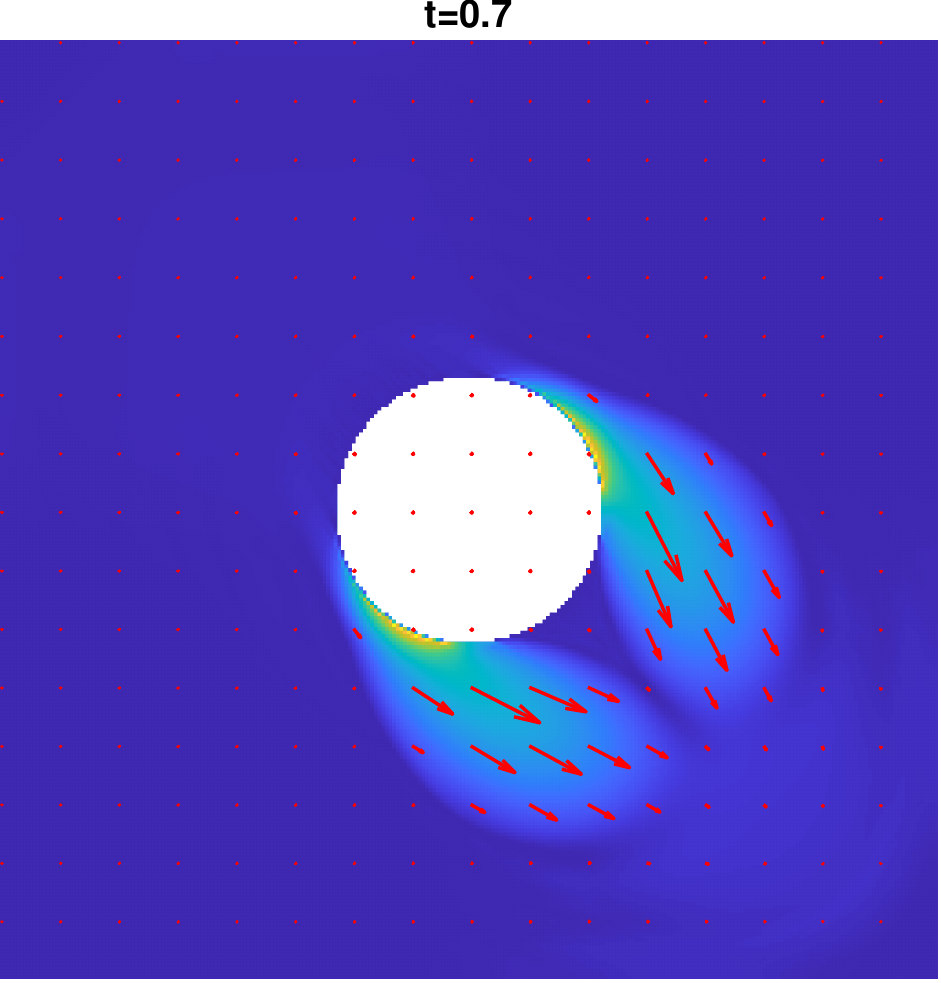}
\includegraphics[width=2.4cm]{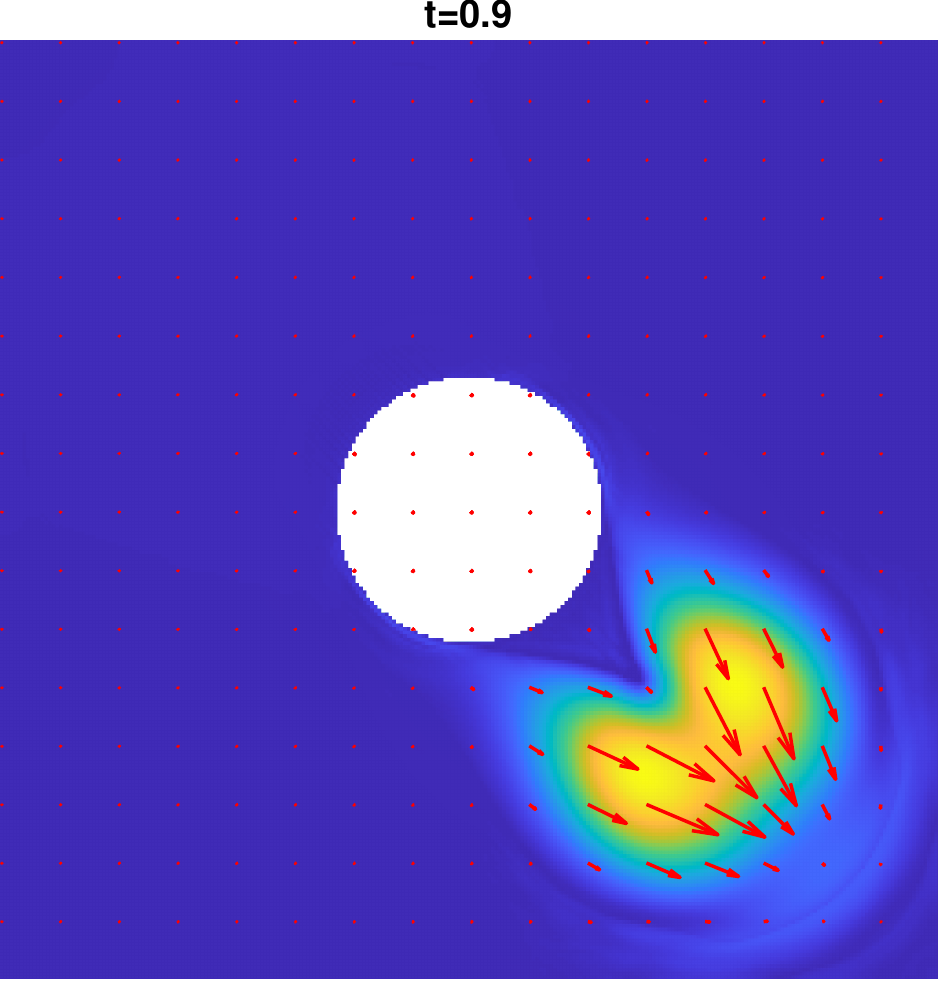}
}
\caption{(a-c): Initial density $\rho_0$, terminal density $\rho_1$ and obstacle region $\Omega_0$ highlighted as white regions. (d-f) Snapshots of $\rho$ at $t=0.1,0.3,0.5,0.7,0.9$. }
\label{fig: num maze snapshot}
\end{figure}

To be precise, letting $\Omega=\left[-\half,\half\right]\times\left[-\half,\half\right]$, we consider MFP problem with objective function 
\begin{equation*}
    \int_0^1\int_\Omega L(\rho(t,\bmx),\bmm(t,\bmx)\deri \bmx \deri t
+ \lambda_Q \int_0^1\int_\Omega \rho(t,\bmx)Q(\bmx) \deri\bmx,
\end{equation*}
Different choices of $\rho_0,\rho_1,Q$ are shown in the first row of \cref{fig: num maze snapshot} and $Q(\bmx) = \begin{cases}
1,\quad \bmx\in \Omega_0\\
0,\quad \bmx \not\in \Omega_0
\end{cases}$ where $\Omega_0$ is the white region.
By setting $\lambda_Q$ to be a very large number (e.g. $\lambda_Q=8\times10^4$ in our implementation), we expect the set $\Omega_0$ to be viewed as an obstacle and the density evolution to circumvent the region. The snapshots of the evolution shown in \cref{fig: num maze snapshot} demonstrate the success of our algorithm that the mass circumvents the obstacles very well.


\subsection{Flexibility}
\label{sec: num comp image}

As one of the greatest advantages, our method enjoys flexibility to  handle different types of objective functions in variational MFP problems.  To show the effectiveness of our algorithm, we apply \cref{alg: fista disct} to the five models listed in \cref{sec: review}.   We can also observe how different objective functions affect the density evolutions. 

\begin{figure}[htbp]
\centering
\includegraphics[width=3cm]{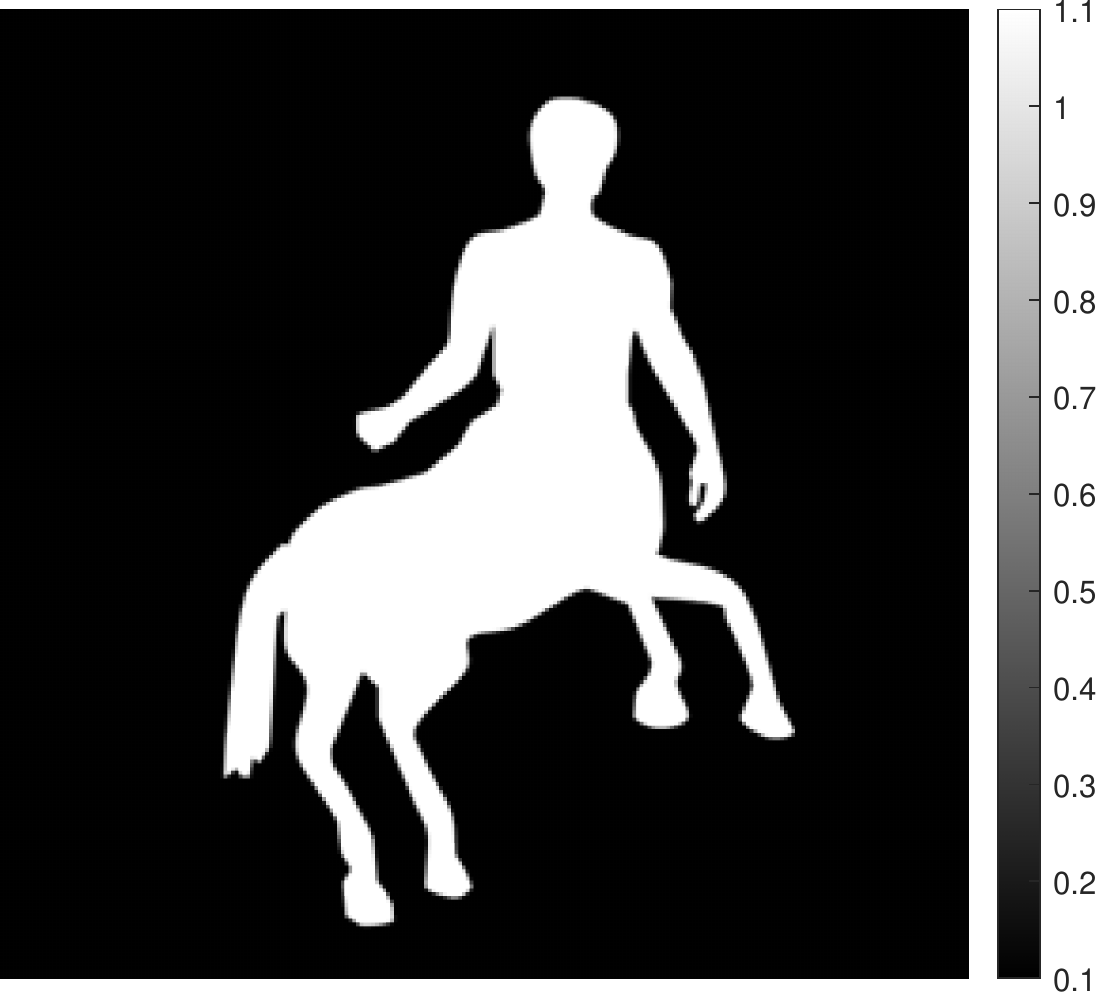}
\includegraphics[width=3cm]{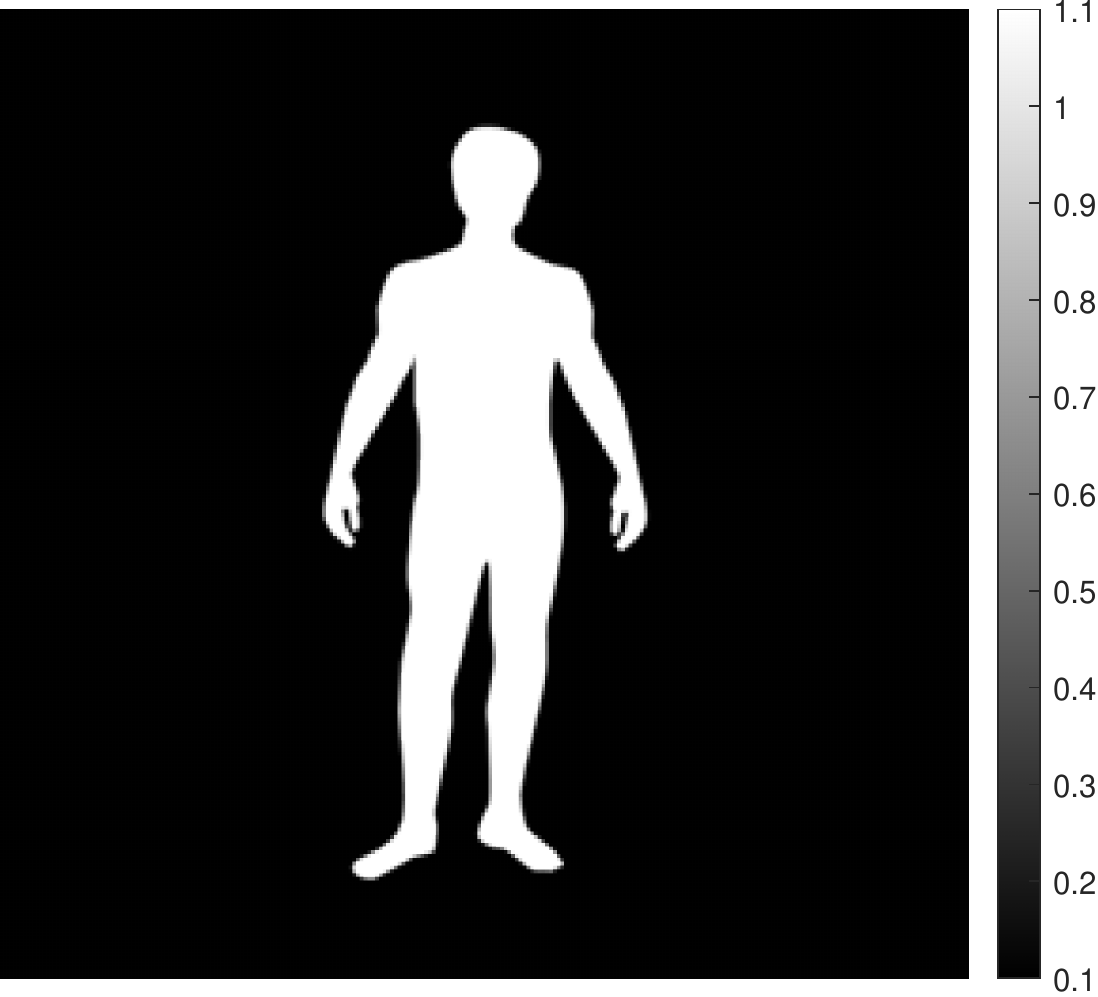}
\includegraphics[width=3cm]{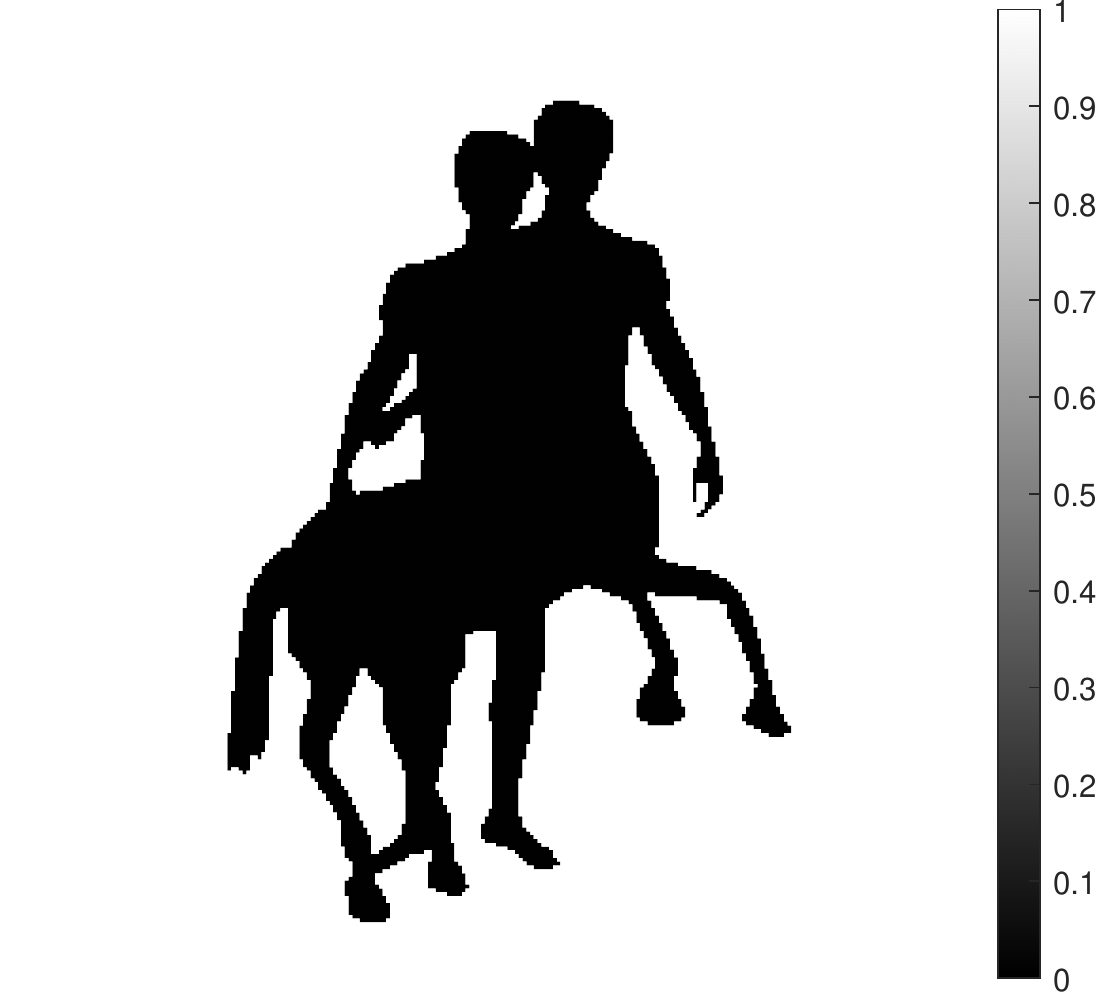}
\includegraphics[width=3cm]{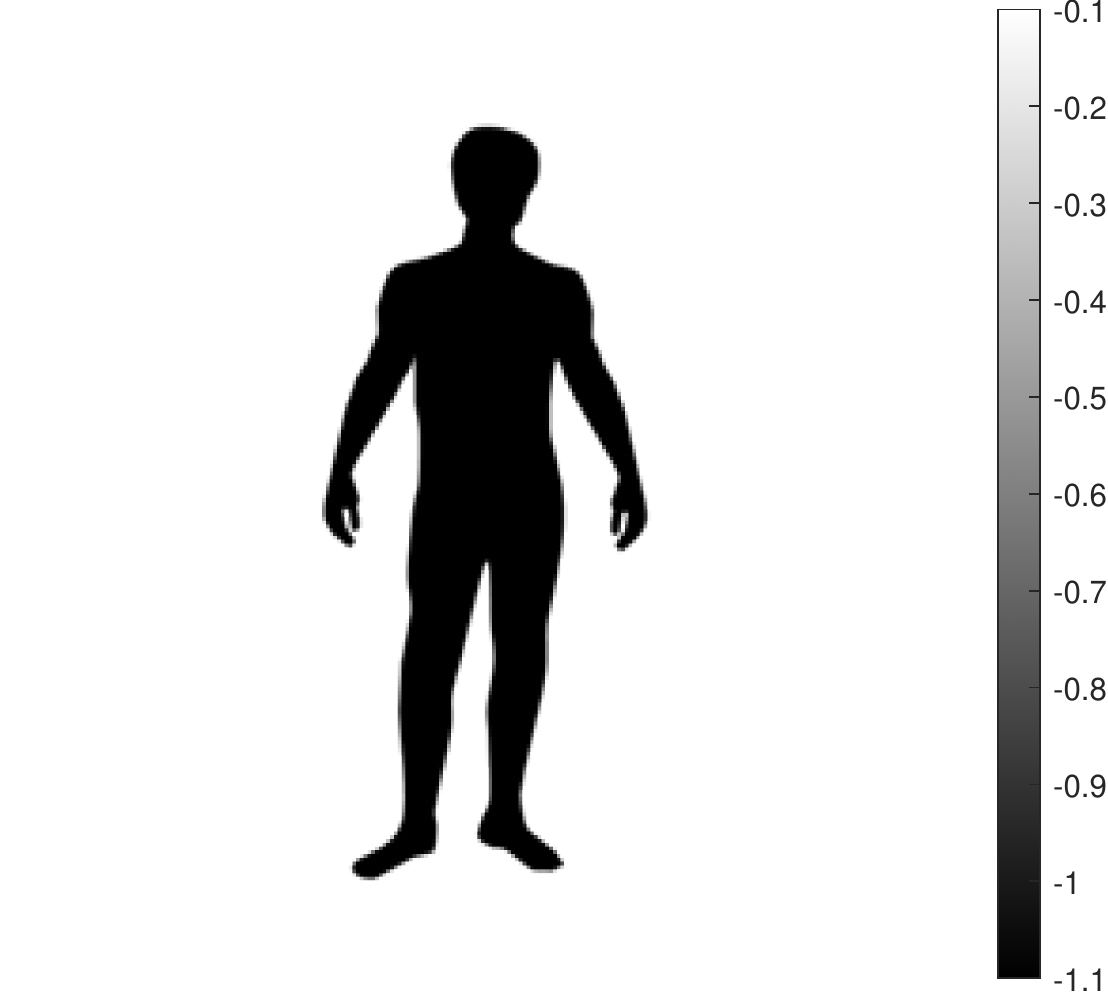}
\caption{From left to right: initial density $\rho_0$, final $\rho_1$, interaction penalty $Q(\bmx)$ and terminal density regularizer $G(\bmx)$}
\label{fig: num image illustration}
\end{figure}

Let $\Omega = [0,1]\times[0,1],$ $\rho_0,\rho_1$ being two images shown in \cref{fig: num image illustration}, $G(\bmx) = -\rho_1(\bmx)$ and 
$Q(\bmx) = \begin{cases}
    0,\quad \rho_0(\bmx) \neq 0 \text{ or } \rho_1(\bmx) \neq 0,\\
    1,\quad \text{ otherwise },
\end{cases}$. 
We consider MFP problem of the following form
\begin{equation*}
    \min_{(\rho,\bmm)\in\calC(\rho_0,\rho_1)} 
  \left\{\begin{aligned}
    &\int_0^1\int_\Omega L(\rho(t,\bmx),\bmm(t,\bmx)) \deri \bmx \deri t   \\
    +\lambda_E&\int_0^1\int_\Omega F_E(\rho(t,\bmx))\deri\bmx\deri t+ \lambda_Q\int_0^1\int_\Omega\rho(t,\bmx)Q(\bmx)\deri\bmx\deri t
    \end{aligned}\right\},
\end{equation*}
We apply the proposed algorithm to the following four MFP models:
\begin{align}
\text{(OT) } 
&\lambda_E=\lambda_Q=0, \label{eqn: num eg ot}\\
 \text{(Model 1) } 
&\lambda_E=0.01,\lambda_Q=0.1,F_E:\bbRplus\to\bbR, \rho\mapsto\begin{cases}
    \rho\log(\rho),\quad \rho>0\\
    0,\quad \rho=0
\end{cases} \label{eqn: num eg model1} \\
 \text{(Model 2) } 
&\lambda_E=0.01,\lambda_Q=0.1,F_E:\bbRplus\to\bbR,\rho\mapsto\frac{\rho^2}{2}, \label{eqn: num eg model2}\\
 \text{(Model 3) } 
&\lambda_E=0.01,\lambda_Q=0.1,F_E:\bbRplus\to\bbR,\rho\mapsto\begin{cases}
    \rho,\quad \rho>0\\
    0,\quad \rho=0
\end{cases} \label{eqn: num eg model3}
\end{align}
and a MFG model 
\begin{equation}
\min_{(\rho,\bmm)\in\calC(\rho_0)} 
\left\{\begin{aligned}
    &\int_0^1\int_\Omega L(\rho(t,\bmx),\bmm(t,\bmx)) \deri \bmx \deri t   \\
    +\lambda_E&\int_0^1\int_\Omega \rho(t,\bmx)\log(\rho(t,\bmx))\deri\bmx\deri t
    + \lambda_Q\int_0^1\int_\Omega\rho(t,\bmx)Q(\bmx)\deri\bmx\deri t \\
     + \lambda_G&\int_\Omega \rho(1,\bmx)G(\bmx)\deri\bmx.
\end{aligned}\right\}
\label{eqn: num eg mfg}
\end{equation}
with $\lambda_E=0.01,\lambda_Q=0.1,\lambda_G=1$. 
It is worth mentioning that to solve model \cref{eqn: num eg ot}-\cref{eqn: num eg model3}, we must rescale $\rho_1,\rho_1$ such that $\int_\Omega\rho_0 = \int_\Omega\rho_1$ but we do not have to rescale $G(\bmx)$ for $\rho_0$ in \cref{eqn: num eg mfg}.

\begin{figure}[h]
\centering  
\subfigure[OT: $F_E(a)=F_0(a) :=0.$]{
\label{fig: num image snapshot m0}
\includegraphics[width=2.4cm]{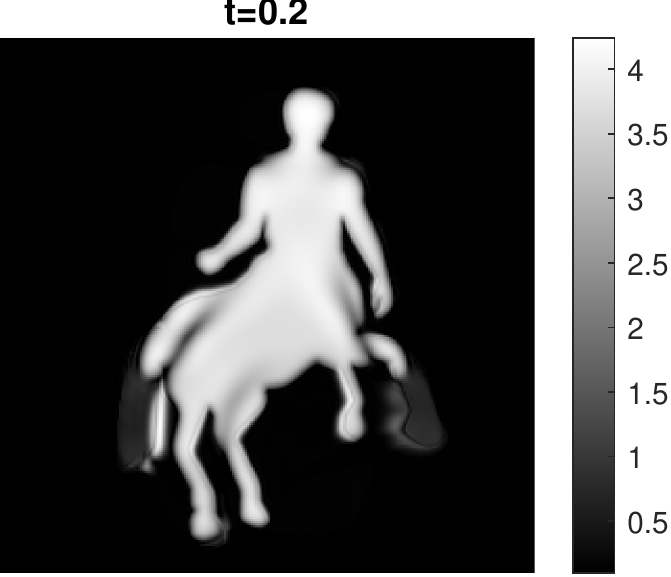}
\includegraphics[width=2.4cm]{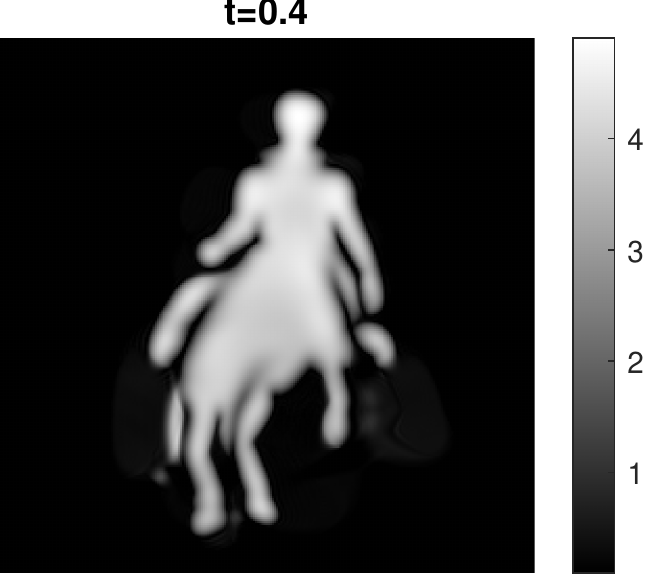}
\includegraphics[width=2.4cm]{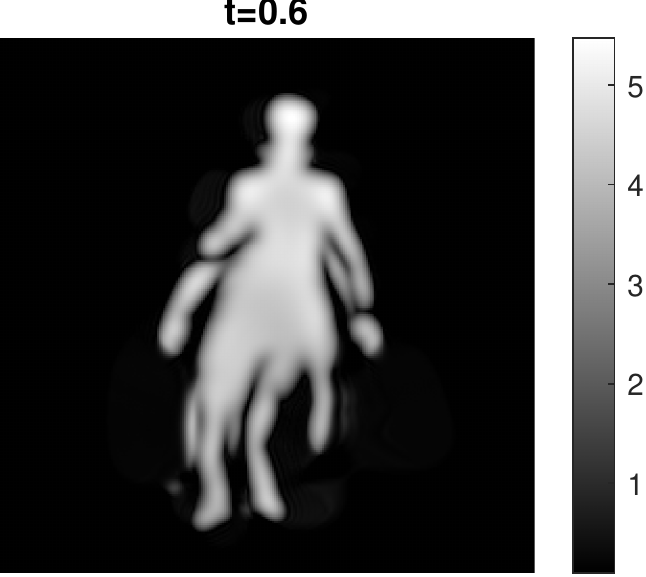}
\includegraphics[width=2.4cm]{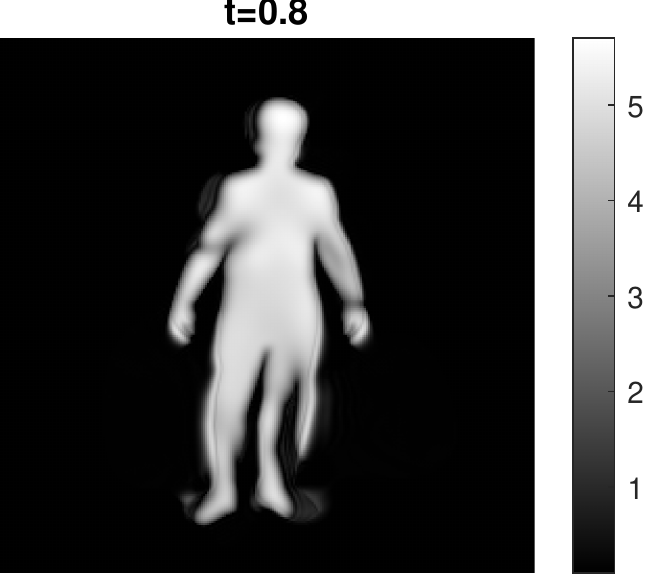}
\includegraphics[width=2.4cm]{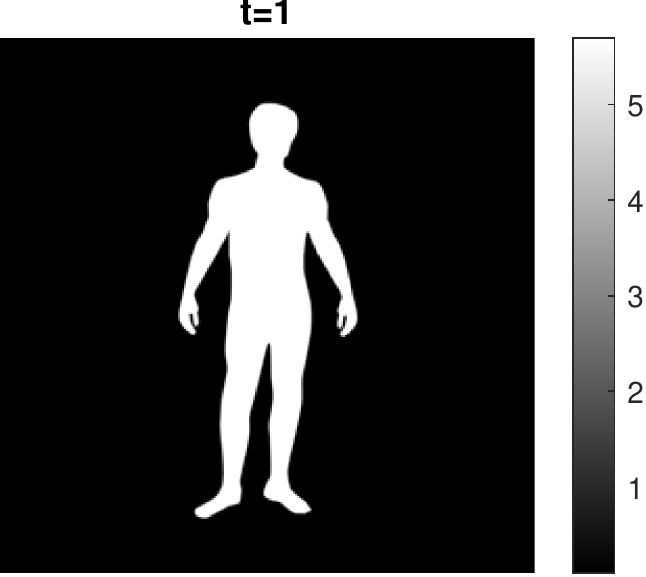}
}
\subfigure[Model 1:$F_E(a)= a\log a, a>0.$]{
\label{fig: num image snapshot m1}
\includegraphics[width=2.4cm]{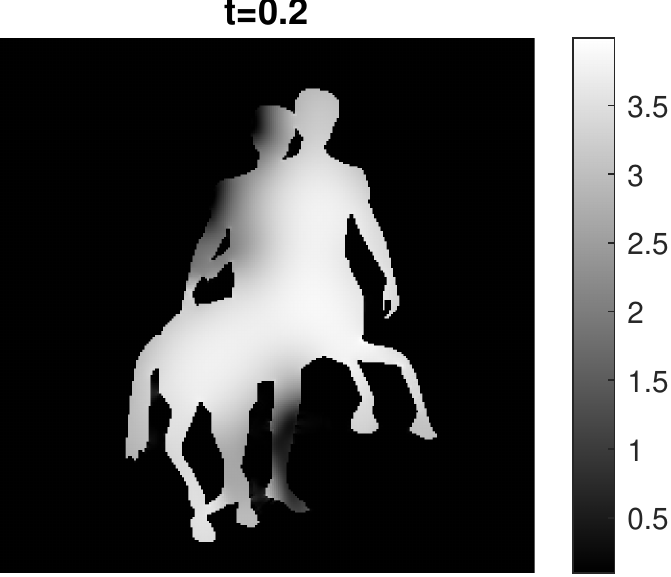}
\includegraphics[width=2.4cm]{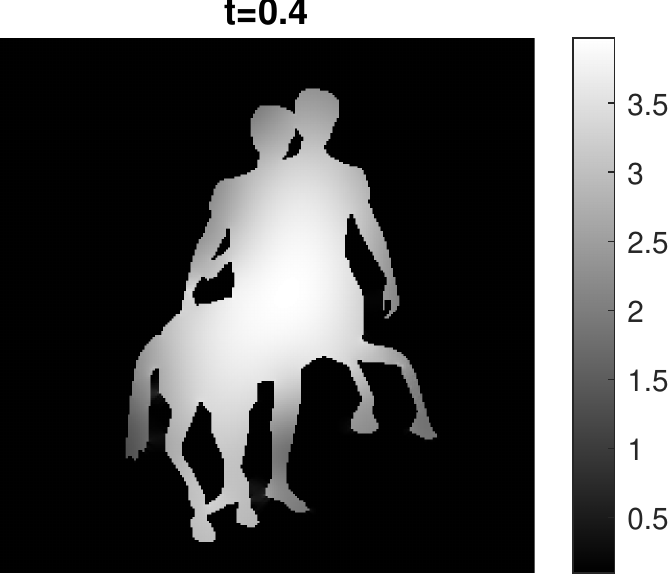}
\includegraphics[width=2.4cm]{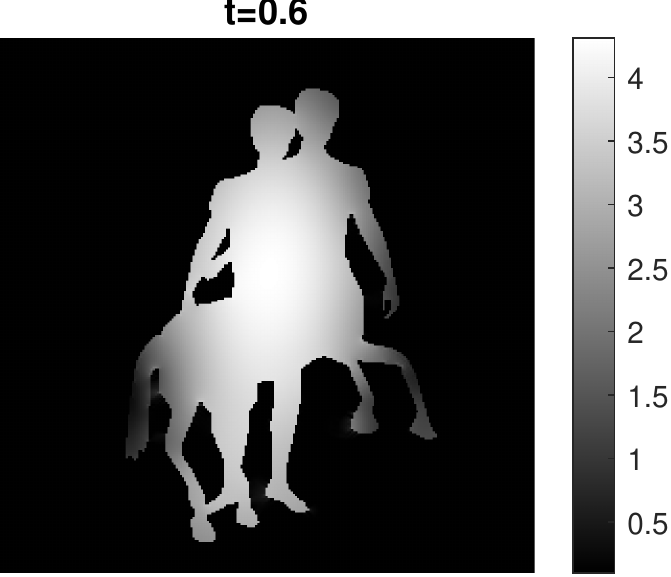}
\includegraphics[width=2.4cm]{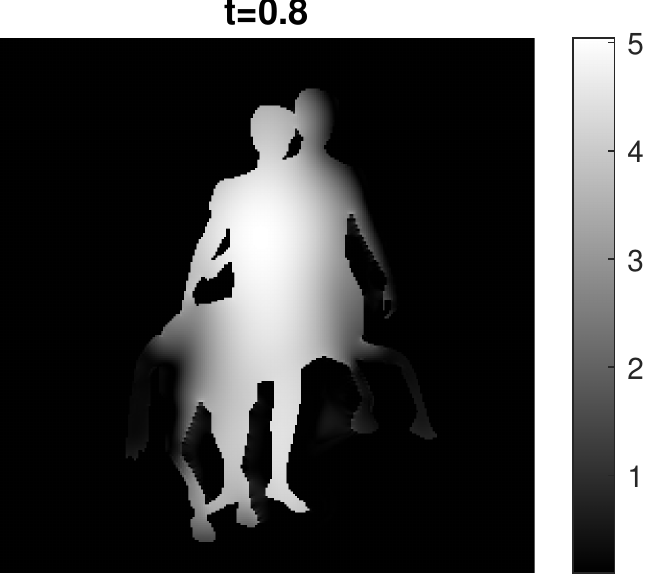}
\includegraphics[width=2.4cm]{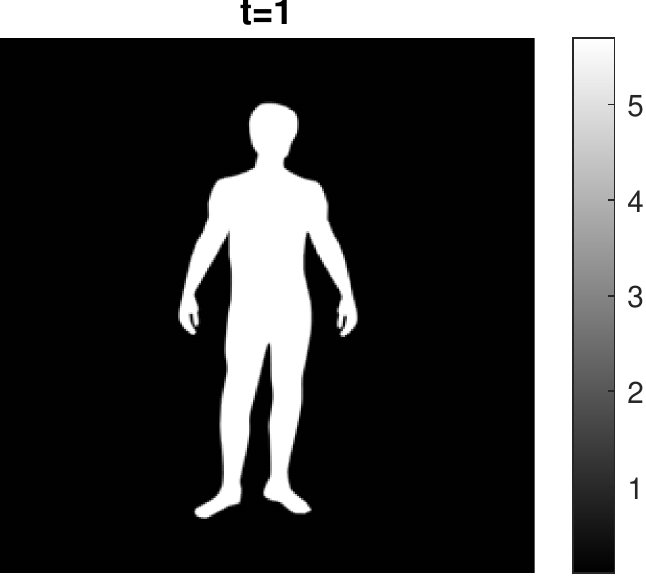}
}
\subfigure[Model 2: $F_E(a)= \frac{a^2}{2}.$]{
\label{fig: num image snapshot m2}
\includegraphics[width=2.4cm]{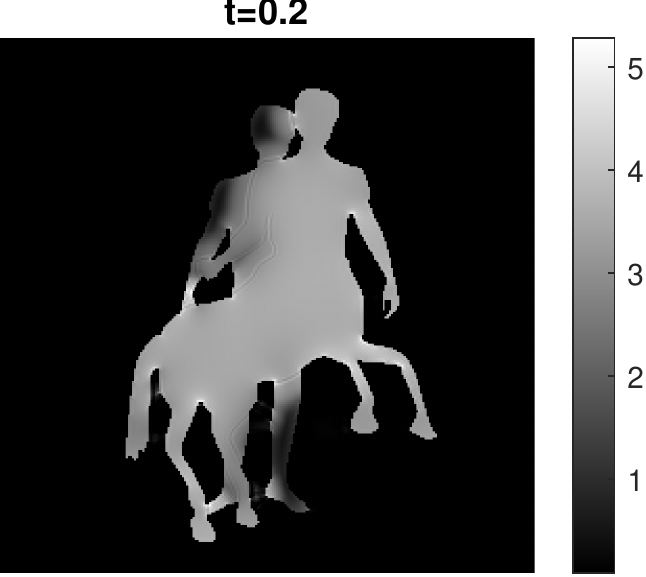}
\includegraphics[width=2.4cm]{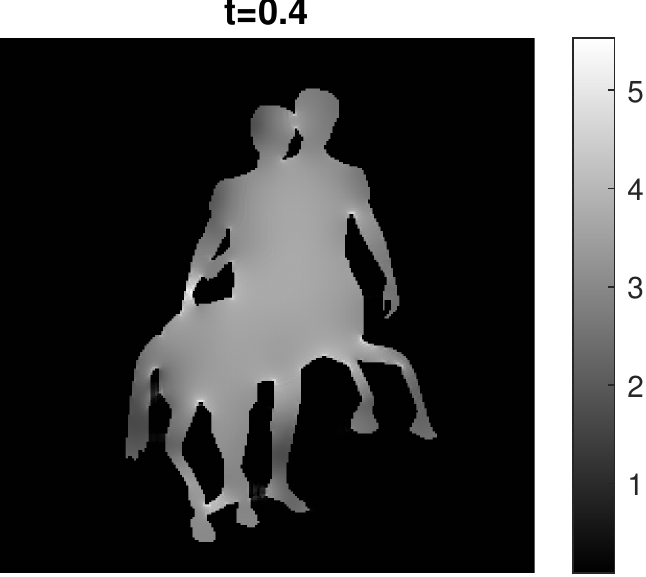}
\includegraphics[width=2.4cm]{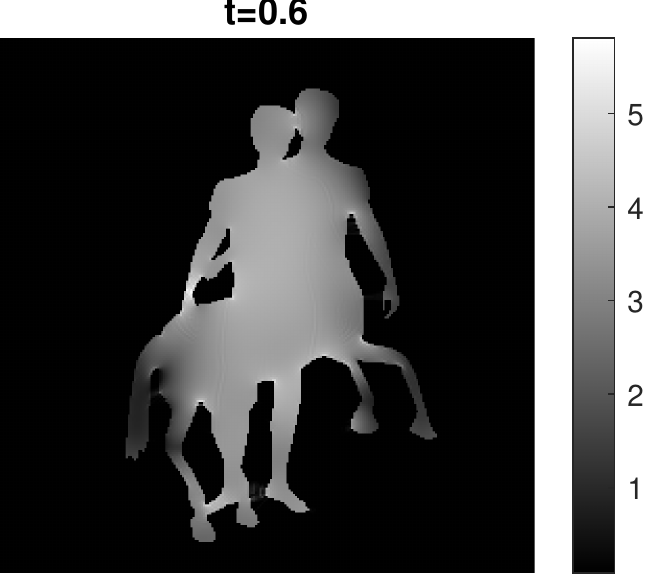}
\includegraphics[width=2.4cm]{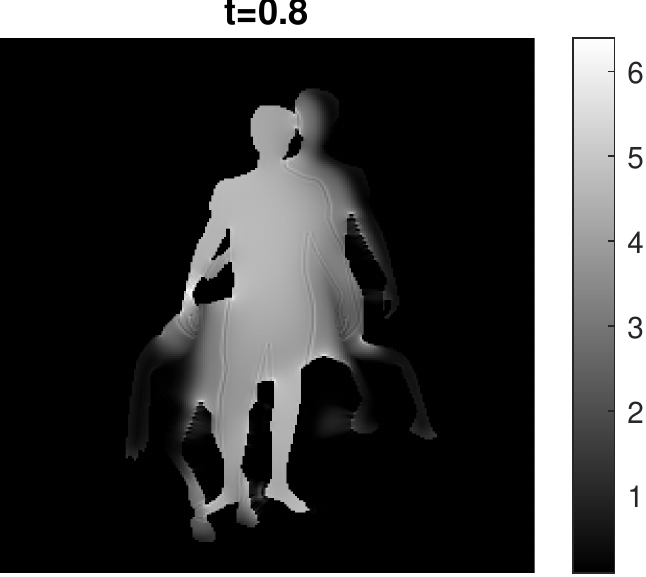}
\includegraphics[width=2.4cm]{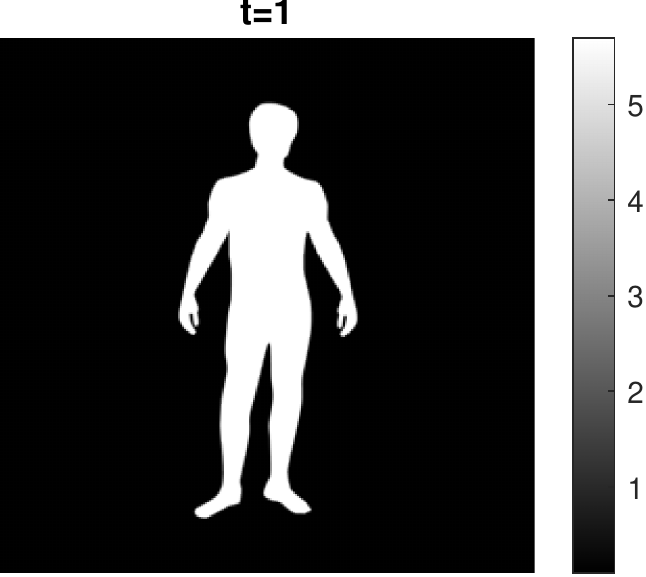}
}
\subfigure[Model 3: $F_E(a)= \frac{1}{a}, a>0.$]{
\label{fig: num image snapshot m3}
\includegraphics[width=2.4cm]{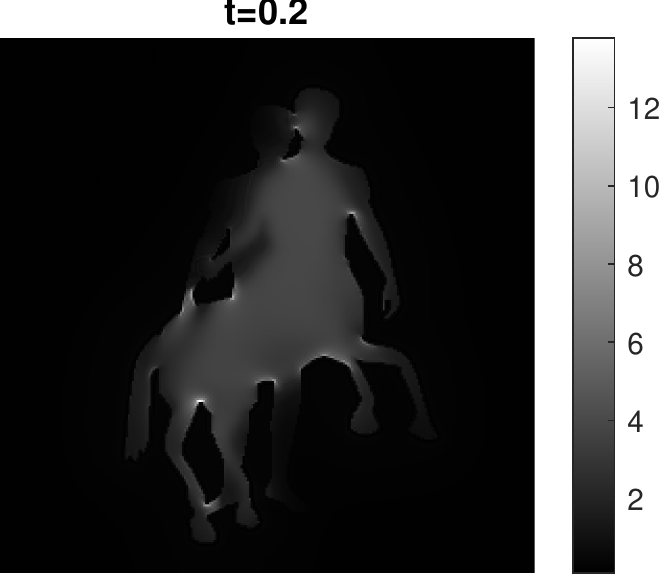}
\includegraphics[width=2.4cm]{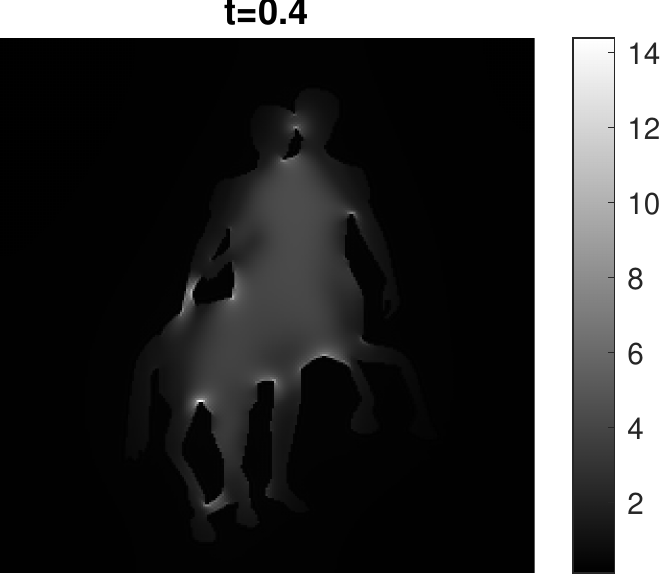}
\includegraphics[width=2.4cm]{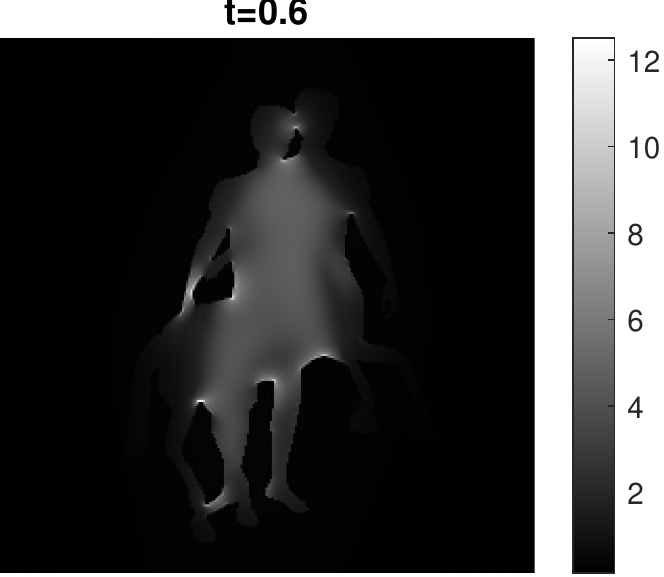}
\includegraphics[width=2.4cm]{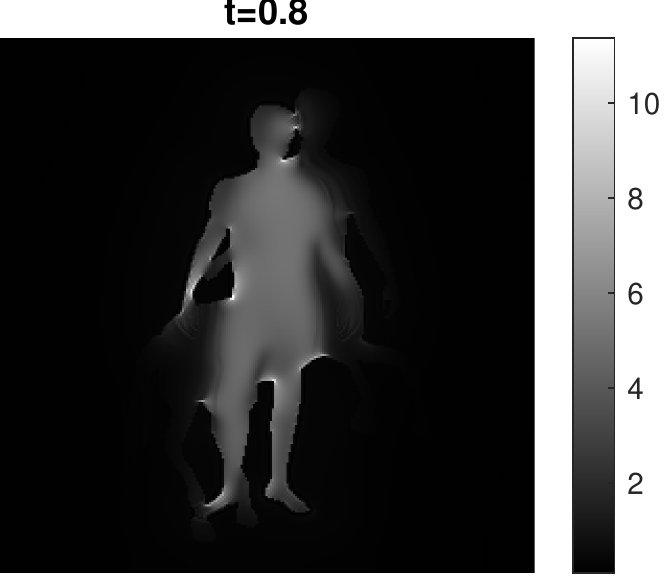}
\includegraphics[width=2.4cm]{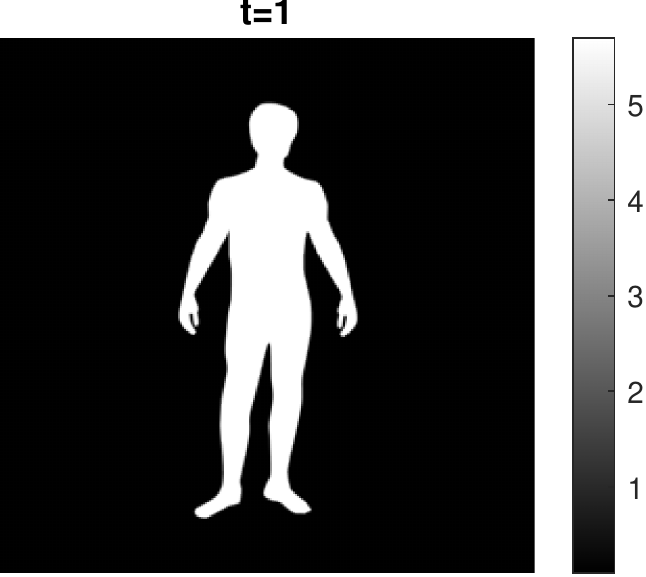}
}
\subfigure[MFG: $F_E(a)= a\log a, a>0.$]{
\label{fig: num image snapshot m4}
\includegraphics[width=2.4cm]{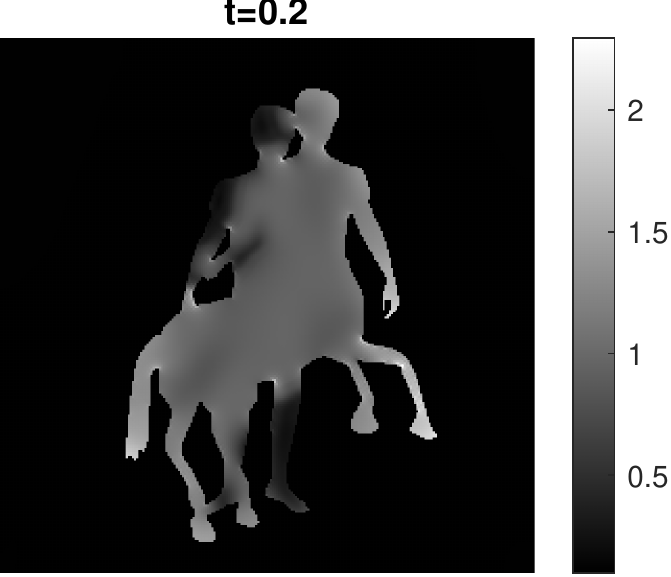}
\includegraphics[width=2.4cm]{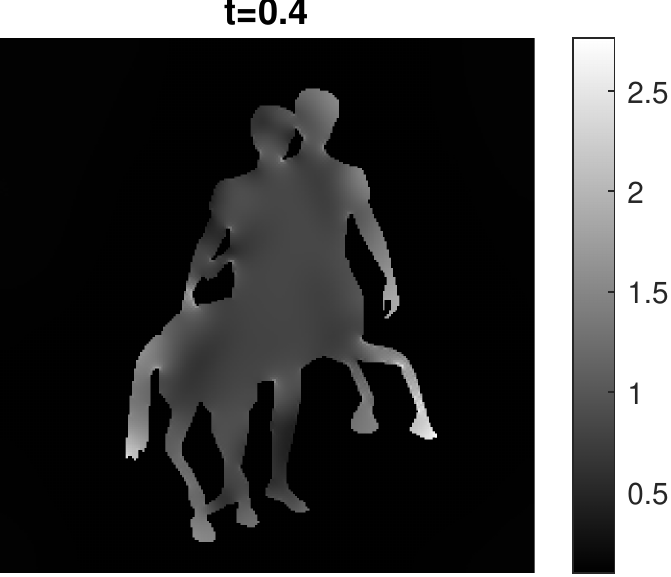}
\includegraphics[width=2.4cm]{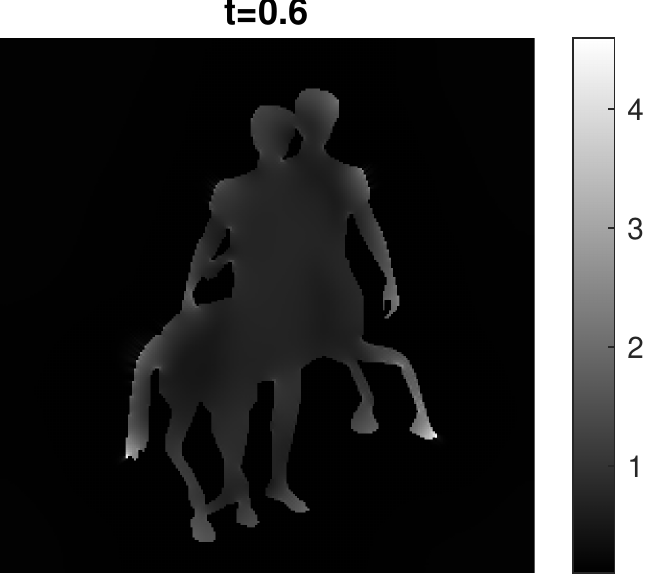}
\includegraphics[width=2.4cm]{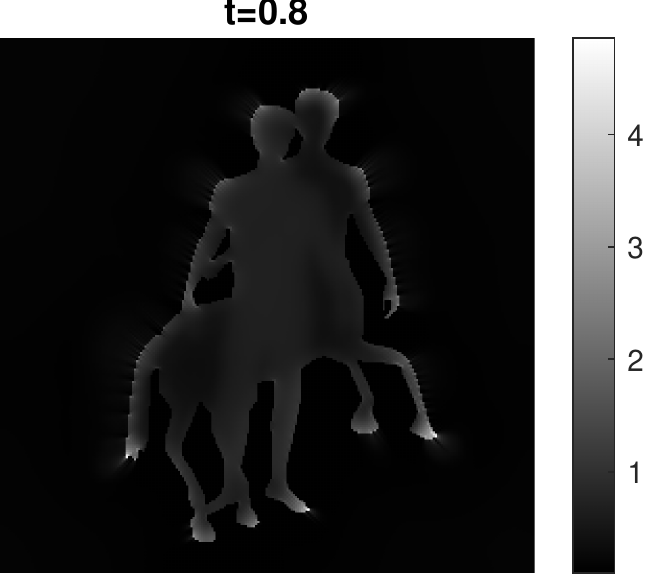}
\includegraphics[width=2.4cm]{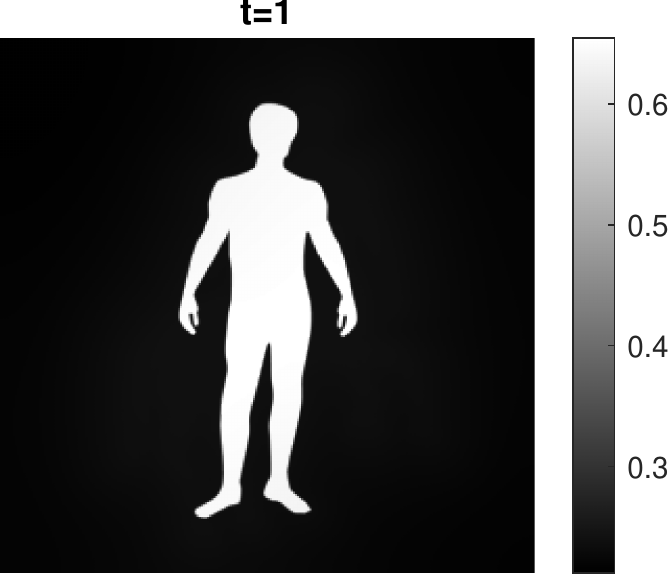}
}
\caption{Snapshot of $\rho$ at $t=0.2, 0.4, 0.6, 0.8, 1$ (from left to right)}
\label{fig: num image snapshot}
\end{figure}

\cref{fig: num image snapshot} show the snapshots of the density evolutions. Since \cref{eqn: num eg model2}-\cref{eqn: num eg mfg} set the space preference to the evolution, the mass evolutions are within the dark region and the optimal transport model \cref{eqn: num eg ot} has a more free evolution style.

Comparing model \cref{eqn: num eg model1},\cref{eqn: num eg model2} with \cref{eqn: num eg model3}, we observe that the mass evolution of model \cref{eqn: num eg model1},\cref{eqn: num eg model2} are dense, while that of \cref{eqn: num eg model3} experiences a congest-flatten process and tends to be sparse. This is compatible with our discussions in \cref{sec: review}.

\section{Conclusion}
\label{sec: summary}

In this paper, we propose an efficient and flexible algorithm to solve potential MFP problems based on an accelerated proximal gradient algorithm. In the optimal transport setting, we can converge faster or nearly as fast as G-prox and approach optimizer with the same accuracy. With multilevel and multigrid strategies, our algorithm can be accelerated up to 10 times without sacrificing accuracy. In broader settings of MFP and MFG, our method is more flexible than primal-dual or dual algorithms as it enjoys the flexibility to handle differentiable objective functions. Theoretically, we, for the first time based on an optimization point of view, analyze the error introduced by discretizing $\rho,\bmm$, and show that under some mild assumptions, our algorithm converges to the optimizer. In the future, we expect to extend the proposed algorithms for non-potential mean field games, which have vast applications in mathematical finance, communications, and data science. 
\bibliographystyle{plain}
\bibliography{references}
\end{document}